\input amstex
\documentstyle{amsppt}
\NoBlackBoxes 
\let\aro=@ 
\input psbox 
\def \scaledpicture #1 by #2 (#3 scaled #4) 
{\dimen0=#1 \dimen1=#2 \divide\dimen0 by 1000\multiply\dimen0 by #4 
\divide\dimen1 by 1000\multiply\dimen1 by #4 
$$\psboxto(\dimen0;\dimen1){#3.ps}$$} 
\mag1000 \topmatter 
\title On the moduli space of certain smooth codimension-one foliations of the 5-sphere by complex surfaces 
\endtitle 
\rightheadtext{moduli space of foliations by complex surfaces} 
\author 
Laurent Meersseman, Alberto Verjovsky 
\endauthor 
\date 4 october, 2007\enddate 
\address 
{Laurent Meersseman}\hfill\hfill\linebreak 
\indent{I.M.B.}\hfill\hfill\linebreak \indent{Universit\'e de 
Bourgogne}\hfill\hfill\linebreak \indent{B.P. 47870}\hfill\hfill\linebreak \indent{21078 Dijon
Cedex France}\hfill\hfill 
\endaddress 
\email laurent.meersseman\@u-bourgogne.fr \endemail 
\address 
{Alberto Verjovsky}\hfill\hfill\linebreak \indent{Instituto de 
Matem\'aticas de la Universidad Nacional Auton\'oma de M\'exico}\hfill\hfill\linebreak\indent{Unidad Cuernavaca}\hfill\hfill\linebreak 
\indent{Apartado Postal 273-3, Admon. de correos 
No.3}\hfill\hfill\linebreak \indent{Cuernavaca, Morelos, 
M\'exico}\hfill\hfill 
\endaddress 
\email alberto\@matcuer.unam.mx 
\endemail 
 
\thanks 
This work was partially supported by CONACyT, grant U1 55084, PAPIIT 
(UNAM), grant IN102108 from Mexico and by joint project {\bf 
CNRS-CONACyT 16895}. 
\endthanks 

 \dedicatory
dedicated to Santiago L\'opez de Medrano on his $65^{\text{th}}$ birthday
\enddedicatory

\abstract In this paper we first determine the set of all possible 
integrable almost CR-struct\-ures on the smooth foliation of $\Bbb S^5$ 
constructed in \cite{M-V}. We give a specific concrete model of each 
of these structures. We show that this set can be naturally 
identified with $\Bbb C\times\Bbb C\times \Bbb C$. We then adapt the 
classical notions of coarse and fine moduli space to the case of a 
foliation by complex manifolds. We prove that the previous set, 
identified with $\Bbb C^3$, defines a coarse moduli space for the 
foliation of \cite{M-V}, but that it does not have a fine moduli 
space. Finally, using the same ideas we prove that the standard 
Lawson foliation on the $5$-sphere can be endowed with almost CR-structures 
but none of these is integrable. This is a foliated analogue to the 
examples of almost complex manifolds without complex structure. 
\endabstract 
 
\endtopmatter 
 
\def\C{{\Bbb C}} 
\def\R{{\Bbb R}} 
\def\Z{{\Bbb Z}}

\document 
\head {\bf 0. Introduction} 
\endhead 
 
In \cite{M-V} the authors have constructed a smooth foliation of the 
sphere ${\Bbb S}^5$ by complex surfaces i.e., an integrable and 
Levi-flat codimension-one almost CR-structure. The underlying smooth 
foliation is a variation of that given by Blaine Lawson in 
\cite{La}, however ours  is topologically different. In particular, 
it contains exactly two compact leaves. 
 
The notion of foliation of a smooth manifold by complex manifolds 
can be seen as a generalization of the notion of complex structure 
on a smooth manifold, which appears as the case of codimension zero. 
The case of codimension-one is of particular interest: the smooth 
manifold has then to be of odd-dimension, so that it is really the 
analogue in the odd-dimensional case of the existence of a complex 
structure. 
 
Keeping this in mind, it is clearly interesting to ask, in the spirit of Kodaira-Spencer and Kuranishi, for a deformation theory of 
such foliations. This is what we intend to do in this paper for the 
example constructed in \cite{M-V}. To be more precise, the purpose 
of this paper is threefold: 
 
\newpage 
 
\roster 
\item To determine the set of all possible integrable almost CR-structures on this foliation. 
This set turns out to be completely determined by the set of complex 
structures of the compact leaves. These leaves are primary Kodaira 
surfaces fibering over the same elliptic curve and we infer from 
this description that our set can be identified to 
$\C\times\C\times\C$. In particular, it is {\it finite-dimensional}. 
We also show that some of these structures admit non-trivial 
CR-automorphisms. 
 
\item To adapt the classical notions of coarse and fine moduli space to the case of 
a smooth foliation by complex manifolds. We prove that the previous 
set can be seen as a coarse moduli space for the foliation of 
\cite{M-V}, but that, due to the existence of non-trivial 
CR-automorphisms, this foliation does not have a fine moduli space. 
 
\item To prove that the standard Lawson foliation 
cannot be endowed with a CR-structure, but only with 
non-integrable almost CR-structures. 
\endroster 
\medskip 
We emphasize that the set we consider is the set of integrable 
CR-structures on a {\it fixed} smooth foliation, modulo foliated 
CR-isomorphisms (see the precise Definition given in Section 2). In 
the same way, the notions of coarse and fine moduli space are 
defined on a {\it fixed} smooth foliation. This induces some 
subtleties. In particular, strictly speaking, there is not a unique 
standard Lawson foliation, but infinitely many. All are {\it 
topologically} isomorphic, but have non-conjugated contracting 
holonomies. The same is true for ``the'' foliation of \cite{M-V}. 
This will be made precise in Section 1.

The first result shows that the CR-structure of our foliation is 
very rigid. In fact it is easy to construct examples of compact 
manifolds carrying a Levi-flat and integrable almost CR-structure whose set 
of integrable almost CR-structures is infinite-dimensional (see Section 2). 
On the other hand, the third result gives an example of a foliation 
which admits non-integrable almost CR-structures and all of whose leaves 
admit separately a complex structure but which cannot be endowed 
with an integrable almost CR-structure. This is the foliated analogue to 
the examples of almost complex surfaces without complex structure 
(see \cite{B-H-P-V, IV.9}). Moreover it proves that the search for 
the existence of integrable and Levi-flat codimension-one 
almost CR-structures on a compact manifold cannot be reduced to local 
analytic questions such as solving a $\bar\partial$-problem along 
the leaves. 
 
\medskip 
{\bf Remark.} We would like to emphasize that, as far as we 
know, the foliation described in \cite{M-V} (as well as the related examples of \cite{M-V}, Section 5) is the only known 
example of a smooth foliation by complex manifolds of complex 
dimension {\it strictly greater than one} on a {\it compact} manifold, which is not 
obtained by classical methods such as the one given by the orbits of 
a locally trivial smooth action of a complex Lie group, the natural 
product foliation on $M\times N$ where $M$ is foliated by Riemann 
surfaces and $N$ is a complex manifold, holomorphic fibrations, or 
trivial modifications of these examples such as cartesian products 
of known examples or pull-backs. Of course, it is very easy to give examples of 
foliations by complex manifolds on open manifolds (in fact even with 
Stein leaves). On the other hand, if a compact smooth manifold has an orientable smooth foliation by 
surfaces then, using a Riemannian metric and the existence of 
isothermal coordinates, we see that the foliation can be considered 
as a foliation by Riemann surfaces. 
\medskip
\newpage

Although this paper obviously depends on the example given in 
\cite{M-V} it can be read independently of that paper since it uses 
different methods and ideas. In this sense this paper is not just an 
addendum of \cite{M-V}.

\medskip
 
We would like to thank Laurent Bonavero, Michel Brion, Marco Brunella, Dominique Cerveau, Andr\'e Haefliger, Janos Kollar, Lucy Moser and Marcel Nicolau for their very valuable 
comments.
We also would like to thank the referee for indicating us with care every point needing clarification
in the first version of the article. 

The first named author would like to thank the IMATE, Cuernavaca, 
Mexico, and the Universitat Autonoma de Barcelona for their hospitality. The second named author would like 
to thank the IRMAR, Universit\'e de Rennes I, and the IMB, Universit\'e de Bourgogne for their hospitality.

\head {\bf 1. Preliminaries} 
\endhead 
 
 In this Section, we fix some notations and recall some facts and constructions around foliations by complex manifolds. Most of this material comes from \cite{M-V}.
 
 \medskip
 
In this article, {\it smooth} means $C^{\infty}$-differentiable. We 
make use of the following notations. 
\medskip 
\noindent (i) we denote by $\Bbb D$ the unit open disk of $\Bbb C$ and by $\overline{\Bbb D}$ its closure. 
 
\noindent (ii) if $X$ is a {\it complex object} (e.g. complex 
manifold, foliation with complex leaves, etc), then $X^{diff}$ 
denotes the underlying compatible {\it smooth object}. 
 
\noindent (iii) for $T_1,\hdots, T_n$ a set of automorphisms of a 
complex (respectively smooth) manifold $X$, we denote by $\langle 
T_1,\hdots, T_n\rangle$ the group generated by $T_1,\hdots, T_n$, 
and by $X/\langle T_1,\hdots, T_n\rangle$ the quotient space of $X$ 
by this group. We use this last notation only in the case where this 
quotient space is a complex (respectively smooth) manifold, i.e. 
when the action is free and totally discontinuous. 
 
\noindent (iv) for a smooth manifold $M$, we denote by $\partial M$ 
the boundary of $M$ and by Int $M$ the interior part of $M$, namely 
the open manifold $M\setminus\partial M$. 
\medskip 
We recall that, if $p:\tilde M \to M$ is a covering map of a smooth 
manifold $\tilde M$ {\it with boundary} onto a smooth manifold $M$ 
with boundary, then $p$ restricted to each component of $\partial 
\tilde M$ is a covering map onto its image, which is a component of 
$\partial M$.

\subhead 1.1. Foliations by complex manifolds 
\endsubhead 
 
Recall that an {\it almost CR-structure} on a smooth manifold $V$ is the 
data of a subbundle $E$ of the tangent bundle $TV$ together with an 
operator $J : E\to E$ acting linearly on the fibers of $E$ and 
satisfying $J^2\equiv -Id$ on every fiber. We denote an almost CR manifold 
by $(V,J)$. In this paper every almost CR-structure is assumed to be 
smooth. Recall also that a {\it CR-map} between $(V,J)$ and 
$(V',J')$ is a smooth map $f$ from $V$ to $V'$ whose differential 
commutes with the almost CR-structures, that is $df\circ J\equiv J'\circ 
df$. 
 
\remark{Remark}
In our paper \cite{M-V}, an almost CR-structure is just called a CR-structure. We would like to thank Claude LeBrun, who pointed out to us that the terminology of almost CR-structure is more appropriate 
(since it is parallel to the
notions of almost complex and complex structures).
\endremark
\medskip
 
We refer to our paper \cite{M-V} and to \cite{Tu} for the notions of 
integrability and Levi-flatness of an almost CR-structure. We just recall 
here the two definitions of \cite{M-V} which are essential for this 
article. 
 
\proclaim{Definition} Let $M$ be a smooth manifold of odd dimension 
and without boundary. A foliation by complex manifolds on $M$ is the 
data of a smooth codimension-one foliation $\Cal F$ of $M$ which is 
endowed with an integrable almost CR-structure whose corresponding 
distribution is the distribution tangent to the leaves of $\Cal F$. 
\endproclaim 
 
\remark{Remark} An integrable almost CR-structure is also called a CR-structure.
\endremark
\medskip

\remark{Remark} Such an almost CR-structure is automatically Levi-flat by 
the Frobenius Theorem. 
\endremark 
\medskip 
Equivalently, a foliation by complex manifolds $\Cal F$ of 
codimension-one on  a smooth manifold $M$ of dimension $2n+1$, can 
be defined by a foliated atlas 
$$ 
\Cal A=\{(U_i,\phi_i)_{i\in I}\quad\vert\quad \phi_i (U_i) \subset 
\Bbb C^n\times \Bbb R\cong \Bbb R^{2n+1}\} 
$$ 
such that the changes of charts 
$$ 
(z,t)\in\phi_i(U_i\cap U_j)\longmapsto \phi_j 
\circ\phi_i^{-1}(z,t):=(\xi_{ij}(z,t),\zeta_{ij}(t))\in\phi_j(U_i\cap 
U_j) 
$$ 
are holomorphic in the tangential direction, {\it i.e.} the map 
$\xi_{ij}$ is holomorphic for fixed $t$. 
\medskip 
 
Let $M$ be a smooth manifold {\it with} boundary. Let $N$ be the 
open manifold obtained by adding the collar $\partial M\times [0,1)$ 
to $M$ equipped with the unique differentiable structure such that 
the natural inclusions of $\partial M\times [0,1)$ and $M$ into $N$ 
are smooth embeddings (see [Hi, Chapter 8, Theorem 2.1]). 
 
\proclaim{Definition} A (codimension-one) {\bf tame} almost CR-structure on 
$M$ is the data of an almost CR-structure on the interior of $M$ and of an 
almost complex structure on the boundary $\partial M$ such that the 
following {\bf gluing condition} is verified. Let $M\subset{N}$ and 
$\partial M\times [0,1)\subset{N}$ be the natural embeddings. Then, 
the almost CR-structure on $M$ extends to a almost CR-structure on $N$ by 
considering on the collar the distribution tangent to the 
submanifolds $\partial M\times \{t\}$, $0\leq t<1$, and equipping 
this distribution with the natural almost complex structure 
inherited from $\partial M$. 
 
A foliation by complex manifolds on $M$ is the data of a smooth 
foliation $\Cal F$ of $M$ of codimension-one which is endowed with 
an integrable tame almost CR-structure. 
\endproclaim 
 
The same remark after the first definition is valid in this case. 
\subhead 1.2. Lawson Foliation 
\endsubhead 
 
Let us now recall Lawson's construction of a smooth codimension-one 
foliation of $\Bbb S^5$. Let 
$$ 
P\ :\ (z_1,z_2,z_3)\in\Bbb C^3\longmapsto z_1^3+z_2^3+z_3^3\in\Bbb C 
$$ 
and let $V=P^{-1}(0)$. The manifold 
$$ 
W:=V\setminus\{(0,0,0)\}=\left\{(z_1,z_2,z_3)\neq0\quad\vert\quad 
z_1^3+z_2^3+z_3^3=0 \right\} 
$$ 
intersects transversally the Euclidean unit sphere in the smooth 
compact manifold $K$. Moreover, it projects onto the projective 
space $\Bbb P^2$ as an elliptic curve $\Bbb E_{\omega}$ of modulus 
$\omega$. This curve admits an automorphism of order three, hence 
$\omega^3=1$. The canonical projection 
$$ 
p:W\to{\Bbb E}_\omega 
$$ 
describes $W$ as a holomorphic principal ${\Bbb C}^*$-bundle over 
the elliptic curve ${\Bbb E}_\omega$, with first Chern class equal 
to $-3$. By passing to the unit bundle, one has that $K$ is a 
principal circle-bundle over a torus with Euler class equal to $-3$ 
(see \cite {Mi2, Lemma 7.1 and Lemma 7.2}). 
 
Another way of seeing $K$ is the following.  Let 
$$ 
A=\pmatrix 1 &-3\\ 
0 &1 
\endpmatrix 
$$ 
then $K$ is the suspension of the unipotent isomorphism induced on 
the two-di\-men\-sio\-nal torus by the matrix $A$. Notice that, in 
this way, we define $K$ as a torus bundle over $\Bbb S^1$, but {\it 
not as a principal torus bundle}. This subtlety will be important in 
the sequel.

The Milnor Fibration Theorem \cite{Mi1, Theorem 4.8} shows that 
$\Bbb S^5\setminus{K}$ fibers over $\Bbb S^1$ and describes $\Bbb 
S^5$ as an open book in the sense of Winkelnkemper \cite{Wi}, with 
fibers diffeomorphic to $P^{-1}(z)$, $z\in\Bbb C^*$. On the other 
hand, $K$ is smoothly embedded in $\Bbb S^5$ with trivial normal 
bundle and therefore has a closed tubular neighborhood diffeomorphic 
to $K\times \overline {\Bbb D}$. Thus it follows that both a closed tubular 
neighborhood $\Cal N$ of $K$ in $\Bbb S^5$ and $\Cal E$, the closure 
of its complement, fiber over the circle by fibrations which are 
also fibrations when restricted to the boundary. Hence they are 
smoothly foliated by using a standard {\it Tourbillonnement} Lemma: 
 
\proclaim{Lemma [La]} Let $M$ be a compact manifold with boundary 
$\partial M$, and suppose there exists a $C^{\infty}$-submersion 
$\psi\ :\ M\to \Bbb S^1$ such that $\psi\vert_{\partial M}$ is a 
submersion of the boundary. Then there exists a codimension-one 
$C^{\infty}$-foliation of $M$. 
\endproclaim 
 
Gluing carefully those pieces together by a diffeomorphism, Lawson 
obtained a smooth codimension-one foliation of $\Bbb S^5$ with a 
unique compact leaf diffeomorphic to $K\times \Bbb S^1$. 
 
\remark{Remark} Notice that there is not a unique Lawson foliation, 
but infinitely many. Indeed, there are infinitely many ways of 
turbulizing, which give rise to non-conjugated holonomies of the 
compact leaf. Therefore, the standard Lawson foliation is unique 
{\it topologically} but certainly not {\it differentiably}. In the 
sequel, we will however still talk of {\it the} standard Lawson 
foliation. This is due to the fact that the result we prove (Theorem 
A) is independent of the choice of the turbulization. 
\endremark 
 
 \subhead 1.3. A foliation of $\Bbb S^{5}$ by complex surfaces
 \endsubhead
 
Let us start with the following Lemma, whose proof is immediate. A 
CR-submersion is a smooth submersion between almost CR-manifolds which is 
also a CR-map. 
 
\proclaim{Lemma} Let $(M,\Cal F)$ be a Levi-flat CR-structure, with or without boundary. Let $p : \tilde M\to M$ be a 
smooth covering. Let $\tilde {\Cal F}$ be the foliation on $\tilde 
M$ whose leaves are the connected components of the pull-back by $p$ 
of the leaves of $\Cal F$. Then, there exists a unique smooth, 
integrable almost CR-structure on $(\tilde M, \tilde{\Cal F})$ such that 
$p$ is a CR-submersion. 
\endproclaim 
 
We call $\tilde {\Cal F}$ the {\it induced pull-back foliation}. 
 
\proclaim{Definition} Let $(M,\Cal F)$ be a Levi-flat CR-structure. A covering map $\pi:\tilde{M}\to{M}$ is called a product 
foliated covering if one of the two following statements is 
satisfied: 
\medskip 
\noindent (i) if $\partial M$ is empty, then $\tilde M$ is 
diffeomorphic to $N\times \R$ (for $N$ a smooth manifold) and the 
induced pull-back foliation $\tilde \Cal F$ has leaves diffeomorphic 
to $N\times\left\{\text{\rom{point}}\right\}$. 
 
\noindent (ii) if $\partial M$ is non-empty, then $\tilde M$ is 
diffeomorphic to $(N\times\R^+)\setminus (A\times\{0\})$ (for $N$ a 
smooth manifold and $A$ a - possibly empty - analytic subset of $N$) and the foliation 
has leaves diffeomorphic to 
$N\times\left\{\text{\rom{point}}\right\}$ in the interior and 
diffeomorphic to $(N\setminus A)\times\{0\}$ on the boundary. 
\endproclaim 
 
 
\remark{Remark} In general, the complex structure on the leaves 
$N\times\left\{\text{point}\right\}$ {\it depends on} $t$ so a product foliated covering {\it is not} CR-isomorphic to a product. 
\endremark 

\remark{Remark} An example of (ii) in the definition above is the 
solid cylinder $\overline{\Bbb D}\times \Bbb R\cong(\Bbb  C\times\Bbb 
R^+)\setminus\{(0,0)\}$ corresponding to the infinite cyclic 
covering of the solid torus $\overline{\Bbb D}\times \Bbb S^1$ foliated by the 
Reeb foliation, see \cite{M-V, Lemma 2}. 
\endremark

\medskip


The detailed construction of the foliation of $\Bbb S^5$ by complex 
surfaces can be found in \cite{M-V}. We denote this foliation by 
$\Cal F_{\Bbb C}$. The idea is to endow suitable coverings of $\Cal 
N$ and of the closure of its open complement in $\Bbb S^5$ with 
trivial foliations by complex surfaces such that the covering 
transformations are CR-isomorphims, i.e. to construct product 
foliated coverings as in the definition above (here we use the 
notation of 1.2). Then, taking the quotient, we obtain foliations by 
complex surfaces of $\Cal N$ and of the closure of its open 
complement in $\Bbb S^5$. Due to the tame condition, these two 
foliations can be glued together \cite{M-V, Lemma 1}. 
 
We content ourselves with describing the two 
coverings and referring to \cite{M-V} for more details. 
 
Let $\lambda$ be a real such that $0<\lambda<1$. 

\remark{Remark}
This is not the same convention as in \cite{M-V} where $\lambda$ is supposed to be strictly greater than one.
\endremark
\medskip

The following 
function 
$$ 
z\in\Bbb C^3\longmapsto \lambda\omega\cdot z\in\C^3 
$$ 
leaves $W$ invariant. The group generated by this transformation 
acts properly and discontinuously on $W$ and the quotient is a 
compact complex manifold diffeomorphic to $K\times \Bbb S^1\simeq 
\partial\Cal N$. Let us call it $\Cal S_{\lambda}$. We remark that 
it is a primary Kodaira surface. 
 
 Let 
$$ 
\widetilde X=\Bbb C^*\times (\Bbb C\times [0,\infty)\setminus 
\{(0,0)\}) 
$$ 
and let $\Gamma$ be the group generated by the commuting 
diffeomorphisms $T$ and $S$ defined as follows: 
$$ 
\eqalign {\forall (z,u,t)\in\widetilde X,\qquad &T(z,u,t)= 
(z,\lambda\omega\cdot u,d(t))\hfill\cr \text{and}\quad &S(z,u,t)= 
(\exp(2i\pi\omega){\cdot}z, (\psi(z))^{-3}{\cdot}u,t)} 
$$ 
where $d$ is a smooth diffeomorphism equal to $t$ when $t\leq 0$ and 
satisfying $d'(t)<1$ when $t>0$ and where $\psi$ is the automorphic 
factor of $W$ as $\C^*$-bundle over $\Bbb E_{\omega}$.

\remark{Remark}
This is not the same convention as in \cite{M-V}. Indeed, it is linked to the previous remark,
for taking $\lambda>1$ (respectively $0<\lambda<1$) implies taking $d'(t)>1$ for $t>0$ (respectively $d'(t)<1$), otherwise the previous
action is not proper.
\endremark
\medskip

Let $\Cal F_i$ be the foliation whose leaves $L_t$ are the level 
sets in $\widetilde X$ of the projection on the third factor. These 
leaves are naturally complex manifolds biholomorphic to $\Bbb 
C^*\times\Bbb C$ for $t>0$ and $\Bbb C^*\times\Bbb C^*$ for $t=0$. 
The group $\Gamma$ preserves the foliation $\Cal F$ and sends one 
leaf biholomorphically onto its image. Then it is proved in 
\cite{M-V} that the quotient of $\widetilde X$ by $\Gamma$ is 
diffeomorphic to $\Cal N$ and thus provides a 
foliation by complex manifolds on this closed set. The boundary leaf 
is biholomorphic to $\Cal S_{\lambda}$. The other leaves are all 
biholomorphic to the line bundle over $\Bbb E_{\omega}$ obtained 
from $W$ by adding a zero section. 
 
On the other hand, let $g:\C^3\times [-1,\infty ) \to \C$ be the 
function defined by 
$$ 
g((z_1,z_2,z_3),t)={z_1}^3+{z_2}^3+{z_3}^3-\phi(t) 
$$ 
where $\phi$ is a smooth function which is zero exactly on the non-positive real 
numbers. 
 
Let $\widehat{\Xi}=g^{-1}(\{0\})$ and $\Xi=\widehat{\Xi}\setminus 
(\{(0,0,0)\}\times [-1,\infty))$. Then, $\Xi$ is a smooth manifold and has 
a natural smooth foliation ${\Cal{F}}_e$ by complex manifolds whose 
leaves ${\{L_t\}}_{t\in\R}$ are parametrized by projection onto the 
factor $[-1,\infty)$. 
 
Let  $G:\Xi\to{\Xi}$ be the diffeomorphism given by 
$$ 
G((z_1,z_2,z_3),t)= 
((\lambda\omega\cdot{z_1},\lambda\omega\cdot{z_2}, 
\lambda\omega\cdot{z_3}),h_{\lambda}(t)) 
$$ 
where $h_{\lambda}$ is a smooth diffeomorphism whose fixed points 
are $0$ and $(-\infty, -1]$. For good choices of $\phi$ and 
$h_{\lambda}$ which are specified in \cite{M-V}, the pair $(\Xi, 
{\Cal{F}}_e)$ is a covering of the closure of $\Bbb 
S^5\setminus \Cal N$ with deck transformation group $\Gamma$, 
generated by $G$. The boundary is a leaf biholomorphic to $\Cal 
S_\lambda$ and the gluing condition is verified. There is another 
compact leaf corresponding to $t=0$. It is also biholomorphic to 
$\Cal S_{\lambda}$. These two compact leaves form the boundary of a 
collar whose interior leaves are all biholomorphic to $W$. Finally, 
the other leaves are all biholomorphic to the affine cubic surface 
$P^{-1}(1)$ of $\C^3$. 
 
\remark{Remark} Observe that the second covering {\it is not} a product foliated covering since it has two 
topologically distinct leaves: $W$ and $P^{-1}(\{t\})$. But it is a union of product foliated coverings. Indeed, the restriction
of $\Xi$ to $[-1,0)$ is a product foliated covering, as well as its restriction to $(-1,0]$ and to $(0,\infty)$.
\endremark

\remark{Remark} The construction recalled above depends on the 
choices of the smooth functions $d$, $\phi$ and $h_{\lambda}$. 
Notice that $d$ and $h_{\lambda}$ define the holonomy of the compact 
leaves. As a consequence, if we construct such a foliation $\Cal F$ 
from $d$ and $h_{\lambda}$ and another one, say $\Cal F'$, from $d'$ 
and $h'_{\lambda}$ with the property that $d'$ (or respectively 
$h'_{\lambda}$) is not smoothly conjugated to $d$ (or respectively 
to $h_{\lambda}$), then $\Cal F$ and $\Cal F'$ are not smoothly 
isomorphic, although they are topologically isomorphic (such maps 
exist, see \cite{Se}). 
 
Nevertheless, the smooth type of the foliation is independent of the 
choice of the parameter $\lambda$ in the following sense. Fix some 
$\lambda$ and some smooth functions $d$ and $h_{\lambda}$ and call 
$\Cal F_{\lambda}$ the resulting foliation. Choose now $0<\mu<1$ 
different from $\lambda$. The function $h_{\lambda}$ has the 
property that it coincides with the parabolic M\"obius 
transformation $t/(1-3(\log \lambda)t)$ near $0$ (see \cite{M-V, p.
925}). There exists a smooth function $f : \R \to \R$ fixing $0$ 
with the property that $f\circ h_{\lambda}\circ f^{-1}$ coincides 
with the parabolic M\"obius transformation $t/(1-3(\log \mu)t)$. It 
is easy to check that this new diffeomorphism can be used as 
$h_{\mu}$. As $h_{\mu}$ is globally conjugated to $h_{\lambda}$, the 
foliation $\Cal F_{\mu}$ obtained from the previous construction 
using the functions $d$ and $h_{\mu}$ is smoothly isomorphic to 
$\Cal F_{\lambda}$. However, they are not isomorphic as 
CR-structures, since different choices of $\lambda$ yield different 
complex structures on the compact leaves (see Section 5). 
 
In the sequel, we still talk of {\it the} foliation of \cite{M-V}, 
since the results we prove (Theorems B, C, D) are independent of the 
particular choices of the functions $d$ and $h_{\lambda}$. It will 
be important however to keep in mind the independence of the 
foliation with respect to $\lambda$, as was indicated above. 
\endremark 
 
\head {\bf 2. Families of complex structures and locally trivial 
CR-fiber bundles} 
\endhead 
 
Let $(V,J)$ be a smooth almost CR-manifold and let $X$ be a smooth 
manifold. The following Definition is a reformulation in the 
language of almost CR-structures of Definition 1.1 of Kodaira-Spencer 
\cite{K-S}. 
 
\proclaim{Definition} A smooth map $\pi : V\to X$ is a smooth family 
of complex structures or a smooth deformation family if 
\medskip 
\noindent (i) It is a smooth submersion with compact fibers. 
 
\noindent (ii) The almost CR-structure $J$ is integrable. 
 
\noindent (iii) The almost CR-structure $J$ is Levi-flat and the associated 
smooth foliation is given by the level sets of $\pi$. 
\endproclaim 
 \remark{Remark}
The manifolds $V$ and $X$ can have boundary (for example, $X$ may be the closed interval $[0,1]$). In this case, we ask $\pi$ to
be a submersion on Int $V$ and also on $\partial V$.
\endremark
\medskip

By Ehresmann's Lemma \cite{Eh}, $\pi$ is a locally trivial smooth 
fiber bundle, therefore the fibers $\pi^{-1}(\{t\})$ are all 
diffeomorphic. Moreover, they are endowed with a complex structure 
obtained by restriction of $J$. We denote by $V_t$ the complex 
manifold corresponding to $\pi^{-1}(\{t\})$. 
 
\proclaim{Definition} A smooth family of complex structures $\pi : 
V\to X$ is a locally trivial CR-fiber bundle if 
\medskip 
\noindent (i) All the fibers $V_t$ are biholomorphic to a fixed 
complex manifold $V_0$. 
 
\noindent (ii) There is an open covering $(U_{\alpha})$ of $X$ and 
CR-isomorphisms $\phi_{\alpha} : \pi^{-1}(U_{\alpha}) \to 
V_{0}\times U_{\alpha}$ (where the CR-structure of $V_{0}\times 
U_{\alpha}$ is given by the complex tangent distribution to 
$V_{0}\times \{x\}$ for $x$ varying in $U_{\alpha}$). 
 
\noindent (iii) There are commutative diagrams 
$$ 
\CD \pi^{-1}(U_{\alpha}) \aro >\phi_{\alpha}>> V_{0}\times 
U_{\alpha}\cr \aro V \pi VV \aro V p VV \cr U_{\alpha} \aro > Id>> 
U_{\alpha} 
\endCD 
$$ 
where $p$ is the natural projection. 
\endproclaim 
 
The following Proposition is the CR-version of a classical result of 
Fischer and Grauert. 
 
\proclaim{Proposition 1 (see [F-G])} A smooth family of complex 
structures $\pi : V\to X$ is a locally trivial CR-fiber bundle if 
and only if all the fibers $V_t$ are biholomorphic. 
\endproclaim 
 
The proof of this result is rigorously identical to that given in 
\cite{F-G}. It uses in an essential way Theorem 6.2 of \cite{K-S}. 
 
In the particular case where $X$ is diffeomorphic to an interval, 
then $\pi$ is globally trivial \cite{St, Theorem 11.4}, that is the following diagram is 
commutative 
 
$$ 
\CD V \aro >>> V_{0}\times X\cr \aro V \pi VV \aro V p VV \cr X \aro 
> Id>> X 
\endCD 
$$ 
 
 
The following immediate Corollary will be used frequently in the sequel. 
 
\proclaim{Proposition 2} Let $(M,\Cal F)$ be a Levi-flat CR-structure 
manifolds. Let $\pi : \tilde M\to M$ be a product foliated covering without boundary. 
Assume that all the leaves of $\tilde M$ are compact and 
are biholomorphic to a fixed compact complex manifold $N$.

Then, $\tilde M$ is 
CR-isomorphic to $N\times \Bbb R$. 
 
\endproclaim 
 
\remark{Remark}
Except for very particular cases, we do not know if such a statement is true in the case with boundary, that is: if all the interior leaves of $\tilde M$ are compact and 
are biholomorphic to a fixed compact complex manifold $N$, and if the boundary leaf is biholomorphic to $N\setminus A$ for $A$ as in the Definition of a product foliated covering, then
 $\tilde M$ is 
CR-isomorphic to $(N\times\Bbb R^+)\setminus (A\times \{0\})$. 

If it was the case, many arguments in
the sequel could be greatly simplified. The difficulty is to prove that the Levi-flat CR-structure on $\tilde M$ can be extended to $A$ on the boundary leaf: this of course can be done {\it along} the boundary leaf; but one has to check that this extension is smooth transversally to the leaves and this is not clear at all. The Example given below even shows that it could be false.
\endremark
\remark{Remark}
Of course, Proposition 2 is valid in the case of a product foliated covering with boundary if the set $A$ is empty.
\endremark 
\medskip

 In the sequel, we will also consider smooth deformation families of {\it non-compact} manifolds. The only difference in the definition is that we have to impose that the
 family is smoothly trivial, since Ehresmann's Lemma is false in the non-compact case. The notion of CR-triviality is then exactly the same as before.
 
One of the main technical difficulties in the sequel however is that statements such as Proposition 1 and 2 are false when the leaves are non-compact. Here is 
an Example showing how subtle is the situation in the non-compact case (the first author wants to thank Marco Brunella for explaining him this example).

\example{Example}
Consider the trivial foliation by Riemann spheres $M=\Bbb P^1\times [0,1]$. Let $s$ be a continuous section from $[0,1]$ to $M$ and let
$X$ be obtained from $M$ by removing the image of this section. 
Then the natural projection map $M\to [0,1]$ is diffeomorphic to the trivial
bundle $\Bbb R^2\times [0,1]\to [0,1]$. So is a deformation family of (non-compact) complex manifolds, each of them being a copy of the
complex plane. Although the leaves are all biholomorphic, and although the family $X$ can be compactified as the trivial deformation family
$M$, it is not CR-isomorphic to the product $\Bbb C\times [0,1]$. 

Assume the contrary.  Consider the CR-isomorphism
$$
X\aro >\simeq >> \Bbb C\times [0,1] 
$$
Let us denote by $i$ this map. 

If we have a look at the following diagram
$$
\CD
X \aro >i >> \Bbb C\times [0,1] \cr
\aro V \text {natural inclusion} VV \aro VV \text{natural inclusion} V\cr
M=\Bbb P^1\times [0,1] \aro >\phi >>\Bbb P^1\times [0,1]
\endCD
$$
then we see that the bottom map $\phi$, which is a priori only defined outside the section $s$, extends as a CR-isomorphism of $\Bbb P^1\times [0,1]$. Indeed, for fixed $t$, the map $\phi_t$ extends
continuously and thus holomorphically to $\Bbb P^1$ by setting $\phi_t(s(t))=\infty$. Hence $(\phi_t)$ is a family of rational maps of degree one of $\Bbb P^1$. Since this family is smooth in $t$ when restricted to $\Bbb P^1
\setminus s$, the coefficients of these rational maps are smooth and the family is smooth in $t$ on the whole $\Bbb P^1$. Now, this CR-isomorphism sends 
the {\it continuous} section $s$ to the {\it smooth} section $\infty\times [0,1]$. Contradiction.

\endexample

\head
{\bf 3. Some examples of deformation families}
\endhead

In this Section, we first recall some basic facts about families of line bundles over elliptic curves, cf \cite{Gu}, \cite{G-H}.
The only part which is not classical (although it is an easy consequence of classical facts) is the Dumping Lemma.
\medskip

Let $\alpha\in\Bbb H$ and let $n\in\Bbb Z$. The subset $\text{Pic}_{n}(\Bbb E_{\alpha})$ of the Picard group of the elliptic curve $\Bbb E_{\alpha}$ is constituted
by elements corresponding to line bundles of Chern number $n$. It has a natural structure of an elliptic curve \cite{Gu, \S 7-8}.
\medskip

This structure of an elliptic curve makes a moduli space of $\text{Pic}_{n}(\Bbb E_{\alpha})$. This means in particular the following.
\medskip

Let $\pi : (\Cal X, J)\to [0,1]$ be a smooth family of deformations of line bundles over $\Bbb E_{\alpha}$ of fixed topological degree $n$, that is (cf \cite{K-S, \S III.7})

\medskip
\noindent (i) The map $\pi$ is 
a trivial smooth bundle, every fiber of $\pi$ is a line bundle of Chern number $n$ over $\Bbb E_{\alpha}$.

\noindent (ii)  There is a commutative diagram
$$
\CD
(\Cal X, J) \aro >\pi>> [0,1] \cr
\aro Vp VV \aro VV Id V\cr
\Bbb E_\alpha\times [0,1] \aro> \text{2nd projection} >>[0,1]
\endCD
$$
\noindent where the restriction of $p$ to a fiber of $\pi$ is the bundle projection.

\noindent (iii) The map $p : (\Cal X,J)\to \Bbb E_\alpha\times [0,1]$ is a smooth fiber bundle with fiber $\Bbb C$ and structural group $\Bbb C^*$.

\example{Example}
Let $\pi : \tilde M\to M$ be a product foliated covering with surfaces as leaves. Assume that we have a smooth CR-embedding 
$$
i : \Bbb E_\alpha\times\Bbb R \to \tilde M
$$
(or $i : \Bbb E_\alpha\times [0,\infty) \to \tilde M$ in the case of a product foliated covering with boundary). We may choose locally defining functions for the submanifolds
$E_t=i_t(\Bbb E_\alpha)$ (that is local holomorphic functions on the leaf $L_t$ of $\tilde M$ whose zero set defines an open set of $E_t$) which depend smoothly on $t$.
Therefore, we may choose for the normal bundles of $E_t$ in $L_t$ cocycles depending smoothly on $t$. Hence we may construct abstractly
a family of deformations of line bundles over $\Bbb E_\alpha$ of fixed topological degree $p : (\Cal X,J)\to \Bbb E_\alpha\times [0,1]$ such that $p^{-1}(\{t\})$
is the normal bundle of $E_t$ in $L_t$.

\endexample
\medskip

 Then the natural map
$$
j\ :\ (\Cal X, J)\longrightarrow \text{Pic}_{n}(\Bbb E_{\alpha})
$$ 
sending a fiber $X_{t}=\pi^{-1}(\{t\})$ onto the element of $\text{Pic}_{n}(\Bbb E_{\alpha})$ characterizing it as a line bundle is a smooth map. Indeed, in such a situation, we have a diagram
of coverings
$$
\CD
\Bbb C\times\Bbb C\times [0,1] \aro >>> (\Cal X, J)\cr
\aro VVV \aro VV p V\cr
\Bbb C\times [0,1] \aro >>> \Bbb E_\alpha \times [0,1]
\endCD
$$
so that we may locally choose for the bundles $p^{-1}(\Bbb E_\alpha\times\{t\})\to 
\Bbb E_\alpha\times\{t\}$ a set of multipliers depending smoothly on $t$. This is enough to show the result (cf \cite{G-H, \S 2.6}).
\medskip

From this, we deduce easily the following Lemma.

\proclaim{Dumping Lemma}
Let $\pi : (\Cal X, J)\to [0,1]$ be a smooth family of deformations of line bundles over $\Bbb E_{\alpha}$ of 
fixed topological degree $n$. Assume that, for every $t_{0}\in [0,1)$, there exists a sequence $(t_{n})$ with $1$ as limit such that $X_{t_{n}}$ is isomorphic 
(as a line bundle) to $X_{t_{0}}$.

Then all the fibers $X_{t}$ of $\pi$ are isomorphic.
\endproclaim

\remark{Remark}
Let $\pi_{0} : X_{0}\to\Bbb E_{\alpha}$ and $\pi_{1} : X_{1}\to\Bbb E_{\alpha}$ be two line bundles of fixed topological degree. Assume that they are
biholomorphic {\it as complex manifolds}. Let $f$ be a biholomorphism between them. Notice that $f$ extends as a biholomorphism between the total spaces of the
associated $\Bbb P^1$-bundles, say $X_0^c$ and $X_1^c$. Now, $f$ must preserve the Albanese varieties of $X_0^c$ and $X_1^c$, that is we have a commutative diagram
$$
\CD X_0^c \aro > f>> X_1^c \cr
\aro VVV \aro VVV\cr
\Bbb E_\alpha \aro >>> \Bbb E_{\alpha}
\endCD
$$
so $f$ maps biholomorphically each fiber of $\pi_0$ onto a fiber of $\pi_1$. It is then straightforward to check that $f$
must be a linear automorphism when restricted to
a fiber if the bundles are holomorphically non-trivial. Hence, in the non-trivial case, the word isomorphic in the previous Lemma could be replaced by biholomorphic.
\endremark

\demo{Proof}
We consider the previously described smooth map:
$$
j\ :\ (\Cal X, J)\longrightarrow \text{Pic}_{n}(\Bbb E_{\alpha})
$$ 
where $n$ is the common Chern number of the fibers of $\Cal X$. Let $t_{0}\in [0,1)$. By assumption, there exists $(t_{n})$ with $1$ as limit such that
$j(X_{t_{n}})=j(X_{t_{0}})$. Hence $j(X_{t_{0}})=j(X_{1})$ by continuity of $j$, that is the map $j$ is constant and all the fibers are isomorphic.
$\square$
\enddemo

\medskip

We finish this part with a short study of CR-suspensions, which give examples of deformation families of both compact and non-compact manifolds. 
We are mainly interested in knowing when two such families are CR-isomorphic.
\medskip

(CR)-suspensions form a simple, but important case of product foliated coverings (with empty boundary). That is, given $L$ a complex manifold and $A$ a biholomorphism
 of $L$, form the smooth manifold $X_A=(L\times \Bbb R)/\sim$ where the equivalence relation $\sim$ is given by
 $$
 (z,t)\sim (w,s) \iff w=A^{p}(z),\quad s=t+p \quad\text{for some }p\in\Bbb Z\ .
 $$
Consider the trivial foliation of $L\times\Bbb R$ by complex leaves $L\times\{pt\}$. It is preserved under the equivalence relation and descends to a foliation
by complex manifolds of $X_A$. Each leaf is biholomorphic to $L$ and the natural projection map
$$
L\times\Bbb R \longrightarrow X
$$
is a foliated product covering. 

Another way of describing the suspension is the following. The foliated manifold $X_A$ is obtained from $L\times [0,1]$ endowed with its trivial foliation 
by gluing $L\times\{0\}$ and $L\times \{1\}$ by $A$.  
If $(A_t)$ is a smooth isotopy of biholomorphisms of $L$ between $A=A_0$ and $A_1$, then the CR-map
$$
(z,t)\in L\times [0,1] \longmapsto (A_t\circ A_0^{-1}(z), t)\in L\times [0,1]
$$
descends as a CR-isomorphism between $X_{A_0}$ and $X_{A_1}$.
\medskip

Conversely, we have the following result.

\proclaim{Proposition 3}
Let $\pi : \tilde M\to M$ be a product foliated covering with $\partial M$ empty. Assume that $\tilde M$ is CR-isomorphic to $L\times\Bbb R$ for some complex
manifold $L$. Moreover, assume that the deck transformation group is isomorphic to $\Bbb Z$ and that it acts without fixing any leaf $L\times\{pt\}$.

Then, there exists a well-defined biholomorphism $A$ of $L$ (up to smooth isotopy) and $\tilde M\to M$ is CR-isomorphic as a covering space to 
$L\times\Bbb R\to X_A$ (where $X_A$ is the suspension of $L$ by $A$), i.e. the following diagram is commutative:
$$
\CD
\tilde M \aro >\text{CR-isomorphism}>> L\times\Bbb R\cr
\aro V \pi VV \aro VVV\cr
M \aro>\text{CR-isomorphism}>> X
\endCD
$$
\endproclaim

We call this biholomorphism the {\it monodromy} of the product foliated covering.

\demo{Proof}
From the hypotheses, for any $t$, there exists $s>t$ such that a fundamental domain for the action is $L\times [t,s]$, once $L\times\{t\}$
is identified with $L\times\{s\}$ by $A_t$ (where $A(z,t)=(A_t(z), f(t))$ is a well-chosen generator of the deck transformation group). 
Hence $(M, \Cal F)$ is obtained as the suspension of $L$ by any $A_t$ for fixed $t$.
$\square$
\enddemo

It should be noticed that there exist foliated product coverings which are not {\it smoothly} isomorphic to a trivial one. Here is an example.

\example{Example}
Let
$$
A=\pmatrix
2 &1 \\
1 &1
\endpmatrix
$$
 and let $M$ be the manifold obtained as the suspension of a real 2-torus $T$ by $A$. The level sets of the suspension map $M\to\Bbb S^{1}$ foliates $M$ by
 copies of $T$. 
 
 We claim that this foliation $\Cal F$ can be turned into a foliation by elliptic curves. Indeed, choose any smooth path $c$ in the upper half-plane $\Bbb H$ such that
 $c^{(n)}(1)=A\cdot c{(n)}(0)$, i.e. 
 $$
 c^{(n)}(1)=\dfrac{2+c^{(n)}(0)}{1+c^{(n)}(0)}
 $$ 
 for all $n\in\Bbb N$.
 
 Consider then the action of $\Bbb Z^{3}$ onto $\Bbb C\times\Bbb R$ given by
 $$
 ((p,q,r),(z,t))\in\Bbb Z^{3}\times (\Bbb C\times\Bbb R)\longmapsto \left (\left (\dfrac{1}{1+c(0)}\right )^{r}\cdot (z+p+qc(t)),t+r\right )\in\Bbb C\times\Bbb R\ .
  $$
 where we extend $c$ to $\Bbb R$ by setting $c(t)=A^{E(t)}\cdot c(t-E(t))$ where $E(t)$ is the integer part of $t$.
 
 The quotient space of $\Bbb C\times\Bbb R$ by this action is exactly $M$ and the trivial foliation by copies of $\Bbb C$ descends to $M$ and turns $\Cal F$ into a foliation
 by elliptic curves as wanted. Identifying $\Bbb S^{1}$ and $\Bbb R/\Bbb Z$ and letting the brackets denote the class of a real in $\Bbb R/\Bbb Z$, we have that the
 leaf over $[t]$ is the elliptic curve of modulus $c(t)$.
 
 Setting $\tilde M=T\times\Bbb R$, we have a product foliated covering $\tilde M\to M$. Since in this construction, the path $c$ ends at a different point from its starting
 point, and thus cannot be constant, we observe that the leaves of $\tilde M$ are not all biholomorphic. So this product foliated covering cannot be CR-trivial.
 
 This is not due to the construction. We claim that there does not exist on $(M,\Cal F)$ a complex structure such that $\tilde M$ is CR-trivial. For if we could find such 
 a structure, then by Proposition 3, the monodromy map
 $$
A=\pmatrix
2 &1 \\
1 &1
\endpmatrix
$$
would be isotopic to an automorphism of the elliptic curve serving as complex leaf of the foliation. But no elliptic curve admits a biholomorphism smoothly isotopic to $A$.
Hence the result. Observe that, in fact, using Proposition 2, we have that $(M,\Cal F)$ cannot be endowed with a complex structure whose 
leaves are all biholomorphic. Observe also that all this implies what was announced as an introduction to this Example: the CR-structure $(M,\Cal F)$ is not {\it smoothly}
isomorphic to a trivial CR-structure.
  \endexample
  \medskip
  
  On the other hand, a uniformization result can be drawn from Proposition 3.
   
   \proclaim{Corollary 1}
   Let $\pi_{i}\ :\ (\tilde M,\tilde{\Cal F})\to (M,\Cal F_{i})$ (for $i=0,1$)
be two product foliated coverings satisfying the hypotheses of Proposition 3. Moreover, assume that:
\medskip
\noindent (i) We have $\Cal F_{0}^{diff}=\Cal F_{1}^{diff}, \ \pi_{0}^{diff}=\pi_{1}^{diff}$.



\noindent (ii) The monodromies are equal (up to smooth isotopy of biholomorphisms).
\medskip

Then $(M,\Cal F_{0})$ and $(M,\Cal F_{1})$ are CR-isomorphic.
\endproclaim

\demo{Proof}
Apply twice Proposition 3. It tells you that both $(M,{\Cal F_{0}})$ and $(M,{\Cal F_{1}})$ are CR-isomorphic to the suspension of
the same complex manifold by the same map. Hence are CR-isomorphic.
$\square$
\enddemo

Let us now have a look to the non-empty boundary case. Let $\pi : (\tilde M,\tilde{\Cal F})\to (M,\Cal F)$ be a product foliated covering with boundary. 
Assume that it is CR-trivial, i.e. that it is CR-isomorphic to $L\times [0,\infty)\setminus A\times\{0\}$ for $L$ a fixed complex 
manifold and $A$ a subset of $L$. Assume also that the deck transformation group is isomorphic to $\Bbb Z$ 
and does not fix any leaf of the interior. This is analogous to the case of suspensions. 

In the previous model, the action of the deck transformation group is generated by a map 
$$
(z,t)\in L\times\Bbb R^+\setminus A\times\{0\}\longmapsto (T(z,t), d(t))\in L\times\Bbb R^+\setminus A\times\{0\}
$$
Since the interior of $\tilde M$ is CR-trivial, it has a well-defined monodromy by Proposition 3. Notice that $d$ has no positive fixed point since the deck transformation group does not fix any leaf of the interior.


We have

\proclaim{Proposition 4}
Let $\pi\ :\ (\tilde M,\tilde{\Cal F})\to (M,\Cal F)$  be a product foliated covering with boundary satisfying the hypotheses above. Assume moreover that
the biholomorphism $T_0(-)=T(-,0)$ induced on the boundary leaf by a generator of the deck transformation group extends as a biholomorphism of $L$.


\medskip
Then there is a commutative diagram
$$
\CD
(\tilde M,\tilde{\Cal F})\aro >\text{CR-isomorphism}>> L\times\Bbb R^+\setminus A\times\{0\}\cr
\aro V \pi VV \aro VVV \cr
(M,\Cal F) \aro >\text{CR-isomorphism}>> X
\endCD
$$
where $X$ is CR-isomorphic to the quotient of $L\times\Bbb R^+\setminus A\times\{0\}$ by the group generated by
$$
(z,t)\in L\times\Bbb R^+\setminus A\times\{0\}\longmapsto (T_0(z),d(t))\in L\times\Bbb R^+\setminus A\times\{0\}
$$
\endproclaim


\demo{Proof}
Since $T_0$ extends as
a biholomorphism of $L$ and since 
the monodromy of Int $\tilde M$ is smoothly isotopic to $T_0$, we may also choose as monodromy of this suspension the map $T_0$.
Observe that the CR-isomorphism sending Int $\tilde M$ to $L\times (0,\infty)$ and conjugating the generator $(z,t)\mapsto (T_t(z),d(t))$ to $(z,t)\mapsto (T_0(z),d(t))$ may be chosen
to extend as the identity on the boundary. Indeed, this extension can be constructed as follows. Set
$$
W=(L\times\Bbb R^+\setminus A\times\{0\})\times [0,1]
$$
and
$$
S\ :\ (z,t,s)\in W\longmapsto (T_{st}(z),d(t),s)\in W\ .
$$
Then $S$ generates a free and proper CR-action on $W$ whose quotient $V$ is diffeomorphic to $M\times [0,1]$ endowed with a Levi-flat CR-structure $\Cal G$ such that
\medskip
\noindent (i) The slice $(M\times\{0\},\Cal G_{\vert M\times\{0\}})$ is CR-isomorphic to $X$.

\noindent (ii) The slice $(M\times\{1\},\Cal G_{\vert M\times\{1\}})$ is CR-isomorphic to $(M,\Cal F)$.
\medskip
Consider the smooth vector field
$$
\xi(z,t,s_0)=\dfrac{d}{ds}\vert _{s=s_0}(T_{st}(z),d(t),s)\ .
\leqno (z,t,s)\in W
$$
It respects the foliation and is holomorphic along the leaves. It is also invariant by the action and descends as a vector field $\zeta$ on $V$ with the same properties with respect to $\Cal G$ this time.
Seeing $V$ as a closed subset of $M\times\Bbb R$, we may extend $\zeta$ as a vector field in $M\times\Bbb R$ and take its flow. The time $1$ flow sends 
CR-isomorphically the slice $(M\times\{0\},\Cal G_{\vert M\times\{0\}})$ to the slice
$(M\times\{1\},\Cal G_{\vert M\times\{1\}})$.
$\square$
\enddemo

As in the case of suspensions, we draw a uniformization result.

\proclaim{Corollary 2}
   Let $\pi_{i}\ :\ (\tilde M,\tilde{\Cal F})\to (M,\Cal F_{i})$ (for $i=0,1$)
be two product foliated coverings satisfying the hypotheses of Proposition 4. Moreover, assume that:
\medskip
\noindent (i) We have $\Cal F_{0}^{diff}=\Cal F_{1}^{diff}, \ \pi_{0}^{diff}=\pi_{1}^{diff}$.


\noindent (ii) The interior leaf of $\tilde{\Cal F_{0}}$ is biholomorphic to the interior leaf of $\tilde{\Cal F_{1}}$.

\noindent (iii) The boundary leaf of $\tilde{\Cal F_{0}}$ is biholomorphic to the boundary leaf of $\tilde{\Cal F_{1}}$.

\noindent (iv) The biholomorphisms $T_0$ and $T_1$ induced on each boundary leaf by a generator of the deck transformation group are holomorphically conjugated.
\medskip

Then $(M,\Cal F_{0})$ and $(M,\Cal F_{1})$ are CR-isomorphic.
\endproclaim

\demo{Proof}
Using Proposition 4, we have that $(M,\Cal F_{0})$ (respectively $(M,\Cal F_{1})$) is CR-isomorphic to the quotient of 
$L\times\Bbb R^+\setminus A\times\{0\}$ by the group generated by $(z,t)\mapsto (T_0(z),d(t))$ (respectively
 $(z,t)\mapsto (T_1(z),d(t))$. By (iv), these maps are CR-conjugated.
$\square$
\enddemo

We add some important remarks and consequences.

\remark{Remark}
Notice that the map $d$ corresponds to the holonomy of the boundary leaf. It is thus a smooth invariant, that is $\Cal F_{0}^{diff}=\Cal F_{1}^{diff}$
implies that the corresponding holonomies are smoothly conjugated and therefore can be assumed to be the same. In the sequel, since we deal with tame foliations,
the function $d$ will be tangent to the identity at $0$. Nevertheless, it is important to notice that the previous Proposition is valid just assuming that $d$ has no other fixed
point that $0$.
\endremark

\remark{Remark}
The tame condition is not stable under CR-isomorphisms. For example, consider as above the quotient of $L\times\Bbb R^+\setminus A\times\{0\}$ by the group generated by $(z,t)\mapsto (T_t(z),d(t))$. 
If the function $d$
is tangent to the identity at $0$, but the function $T$ is not (in $t$), then the product foliated covering could not be tame. 
\medskip

However, Proposition 4 shows that it is CR-isomorphic to the quotient of
$L\times\Bbb R^+\setminus A\times\{0\}$ by the group generated by $(z,t)\mapsto (T_0(z),d(t))$, which is obviously tame.
 Indeed, this is true in greater generality and we have the following corollary.
\endremark

\proclaim{Corollary 3}
Let $\pi : (\tilde M,\tilde{\Cal F})\to (M,\Cal F)$ be a product foliated covering with boundary. Assume that the deck transformation group is generated by $\Bbb Z$ and does not fix any interior leaf. Assume moreover that the holonomy of the boundary leaf is tangent to
the identity. Then $(M,\Cal F)$ is CR-isomorphic to a tame foliation.
\endproclaim

\demo{Proof}
This is just a reparametrization argument. We may use the following model for $(\tilde M,\tilde \Cal F)$:
$$
(\tilde M, \tilde\Cal F)=(L^{diff}\times [0,\infty)\setminus A^{diff}\times\{0\}, J)
$$ 
with action generated by a CR-map:
$$
(z,t)\in L^{diff}\times\Bbb R^+\setminus A^{diff}\times\{0\}\buildrel T\over\longmapsto (T_t(z), d(t))\in L^{diff}\times\Bbb R^+\setminus A^{diff}\times\{0\}
$$
Let  $\theta$ be a diffeomorphism of $\Bbb R^+$ flat at $0$ (for example, $\theta(t)=\exp(1/t)$ for $t>0$).
Call $L_t$ the leaf corresponding to $t$ and $J_t$ the complex structure on $L_t$. Perform the diffeomorphism
$$
(z,t)\in L^{diff}\times\Bbb R^+\setminus A^{diff}\times\{0\}\buildrel\Phi\over\longmapsto (z,\theta^{-1}(t))\in L^{diff}\times\Bbb R^+\setminus A^{diff}\times\{0\}
$$
and call $J'$ the CR-structure on $\tilde M$ which makes a CR-isomorphism of this diffeomorphism. Observe that
\medskip
\noindent (i) For all $t$, we have $J'_t\equiv J_{\theta(t)}$.

\noindent (ii) The map $T$ is conjugated through $\Phi$ to the map $S(z,t)=(T_{\theta(t)}(z),\theta^{-1}\circ d\circ\theta(t))$.
\medskip
This shows that the quotient of $(\tilde M, J')$ by the action generated by $S$, which is CR-isomorphic to $(M,\Cal F)$, is tame.
$\square$
\enddemo

 
\head {\bf 4. The compactification Lemma} 
\endhead 
 
Let us start by remarking that a non-compact leaf of $\Cal F_{\Bbb 
C}$ has either one topological end (if the leaf is a smooth affine 
cubic or a line bundle over an elliptic curve) or two topological 
ends (if the leaf is a principal $\C^*$-bundle over an elliptic 
curve). Here we represent ends of a non-compact leaf as descending 
sequences $\Cal U_1\supseteq\Cal U_2\supseteq\cdots$ of open subsets 
of the leaf such that 
$\overset{.}\to{\partial}\overset{.}\to{\bar{\Cal U_i}}$ is compact 
and $\underset{i\geq1}\to\cap\overset{.}\to{\bar{\Cal 
U_i}}=\emptyset$ where $\overset{.}\to{\bar{\Cal U_i}}$ and 
$\overset{.}\to{\partial}$ are, respectively,  the closure and 
boundary of ${\Cal U_i}$  {\it in the leaf} with its topology as a 
manifold. The ends of the non-compact leaves ``accumulate'' 
respectively to one or two of the compact leaves. In other words, an 
end of a non-compact leaf spirals around a compact leaf, or, more 
precisely,  $\underset{i\geq1}\to\cap{\bar{\Cal U_i}}=L_0$, where 
$L_0$ is a compact leaf and ${\bar{\Cal U_i}}$ is the closure of 
$\Cal U_i$ in $\Bbb S^5$. This imposes certain recurrence of the 
complex structure near the ends. Thus one may expect that the 
complex structure at infinity of the non-compact leaves is fixed by 
the complex structure of the compact leaves. The precise formulation 
of this rigidity property takes the form of the compactification 
Lemma stated below. In fact, this is the central idea of this paper. We need one more definition before stating this Lemma.

\proclaim{Definition} Let $X$ be a complex manifold of dimension 
$n$. Let $H$ be a compact complex manifold of dimension $n-1$. Let 
$E$ be an end of $X$. Then we say that $X$ admits a partial 
holomorphic compactification at $E$-infinity by adding $H$ if there 
exists a structure of complex manifold on the disjoint union 
$X\sqcup H$ such that 
\medskip 
\noindent (i) The natural injections $X\to X\sqcup H$ and $H\to 
X\sqcup H$ are holomorphic. 
 
\noindent (ii) The submanifold $H$ of $X\sqcup H$ is the limit set 
of $E$. 
\endproclaim 
 
The second point means the following: representing $E$ by a sequence 
$\Cal U_1\supseteq\Cal U_2\supseteq\cdots$ as before, we have that 
$\underset{i\geq1}\to\cap{\bar{\Cal U_i}}=H$, where $\bar{\Cal U_i}$ 
denotes the closure of $\Cal U_i$ in $X\sqcup H$ equipped with the 
topology coming from its structure of a complex manifold. 
 
\remark{Remark} If $X\sqcup H$ is compact, then this implies that 
$X$ has just one end and we say that $X$ admits a holomorphic 
compactification by adding $H$. 
\endremark 
\medskip 
 
Here is an example. Let $H$ be a compact complex manifold and let 
$X$ be a principal $\C^*$-bundle over $H$. Then $X$ admits a partial 
holomorphic compactification at $0$-infinity by adding a zero 
section $H$. In this case $X\sqcup H$ is the associated line bundle 
over $H$. 
\medskip 

Let $(M,\Cal F_0)$ be a foliation by complex manifolds with $\partial 
M\not =\emptyset$. Let $L_0$ be a leaf of $\text{Int}(M)$ such that 
the limit set,  $\underset{i\geq1}\to\cap{\bar{\Cal U_i}}=C$ 
corresponding to a given end $\Cal U_1\supseteq\Cal 
U_2\supseteq\cdots$ of $L_0$ is a compact connected component $C$ of 
$\partial M$.

Now, let $(M,\Cal F_1)$ be another foliation by complex leaves of $M$ such that
\medskip
\noindent (i) The underlying smooth foliations are equal, that is $(\Cal F_0)_{diff}=(\Cal F_1)_{diff}$.

\noindent (ii) On $C$, the complex structures $J_1$ (respectively $J_0$) induced by $\Cal F_1$ (respectively $\Cal F_0$) agree,
that is $(J_0)_{\vert C}\equiv (J_1)_{\vert C}$.
\medskip
We want to compare $L_0$ to $L_1$ (the corresponding leaf of $\Cal F_1$) {\it as abstract complex manifolds}. Since they are of course
diffeomorphic, we want to compare their complex structures. A priori, there is no reason that they are biholomorphic. However, since
the common limit of $L_0$ and $L_1$ is $C$, condition (ii) above means that these complex structures are in a sense close near $C$.

We are now in position to prove the Compactification Lemma. Roughly 
speaking, it states that, if $L_0$ can be compactified 
holomorphically at $C$-infinity, then so does $L_1$, since their 
complex structures are asymptotic near $C$.

\proclaim{Compactification Lemma}
With the hypotheses above, $L_1$ admits a partial holomorphic compactification at $C$-infinity by adding a compact complex
manifold $H$ if and only if $L_0$ admits a partial holomorphic compactification at $C$-infinity by adding a compact complex
manifold $H$.
\endproclaim

\demo{Proof}
The statement is clearly symmetric, so assume that $L_0$ admits a partial holomorphic compactification at $C$-infinity by adding a compact complex
manifold $H$.


For $x\in M$, let us 
denote by $J_0(x)$ (respectively $J_1(x)$) the almost complex operator of $\Cal F_0$ at $x$ (respectively of $\Cal F_1$). 
Let $L=(L_0)^{diff}=(L_1)^{diff}$ and $\Cal F=(\Cal F_0)^{diff}=(\Cal F_1)^{diff}$. We consider them as smooth sections of the vector bundle End $(T\Cal F)\to M$.
These operators are in general not equal, but they agree on $C$. As 
they are smooth and as $C$ is compact, we deduce that, for every neighborhood of the zero section of End $(T\Cal F)\to M$, there exists a  neighborhood 
$V\subset M$ of $C$ such that $f(V)\subset W$,
where $f$ is defined as the difference $J_0-J_1$.

Looking at the injections $L\hookrightarrow M$ and End $(TL)\hookrightarrow \text{End }(T\Cal F)$, this means that 
for every neighborhood of the zero section of End $(TL)\to L$, there exists a  neighborhood 
$V$ of the end $E$ of $L$ such that $f(V)\subset W$.
 
Consider now the inclusion diagram
$$
\CD 
L \aro >>> L^c=L\sqcup H^{diff} \cr
\aro V f VV \aro VVV\cr
\text{End }(TL) \aro >>> \text{End }(TL^c)
\endCD
$$

Setting a riemannian metric $\Vert - \Vert$ on $\text{End }(TL^c)$, we have, using the compacity of $H$, that, for all $\epsilon>0$,
there exists a neighborhood $V_{\epsilon}$ of $E$ in $L$ such that
$$ 
\sup_{x\in V_{\epsilon}} \Vert J_0(x)-J_1(x) \Vert \leq \epsilon 
$$ 
So $f$ extends continuously as the zero section over $H^{diff}$.

Now, assume that $J_0$ extends smoothly to $L^c$. We claim that $J_1$ extends {\it continuously} to 
$L^c$ with $J_0\equiv J_1$ on $H^{diff}$.

Indeed, take local coordinates in a neighborhood $U$ of $x\in H^{diff}$. We assume that End $(TL^c)$ is trivial over $U$, hence $J_0$ and $J_1$ are now functions with values in a euclidean space. We choose as norm $\Vert - \Vert$ over $U$ the euclidean norm. For $\epsilon>0$ and for $y$ sufficiently near $x$, we have
$$
\eqalign{
\Vert J_1(y) - J_0(x)\Vert &\leq \Vert J_1(y) - J_0(y)\Vert
+\Vert J_0(y) -  J_0(x)\Vert\cr
&\leq \Vert J_1(y) -  J_0(y)\Vert+\epsilon \quad\text{ by continuity of }J_0\cr
&\leq 2\epsilon
}
$$
and the claim is proved.

This is however not enough; we would like to show that $J_1$ extends {\it smoothly} to $L^c$. To do this, we use the fact that 
our foliations by complex leaves are by definition tame. This means that we may extend $\Cal F_0$ and $\Cal F_1$ to 
$W=M\sqcup C\times (0,1]$ by stating that the complex leaves of both $\Cal F_0$ and $\Cal F_1$ on $C\times (0,1]$ are $C\times \{t\}$.
So the complex structures $J_0$ and $J_1$ may be assumed equal on the collar $C\times [0,1]$.

We will now proceed by induction and repeat essentially the same argument in every $k$-jet bundle of sections of End $(TL^c)$ for $k\geq 0$.

We just do the case $k=1$. Consider the map 
$$
f_1\ :\ x\in M\longmapsto (x,j^1(f(x))\in J^1(M,\text{End }(T\Cal F))
$$
where $J^1(M,\text{End }(T\Cal F))$ is the bundle of $1$-jets of sections of End $(T\Cal F)$ and $j^1(f)$ is the $1$-jet of $f$.

Since $J_0$ is equal to $J_1$ not only on $C$ but on a collar of the boundary, we have that $j^1(f)$ is the zero section over $C$.

As before, considering the restriction of the situation to $L$ and seeing $f_1$ as a continuous map from $L$ to $J^1(L,\text{End }(TL))$, we have that for every neighborhood of the zero section of $J^1(L, \text{End }(TL))\to L$, there exists a  neighborhood 
$V$ of the end $E$ of $L$ such that $f_1(V)\subset W$.

Considering the inclusion of $L$ into $L^c$ and the corresponding inclusion for the jet bundles, and setting a riemannian metric
$\Vert - \Vert_1$ on $J^1(L^c,\text{End }(TL^c))$, we have, using the compacity of $H$, that, for all $\epsilon>0$,
there exists a neighborhood $V_{\epsilon}$ of $E$ in $L$ such that
$$
\forall k\geq 0,\qquad \sup_{x\in V_{\epsilon,k}} \Vert J_0(x)-J_1(x) \Vert_1 \leq \epsilon 
$$

So $f_1$ extends continuously as the zero section over $H^{diff}$. 
Now, take local coordinates in a neighborhood $U$ of $x\in H^{diff}$. We assume that End $(TL^c)$ is trivial over $U$, hence $J_0$ and $J_1$ are now functions with values in a euclidean space. We choose as norm $\Vert - \Vert_1$ over $U$ the maximum at one point of the euclidean norms of the function and all its first-order derivatives. The same sequence of inequalities as above but with the norm $\Vert -\Vert_1$ this time shows that $j^1(J_1)$ extends continuously on $L^c$ as $j^1(J_0)$ on $H^{diff}$.
It is now easy to deduce that $J_1$ admits a $C^1$ extension to $L^c$ as $J_0$ on $H^{diff}$. By induction, this extension is in fact smooth.

Now the almost complex operator $J_1$ on the whole $L^c$ is automatically integrable, since it is on the open and dense subset
$L\subset L^c$ so the Newlander-Nirenberg Theorem provides us with a complex structure.
$\square$ 
\enddemo

\remark{Remark}
Observe that the tame condition can be replaced by the somewhat weaker condition: on $W=M\sqcup C\times (0,1]$, both $J_0$ and $J_1$ extend
in such a way that they are equal on the collar $C\times [0,1]$. This allows to use the Compactification Lemma in some cases where the holonomy is not flat.
\endremark

\remark{Remark}
The Compactification Lemma compares two complex leaves of two different foliations as abstract complex manifolds. This implies that the 
compactification used may be arbitrary, i.e. does not depend on the foliations themselves. In particular, if $L_0$, as an abstract
complex manifold, admits various partial holomorphic compactifications at $C$-infinity which are topologically distinct, the Lemma works for every compactification and
$L_1$ will admit exactly the same number of partial holomorphic compactifications at $C$-infinity.
\endremark
\medskip

Here is an application of the Compactification Lemma. It shows how this Lemma can be used in some cases
to determine the biholomorphism type of the interior leaves of a foliation by complex leaves. This type of argument will be used many times 
in the next Section.

\example{Example}
Consider the solid torus endowed with the classical Reeb foliation. This foliation may be turned into a foliation $\Cal F_\tau$
by complex leaves, with boundary leaf biholomorphic to an arbitrary elliptic curve $\Bbb E_{\tau}$ and with all interior leaves bihilomorphic to $\Bbb C$
(see \cite{M-V}).

We want to prove, using the Compactification Lemma, that there does not exist on this fixed smooth foliation a complex structure with at least one interior leaf biholomorphic 
to the unit disk. Let $\Cal F$ be any complex structure on this fixed smooth foliation. The boundary leaf must be an elliptic curve, so this
foliation agrees on the boundary with some $\Cal F_\tau$. 

Let $L$ be an interior leaf of $\Cal F$ and let $L_\tau$ be the corresponding leaf of $\Cal F_\tau$. These leaves have one end whose 
limit set is the boundary leaf $\Bbb E_\tau$. Now, the leaf $L_\tau$ is biholomorphic to $\Bbb C$, so admits a holomorphic 
compactification as the Riemann sphere. 

The Compactification Lemma tells us then that $L$ also admits a holomorphic one-point compactification. By the uniformization Theorem, this
compactification is the Riemann sphere and $L$ is then a copy of $\Bbb C$. Hence, $\Cal F$ does not have any leaf biholomorphic to a disk.
\endexample
\medskip

The Compactification Lemma compares two different foliations. However, it is possible to use it with just one foliation.
We now explain this point.
\medskip
 
Let $(M,\Cal F)$ be a foliation by complex manifolds with $\partial 
M\not =\emptyset$. Let $L$ be a leaf of $\text{Int}(M)$ such that 
the limit set,  $\underset{i\geq1}\to\cap{\bar{\Cal U_i}}=C$ 
corresponding to a given end $\Cal U_1\supseteq\Cal 
U_2\supseteq\cdots$ of $L$ is a compact connected component $C$ of 
$\partial M$. There exists a neighborhood $V$ of $C$ and a global 
submersion from $V$ to $C$ which is the identity on $C$. For 
example, we may use a collar of $C$ in $M$ and define

$$ 
(z,t)\in C\times\Bbb R^+\simeq V\longmapsto z\in C \ .
$$ 
 
Reducing $V$ if necessary, we may assume that the previous 
submersion is a surjective local diffeomorphism when restricted to a 
leaf (intersected with $V$). In particular, $L\cap V$ is locally 
diffeomorphic to $C$. Call $(L\cap V)^{pb}$ the manifold 
diffeomorphic to $L\cap V$ with the complex structure induced from 
that of $C$ by pull-back and call $\Cal F^{pb}$ the corresponding 
complex structure on $V$. 
 
Notice that the complex structure of $(L\cap V)^{pb}$ is independent 
of the choice of the submersion. Indeed, for two different choices, 
the foliations $\Cal F^{pb}$ are CR-isomorphic since the submersions 
are isotopic by an isotopy which is the identity on $C$. 
\medskip

\proclaim{Corollary to Compactification Lemma} 
Assume that $(L\cap V)^{pb}$ 
admits a partial holomorphic compactification at $C$-infinity by 
adding a compact complex manifold $H$ of codimension-one. Then $L$ 
admits a partial holomorphic compactification at $C$-infinity by 
adding $H$. 
\endproclaim 
 
\demo{Proof} 
It is enough to prove that the pull-back foliated complex structure $\Cal F^{pb}$ is a foliation by complex manifolds, that is that it is
tame. The Compactification Lemma then applies. Now, since $\Cal F$ is tame by definition, we have that $(\Cal F)^{diff}$ is flat at the 
boundary (i.e. that the holonomies of the boundary components are smooth flat functions). This implies that the submersion used to define the
pull-back foliation is flat at the boundary, hence the pull-back complex structure is tame.
$\square$
\enddemo 
 
In the sequel, by Compactification Lemma, we will always mean the Corollary to the Compactification Lemma. We will also use the 
Compactification Lemma to provide uniform compactification of a ``tube'' of leaves.

\proclaim{Uniform Compactification Lemma}
In the same situation as in the previous Corollary, assume that $L$ has trivial holonomy. Choose a closed transverse section $s\simeq [-1,1]$
to $L$ and consider the tube of leaves $\Cal L\simeq_{diff} L^{diff}\times [-1,1]$ intersecting $s$. Assume that all the leaves in the tube accumulates uniformly onto $C$. Assume also 
that $(L\cap V)^{pb}$ admits a partial holomorphic compactification at $C$-infinity by 
adding a compact complex manifold $H$ of codimension-one. Then $\Cal L$ admits a partial CR compactification at $C\times [-1,1]$-infinity
by adding a compact CR manifold $H\times [-1,1]$.
\endproclaim

We omit the exact Definition of partial CR compactification. It should be clear from the Definition of partial holomorphic 
compactification.

\demo{Proof}
The proof is identical to the previous one. One just has to notice that the estimates work in a whole neighborhood of $C$ in $M$, so are valid
not only for a single leaf but for 
a tube of leaves.
$\square$
\enddemo

\example{Example}
We go back to the previous Example, where we show that a foliation by Riemann surfaces on a Reeb component has all interior leaves
biholomorphic to $\Bbb C$. A construction of the smooth underlying foliation is as follows (cf \cite{M-V}, Lemma 2). Let 
$$ 
X=\R^2 \times [0,\infty)\setminus \{(0,0,0)\} 
$$ 
foliated by the level sets of the projection onto the second factor. 
Let 
$$ 
h\ :\ (x,y,t)\in X \longmapsto (\alpha x-\beta y,\beta x+\alpha y, d(t))\in 
X 
$$ 
where $0<\alpha^2+\beta^2<1$ and where $d$ is a smooth function of 
$\R$ into $\R^+$, which is a contracting diffeomorphism of $\Bbb 
R^+$ and which is the identity on $\R^-$. Then the quotient of $X$ by the 
group generated by $h$ is a Reeb component.

The foliations $\Cal F_{\tau}$ above can be constructed by endowing $X$ with the trivial foliation 
$$
\C \times [0,\infty)\setminus \{(0,0,0)\} \ .
$$

Indeed, $h$ becomes multiplication by the complex number $\exp(2i\pi\tau)=\alpha+i\beta$ and
is a biholomorphism of the leaves. The foliation descends to the foliation $\Cal F_{\tau}$ previously described. 

We claim that, given {\it any} foliation by Riemann surfaces on this Reeb component, it is CR-isomorphic {\it in the interior} to some $\Cal F_{\tau}$.
We now want to use the uniform Compactification Lemma to prove this assertion (compare the following argument with the Example at the
end of Section 2). Endow the Reeb component with an arbitrary complex structure and endow the covering $X$ of the Reeb component by the 
pull-back complex structure. We already know that such a structure coincides on the boundary with some $\Cal F_{\tau}$ and that all the interior leaves of $X$ are 
biholomorphic to $\Bbb C$. Now the uniform Compactification Lemma tells us that a tube $\Cal L$ 
of interior leaves can be uniformly compactified 
as a closed manifold foliated by Riemann spheres. We may thus partially compactify $X$ as a product foliated covering $\bar X$ with compact interior leaves. 
But then Proposition 2 implies that the interior $\bar X$ is CR-isomorphic
to the trivial family $\Bbb P^1\times (0,\infty)$. Since the compactification is uniform, this CR-isomorphism sends 
Int $(\bar X\setminus X)$ to a {\it smooth} section $s$ of $\Bbb P^1\times (0,\infty)$. Since the automorphism group of $\Bbb P^1$ is transitive,
there exists a CR-isomorphism of $\Bbb P^1\times (0,\infty)$ sending $s$ onto the infinite section $\infty\times (0,\infty)$. Composing these
isomorphisms, this gives a CR-isomorphism between $\Bbb P^1\times (0,\infty)$ and Int $\bar X$ sending $\Bbb C\times (0,\infty)$ onto Int $X$.
We conclude by Corollary 1.
\endexample

Let us have a closer look to the situation of a product foliated covering $\pi:\tilde M\to M$ with boundary and deck transformation group isomorphic to $\Bbb Z$. It is diffeomorphic to
$L^{diff}\times [0,\infty)\setminus A^{diff}\times\{0\}$ 
and a generator of the deck transformation group has the following form
$$
(z,t)\in L^{diff}\times\Bbb R^+\setminus A^{diff}\times\{0\}\longmapsto (T(z,t), d(t))\in L^{diff}\times\Bbb R^+\setminus A^{diff}\times\{0\}
$$
As usual, assume that $d$ has no positive fixed point. As a consequence of the uniform compactification Lemma, we have:

\proclaim{Proposition 5}
Let $\pi:\tilde M\to M$ be a product foliated covering satisfying all the hypotheses described just above. Assume moreover that it is tame and that the boundary leaf $\partial\tilde M$ 
admits a partial holomorphic compactification at infinity by adding a compact curve $C$. 

Then, the whole $\tilde M$ admits a uniform compactification at infinity by adding $C$; that is, setting 
$$
\tilde M=(L^{diff}\times [0,\infty)\setminus A^{diff}\times\{0\}, J)
$$
then $J$ extends as a complex structure on $\tilde M\sqcup C\times\Bbb R^+$.
\endproclaim

\remark{Remark}
The boundary leaf $\partial \tilde M$ has two ends, being a $\Bbb Z$-covering of a compact manifold without boundary. On the other hand,
it is diffeomorphic to $L^{diff}\setminus A^{diff}$ and every other leaf is diffeomorphic to $L^{diff}$. We implicitely assume that
$L^{diff}$ has just one end, corresponding to the end of $\partial\tilde M$ called infinity in the statement of Proposition 5.
\endremark

\demo{Proof}
We take the same notations as in the proofs of the compactification Lemmas and consider $W=M\sqcup \partial M\times (-1,0)$. Observe that the covering
$\pi$ extends to a covering from $\tilde W=\tilde M\sqcup \partial\tilde M\times (-1,0)$. Consider also the partially compactified space
$\bar W=\tilde W\sqcup C\times (-1,\infty)$. Fix $\alpha>0$ and riemannian metrics on the bundles of $k$-jets of sections of
End $(T\bar W)$. From the inequalities used in the proof of the Compactification Lemma,
and taking the pull-back by the covering, we see that there exists, for all $\epsilon>0$ and for all $k\in\Bbb N$, a neighborhood $\tilde W_{\epsilon, k}$ of the end at infinity of 
$L^{diff}\times [0,\alpha]\setminus A^{diff}\times \{0\}$ such that
$$
\sup_{x\in \tilde W_{\epsilon, k}} \Vert \tilde J_0(x)-\tilde J^{pb}(x)\Vert_k\leq \epsilon
$$
where $\tilde J_0$ (respectively $\tilde J^{pb}$) are the pull-back complex structures of $J_0$ (respectively $J^{pb}$).

By assumption, $\tilde J^{pb}$ extends at infinity on $C\times \{0\}$, and thus uniformly on $C\times [0,\alpha]$ since $(\tilde M^{diff}, \tilde J^{pb})$ is CR-trivial by definition.
Now, the previous inequalities imply that this is also true for $\tilde J_0$ by arguing as in the proof of the Compactification Lemma.
$\square$
\enddemo

\remark{Remark}
As above, observe that the tame condition can be replaced by the somewhat weaker condition: on $W=M\sqcup \partial M\times (0,1]$, both $J_0$ and $J^{pb}$ extend
in such a way that they are equal on the collar $\partial M\times [0,1]$.
\endremark

\remark{Remark}
This Proposition can be considered as an extension result. In fact, we already know from the uniform compactification Lemma that the interior of $\tilde M$ admits the desired compactification.
Hence Proposition 5 states that this compactification can be extended uniformly to the boundary leaf. This is done following the rough argument that the added curves $C$ in the interior of
$\tilde M$ ``converge'' through the action of the deck transformation group onto the curve $C$ at the boundary. Hence the compactification is uniform.
  
To compare with, consider now the case of the other end of $\partial\tilde M$, i.e. assume that $\partial\tilde M=L\setminus A$. 
Then, as said before, it is not clear that the CR-structure of $\tilde M$ extends smoothly to $A$. 
\endremark

 

\head {\bf 5. The Lawson foliation does not admit any integrable 
almost CR-structure} 
\endhead 
 
\proclaim{Theorem A} The Lawson foliation of $\Bbb S^5$ can be 
endowed with a compatible almost CR-structure, however this structure can 
never be integrable. 
\endproclaim 
 
\demo{Proof} Let us first prove that the Lawson foliation admits 
(non-integrable) Levi-flat CR-structures. Notice that the inclusion 
of $\Bbb S^5$ as the unit Euclidean sphere of $\Bbb C^3$ defines the 
canonical integrable almost CR-structure on $\Bbb S^5$. This structure is 
not Levi-flat since the corresponding distribution is the canonical 
contact structure of $\Bbb S^5$. At $z\in\Bbb S^5\subset \Bbb C^3$, 
it is equal to 
$$ 
\{w\in\Bbb C^3\quad\vert\quad \langle z,w\rangle=0\} 
$$ 
where the angles denote the standard hermitian product of $\Bbb 
C^3$. 
 
This contact structure is orthogonal to the contact flow which gives 
the Hopf fibration of  $\Bbb S^5$ by circles. It is generated by the 
unit vector field 
$$ 
v \ : \ z\in\Bbb S^5\subset\Bbb C^3 \longmapsto v(z)=iz\ . 
$$ 
 
It is enough to prove that the distribution $\Cal H$ tangent to the 
Lawson foliation is homotopic to the distribution of this contact 
structure. This implies that $\Cal H$, as an abstract vector bundle over $\Bbb S^5$, is homotopic to a complex vector bundle. 
Therefore, it is a complex vector bundle by \cite{St, Theorem 11.5}, i.e. it has an almost CR-structure.
\medskip

We claim that there are only two homotopy classes of unit vector fields over $\Bbb S^{5}$ and that they are represented by $v$ from the one hand, and by $-v$ from the
other hand. As a consequence, every $4$-planes orientable distribution is homotopic to the distribution of the contact structure, once we take on it the "good" orientation. So
admits an almost CR-structure.
\medskip

Let us first prove that there are only two homotopy classes of unit vector fields over $\Bbb S^{5}$. Indeed, this is exactly the number of homotopy classes of sections
of the unit tangent bundle $U$ of $\Bbb S^5$. Since it has a section (the previous vector field $v$), the homotopy sequence of the fibration splits at each step:
$$
\pi_{i}(U)=\pi_{i}(\Bbb S^{5})\oplus\pi_{i}(\Bbb S^{4})
\leqno i>0
$$
In particular, for $i=5$, this gives
$$
\pi_{5}(U)=\Bbb Z\oplus\Bbb Z_{2}
$$  
Given an element $j$ of $\pi_{5}(U)$, notice that its component in $\Bbb Z$ is the degree of the composition map
$$
\Bbb S^5\aro > j >> U\aro >\text{bundle projection}>>\Bbb S^{5} 
$$
Therefore, a homotopy class of sections of $U$ is exactly an element of $\pi_{5}(U)$ whose component in $\Bbb Z$ is one and the claim is proved.
\medskip

To finish with, it is enough to prove that $v$ and $-v$, the two unit vector fields of $\Bbb S^{5}$ we know, are not homotopic. Assume the contrary. Then the contact
distribution of $\Bbb S^{5}$ and the same distribution with the orientation reversed would be homotopic. Consider the previously described action of $\Bbb S^{1}$ onto
$\Bbb S^{5}$. It leaves both distributions invariant. By \cite{Br, Chapter VI, Theorem 3.1}, we may assume that the homotopy between these two distributions is equivariant. So there is an equivariant oriented isomorphism
between these two distributions \cite{Br, Chapter II, Theorem 7.4}. Hence it descends
to an oriented isomorphism between the tangent bundle of the complex projective space $\Bbb P^2$ and the tangent bundle of $\overline{\Bbb P^{2}}$, the manifold obtained
from $\Bbb P^{2}$ by reversing its orientation. Taking account of what we said above, this would imply that  $\overline{\Bbb P^{2}}$ has an almost complex structure.
But this is known to be false, by use of Wu's theorem characterizing homologically the existence of an almost complex structure on a real $4$-manifold, see \cite{B-H-P-V}.
This finishes the proof.

\medskip 
 
Let us prove now that the almost CR-structures compatible with the Lawson 
foliation can never be integrable. Assume the contrary, i.e. assume 
the existence of such an integrable almost CR-structure. In fact, we will 
prove that it cannot exist on the interior part of the Lawson 
foliation (that is, with the notations of 1.2, the neighborhood 
$\Cal N$ of $K$). Recall that the boundary $\partial\Cal N\simeq K\times \Bbb S^1$ is a leaf and that the non-compact leaves of $\Cal N$ are all 
diffeomorphic to 
$$ 
L=\R^2/\Z^2 \times \Bbb D 
$$ 
and are equipped with a complex structure $J$ by our assumption. 
\medskip 
 
The proof will take the form of several Lemmas. We will describe 
some properties of $J$ imposed by the topology of the foliation; 
these properties will lead to a contradiction. The outline of the 
proof is the following: in Lemmas 1 and 2, we describe the explicit 
biholomorphism type of the boundary leaf, using the Enriques-Kodaira classification. Then, in Lemmas 3 and 4, we prove that the non-compact
leaves are all biholomorphic to the product of $\Bbb C$ by a fixed elliptic curve. We use the Compactification Lemma to obtain such a result. 
Finally, we exhibit an associated CR-trivial product foliated covering and compute a generator of its deck transformation group. This gives us an automorphism of the non-compact
leaf.
 
The contradiction comes now from the fact that the previous 
biholomorphism is not compatible with the complex structure of a non-compact leaf; the map which should be a biholomorphism of the non-compact leaf does not belong to the automorphism group of the leaf for homological reason (its 
action on the homology groups is different from the action induced by any biholomorphism).
\medskip 
We first need to review Lawson's construction with a little more 
care. We refer to \cite{La} for details. 
\medskip 
Recall the bundle map 
$$ 
K\longrightarrow \Bbb S^1_1\times \Bbb S^1_2 
$$ 
induced from the map $W\to \Bbb E_{\omega}$ by passing to the 
associated unit bundle (see Section 1.2). 

\remark{Remark}
It will be important in the sequel to distinguish the different $\Bbb S^1$-factors, so we label them.
\endremark
\medskip

Define the map 
$$ 
\pi \ :\ K\times\Bbb S^1_3\longrightarrow \Bbb S^1_1\times \Bbb S^1_2 
$$ 
by composition of the natural projection and of the previous map. 
The important fact is that $\pi$ defines $K\times\Bbb S^{1}$ as a 
{\it principal} torus bundle over $\Bbb S^1\times\Bbb S^1$. 
 
\remark{Remark} Recall the other bundle map 
$$ 
s \ :\ K\longrightarrow \Bbb S^1_1 
$$ 
obtained from the description of $K$ as a suspension. Then the map 
$$ 
(s,Id)\ :\ K\times\Bbb S^1_3\longrightarrow \Bbb S^1_1\times\Bbb S^1_3 
$$ 
is a smooth submersion with compact and connected fibers, so is a 
locally trivial smooth fiber bundle by Ehresmann's Lemma ; and it is 
not isomorphic to $\pi$. But the key point for us is that this 
bundle is {\it not principal}, as previously remarked for the bundle 
map $s$. 
\endremark 
\medskip 
 
Starting with the map $\pi$, Lawson composes it with the natural 
projection onto one of the $\Bbb S^1$-factors and extends it as a 
submersion 
$$ 
p : K\times\overline{\Bbb D}_4\longrightarrow \Bbb S^1_1 
$$ 
 
He foliates the interior $K\times\Bbb D$ by the level sets of 
$p$. Notice that the choice of the projection does not matter. 
It is easy to prove that exchanging the factors of the base of $\pi$ lifts to a diffeomorphism of the total space
$K\times\Bbb S^{1}$ of $\pi$, hence the foliations obtained by the two different projections are diffeomorphic. Notice also that the two different 
projections come from projections $K\to\Bbb S^1$ both of which 
define $K$ as the suspension of a torus by the matrix $A$. {\it So 
we make the assumption that the image of $p$ is the same as the 
image of the suspension map $s$.} 
\medskip

Consider now the following commutative diagram 
\medskip 
$$ 
\CD 
K\times \overline{\Bbb D}_4\setminus\{0\} \aro >\simeq>> K\times\Bbb S^1_4\times (0,1] \\ 
\aro V Id VV \aro VV (p_{\vert \partial}, Id) V\\ 
K\times \overline{\Bbb D}_4\setminus\{0\} \aro > \bar p>> \Bbb S^1_1\times (0,1]\simeq \overline{\Bbb D}_1\setminus\{0\} \\ 
\aro V Id VV \aro VV \text{natural projection} V \\ 
K\times \overline{\Bbb D}_4\setminus\{0\} \aro > p >> \Bbb S^1_1 
\endCD 
$$ 
\medskip 
where $\bar p$ is defined from the other arrows. 
 \medskip

In this setting, the leaves of the foliation restricted to 
$K\times{\Bbb D}_4\setminus\{0\}$ are also given by the inverse images 
$$ 
\bar p^{-1}(\{\exp{i\theta}\}\times (0,1]) \eqno \theta\in\Bbb R 
$$ 
 
Now Lawson turbulizes this foliation by considering as leaves the 
inverse images $\bar p^{-1}(\Cal C_{\theta})$, where $\Cal C_\theta$ 
spirals in the disk $\Bbb S^1\times (0,1]$ as shown in the following 
picture. 
 
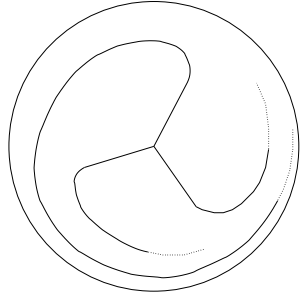
\hfil\scaledpicture 6.1in by 6.1in (spirale scaled 250) \hfil

This turbulized foliation extends as a foliation of $K\times\overline{\Bbb D}$ 
such that the boundary $K\times\Bbb S^1$ is a leaf. 
 \medskip
 
We may now start with the proof 
 
\proclaim{Lemma 1} The compact leaf $(K\times \Bbb S^1,J)$ is a 
primary Kodaira surface, that is a principal holomorphic fiber 
bundle over an elliptic curve with fiber an elliptic curve. 
\endproclaim 
 
We denote by $\Cal S$ this complex compact surface, by $\Bbb 
E_{\alpha}$ the base space of the associated bundle map and by $\Bbb 
E_{\beta}$ the fibers of this map. 
 
\demo{Proof of Lemma 1} From the construction of $\Cal F_{\Bbb C}$, 
we know that the smooth model of $\Cal S$, that is $K\times \Bbb 
S^1$ admits a structure of a primary Kodaira surface. This implies 
that the Chern numbers $c_1^2$ and $c_2$ of this structure are zero 
and that the first Betti number is $3$ (a fact that can also be 
easily recovered from the description of $K$ as a suspension given 
in 1.2). As these numbers are topological invariants, they keep the 
same values for $\Cal S$. On the other hand, the universal covering 
of $K\times\Bbb S^1$ is $\Bbb R^4$, which implies, using for example 
the long exact sequence in homotopy of a fibration, that the second 
homotopy group of $\Cal S$ is zero. Hence $\Cal S$ is minimal. The 
Enriques-Kodaira classification ([B-H-P-V], p.244) shows that $\Cal 
S$ is a primary Kodaira surface or a minimal properly elliptic 
surface. By \cite{F-M, Theorem S3, (ii)}, a smooth manifold cannot 
admit at the same time a complex structure of Kodaira dimension zero 
and another one of Kodaira dimension one. $\square$ 
\enddemo

\proclaim{Lemma 2} The complex structure of $\Cal S$ is compatible 
with $\pi$, that is there exists a structure of elliptic curve $\Bbb 
E_{\alpha}$ on the base space of $\pi$ such that $\pi$ becomes a 
holomorphic principal fiber bundle from $\Cal S$ to $\Bbb 
E_{\alpha}$. 
\endproclaim

\demo{Proof of Lemma 2} By Lemma 1, there exists a holomorphic 
principal elliptic fiber bundle $\Cal S\to \Bbb E_{\alpha}$, so that 
it is enough to prove that this bundle is smoothly isomorphic (that 
is isomorphic as $C^{\infty}$ principal bundles) to the bundle 
$\pi$: endowing $K\times\Bbb S^1$ with the complex structure 
obtained by pull-back by this isomorphism, it becomes a holomorphic 
principal elliptic bundle which is complex isomorphic to $\Cal S\to 
\Bbb E_{\alpha}$. 
 
A principal elliptic fiber bundle over an elliptic curve $\Bbb 
E_{\alpha}$ is obtained from a $\Bbb C^*$-principal bundle over 
$\Bbb E_{\alpha}$ by taking the quotient by a group acting as a 
complex homothety in the fibers. Notice that this description fits 
not only to the case of a holomorphic bundle but also to the case of 
a smooth bundle. Indeed a smooth principal elliptic bundle over 
$\Bbb E_\alpha$ can be thought of as a smooth bundle with complex 
fibers: fixing a structure of an elliptic curve on the fibers, it is 
preserved by the structural group which, by definition, consists 
only of translations. Then such an elliptic bundle is obtained from 
a smooth $\Bbb C^*$-bundle (that is a smooth locally trivial bundle 
over $\Bbb E_\alpha$ with fiber $\Bbb C^*$ and structural group 
$\Bbb C^*$) by taking the quotient by a group acting as a complex 
homothety in the fibers. In the sequel of the proof, by elliptic 
bundle (respectively $\Bbb C^*$-bundle), we mean {\it smooth} ones. 
 
To this $\Bbb C^*$-bundle, we may associate its unit bundle, which 
is an oriented principal circle bundle over $\Bbb E_{\alpha}$. Two 
such elliptic fiber bundles are isomorphic if and only if the 
corresponding $\Bbb C^*$-bundles are isomorphic and this occurs if 
and only if the associated oriented circle bundles are isomorphic. 
Finally this is the case if and only if the Euler numbers are equal. 
 
On the other hand, from the description of $K$ given in Section 1.2 
and in \cite{M-V}, it is straightforward to check that, for 
$n\in\Bbb Z$, the isomorphism class of circle bundles of Euler 
number $n$ can be represented, as an oriented manifold, as the 
suspension of a real torus by the matrix 
$$ 
A_n=\pmatrix 1 &n \\ 
0 &1 
\endpmatrix 
$$ 
 
From this, the first homology group of a circle bundle of Euler 
number $n$ is isomorphic to $\Bbb Z\oplus\Bbb Z\oplus\Bbb Z_{\vert 
n\vert}$. Hence if two such circle bundles are diffeomorphic {\it as 
a manifold}, then their Euler numbers are equal up to sign. In 
particular the Euler numbers of $\pi$ and of $\Cal S\to \Bbb 
E_{\alpha}$ differ at most by a sign. Now, changing the orientation 
of a circle bundle changes the sign of its Euler class, so we may 
assume that the Euler numbers of $\pi$ and of $\Cal S\to \Bbb 
E_{\alpha}$ are the same by choosing the ``right'' orientation on 
$\pi$. 
 
From all that preceeds, it follows that $\Cal S\to \Bbb E_{\alpha}$ 
is smoothly isomorphic to $\pi$, and the Lemma is proved. $\square$ 
\enddemo

\proclaim{Lemma 3} Let $L$ be a complex non-compact leaf of $\Cal 
N$. Then $L$ admits a holomorphic compactification as a ruled surface of genus $1$ by adding an elliptic curve $\Bbb E_{\beta}$. 
\endproclaim 
 
\demo{Proof of Lemma 3} As explained above, the Lawson foliation has 
leaves $\bar p^{-1}(\Cal C_{\theta})$, where $\Cal C_\theta$ is the 
curve of $\Bbb S^1\times (0,1]$ previously drawn. Now, we may view the open set given by the intersection
of a curve $\Cal C_{\theta}$ with $\Bbb S^1_1\times (1/2,1)$ as an open set of a $\Bbb Z$-covering of the 
circular boundary $\Bbb S^1_1\times\{1\}$. From this, we infer that a 
non-compact leaf $L^*$ of $K\times\Bbb S^1_4\times (1/2,1)$ is an open set of the $\Bbb 
Z$-covering of $K\times\Bbb S^1_4$ obtained by unrolling the circle $\Bbb S^1_1$ following the diagram 
$$ 
\CD 
\Cal C_{\theta}\cap \Bbb S^1_1\times (1/2,1)\aro >>> \Bbb S^1_1 \\ 
\aro A \bar p AA \aro AA p_{\vert\partial} A \\ 
L^* \aro >>> K\times\Bbb S^1_4 
\endCD 
$$ 
 
Notice that $L^*$ being the intersection of a leaf $L\simeq \Bbb 
R^2/\Bbb Z^2\times\Bbb D$ with $K\times\Bbb S^1\times (1/2,1)$, it is 
diffeomorphic to $\Bbb R^2/\Bbb Z^2\times\Bbb S^1\times \Bbb D\setminus\{0\}$. 
 
Let us pass now to the complex world. Then, putting on $L^*$ its 
pull-back structure, we obtain an open set of a holomorphic $\Bbb Z$-covering 
satisfying the diagram (see \cite{M-V, p.921--922} for the 
holomorphic triviality of the pull-back bundle) 
$$ 
\CD 
(L^*)^{pb}\subset\Bbb C^*\times\Bbb E_\beta \aro > \Bbb Z-\text{covering} >> \Cal S \\ 
\aro VVV \aro V \pi VV \\ 
\Bbb C^* \aro > \Bbb Z-\text{covering} >> \Bbb E_\alpha 
\endCD 
$$ 
 
Since $L^*$ is asymptotic to $\Bbb C^*\times \Bbb E_\beta$ at $\infty$, the compactification Lemma ensures us that $L^*$ can be compactified 
by adding an elliptic curve $\Bbb E_\beta$ in two ways (exactly as it is the case for $\Bbb C^*\times \Bbb E_\beta$). 
We {\it choose} the compactification as follows (cf the remark following the Compactification Lemma)
$$ 
(L^*)^c=(\Bbb R^2/\Bbb Z^2\times\overline{\Bbb D}\setminus\{0\})\cup (\Bbb 
R^2/\Bbb Z^2\times\overline{\Bbb D}) 
$$ 
so that the complex model $(L,J)$ of a non-compact leaf of $\Cal N$ 
admits a holomorphic compactification 
$$ 
L^c=(\Bbb S^1\times\Bbb S^1\times\overline{\Bbb D})\cup (\Bbb S^1\times \Bbb S^1\times\overline{\Bbb 
D})=\Bbb S^1\times\Bbb S^1\times \Bbb S^2 
$$ 
equipped with an extension $J^c$ of $J$. 
 
Using the Enriques-Kodaira classification, we immediately obtain 
that $(L^c, J^c)$ is a ruled surface of genus $1$ or an elliptic surface. Since $L^c$ is diffeomorphic to the rational surface formed by the product of an elliptic curve by $\Bbb P^1$, we conclude from \cite{F-M, Chapter II,
Theorem S.3. (i)} that we are in the first case.
 $\square$ 
 
 
 
 
 
 
 
 
 
\enddemo 

Here comes the most difficult Lemma. 

\proclaim{Lemma 4}
The interior leaves are all biholomorphic to $\Bbb C\times\Bbb E_{\beta}$
\endproclaim

\demo{Proof of Lemma 4}
 Recall the map (see 1.2) 
$$ 
A \ : \ [x,y]\in\Bbb R^2/\Bbb Z^2\longmapsto [ x-3y,y]\in\Bbb 
R^2/\Bbb Z^2\ . 
$$ 

Consider now the covering
$$
\CD
\Bbb R^{2}/\Bbb Z^{2}\times\Bbb R\times\overline{\Bbb D}_{4}\setminus\{0\} \aro >c>> K\times \overline{\Bbb D}_{4}\setminus\{0\}\\
\aro V q VV \aro VV \bar p V \\
\Bbb R\times (0,1] \aro >>> \overline{\Bbb D}_{1}\setminus\{0\}\simeq \Bbb S^{1}_{1}\times (0,1]
\endCD
$$
whose deck transformation group is generated by the map
$$
([x,y],t,z)\in\Bbb R^{2}/\Bbb Z^{2}\times\Bbb R\times\overline{\Bbb D}_{4}\setminus\{0\} \aro > T>> (A[x,y],t+1,z)\in
\Bbb R^{2}/\Bbb Z^{2}\times\Bbb R\times\overline{\Bbb D}_{4}\setminus\{0\}\ .
$$
The pull-back by $c$ of the Lawson foliation restricted to $K\times \overline{\Bbb D}^{4}\setminus\{0\}$ is given by the inverse image by $q$
of the following foliation of the strip $\Bbb R\times (0,1]$
\medskip
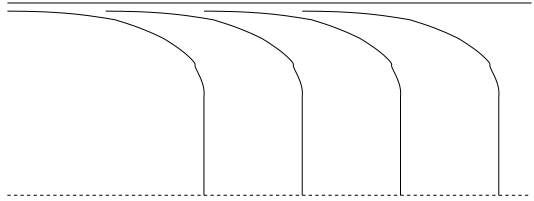
\hfil\scaledpicture 8.0in by 2.9in (strip scaled 350) \hfil 

\medskip
\noindent invariant by the horizontal translation $(t,s)\mapsto (t+1,s)$.
As usual, we consider the pull-back complex structure on this foliation. 


We know from the previous Lemmas that the boundary leaf is biholomorphic to $\Bbb E_\beta\times\Bbb C^*$. Notice that
the leaves we consider have now two ends, since we restrict the foliation to $K\times\overline{\Bbb D}\setminus\{0\}$. We are interested in the end which accumulates onto the boundary.
The restriction of the map $c$ of the previous covering to $q^{-1}(\Bbb R\times (1/2,1])$ sends the foliation restricted to $q^{-1}(\Bbb R\times (1/2,1])$ to the initial foliation with leaves $\Cal C_\theta\cap (1/2,1)$.
In other words, we may make use of the Compactification Lemma on this covering
to partially compactify the leaves as open subsets of ruled surfaces of genus $1$ as in Lemma 3; but in this new context, we may also make use of the uniform Compactification Lemma
to partially compactify a tube of leaves in a uniform way. Observe that the elliptic curve $\Bbb E_\beta$ we add has to be a section for topological reasons.

\remark{Remark}
It should be pointed out that the covering we use here {\it is not} a product foliated covering. Nevertheless, this covering has a special property that we will use. 
Indeed, as said before, the leaves have two ends, but in the quotient
{only one end} accumulates onto the boundary.  And what is more important, this end accumulates onto the boundary not only in
the quotient space but also in the total space of the covering. The important consequence is that we may assume that the tube of leaves we compactify contains the boundary leaf, even if we
are not under the hypotheses of Proposition 5.
\endremark
\medskip

We thus form a compactified tube of leaves $\Xi\simeq \Bbb R^2/\Bbb Z^2\times \Bbb D\times\Bbb [0,1]$ such that

\medskip
\noindent (i) the leaf $X_1$ (boundary leaf) is $\Bbb C\times\Bbb E_{\beta}$.

\noindent (ii) the other leaves are open dense subsets of ruled surfaces with a holomorphic section $\Bbb E_\beta$.

\noindent (iii) there is a CR-injection
$$
\CD
\Bbb E_\beta\times [0,1] \aro >>> \Xi \\
\aro VVV \aro VVV\\
[0,1] \aro > Id>> [0,1]\\
\endCD
$$
\medskip

Notice that the last point is a direct consequence of the fact that the compactification is uniform. Consider now the normal bundles of $\Bbb E_\beta$ in
each fiber. Point (iii) implies that they fit into a deformation family $\Cal Y\to [0,1]$ of line bundles over $\Bbb E_\beta$ (in the sense of Section 3), all of which are topologically
trivial. Using the fact that $T$ sends a leaf of $\Cal X$ CR-isomorphically to another leaf which is closer from the boundary, we have that
the family $\Cal Y$ satisfies the hypothesis of the Dumping Lemma, with boundary fiber $\Cal Y_1$ isomorphic to $\Bbb C\times\Bbb E_\beta$.

Going back to the Lawson foliation, all this means that an interior leaf is obtained from a ruled surface of genus $1$ by removing a section
with {\it holomorphically trivial normal bundle}. Therefore the ruled surface is a product and the interior leaves are biholomorphic to $\Bbb C\times \Bbb E_\beta$.
$\square$
\enddemo

 We are now in position to terminate the proof.

 

The open set Int $\Cal N$ identifies with 
$K\times \Bbb D_4$. From 1.2, we have a $\Bbb Z$-covering 
$$ 
\R^2/\Z^2 \times \R \times \Bbb D_4\longrightarrow K\times \Bbb D_4 
$$ 
whose deck transformations group is generated by 
$$ 
h\ :\ ([x,y], t,w)\longrightarrow (A [x,y],t+1, w) 
$$ 
Consider the trivial foliation of $\R^2/\Z^2 \times \R \times \Bbb 
D_4$ given by the level sets of the projection onto the $\R$-factor. 
This foliation is invariant by $h$ and, from what preceeds, descends 
on $K\times \Bbb D_4$ as the foliation used by Lawson before 
turbulization (recall that we made the assumption that the images of 
$p$ and of $s$ are the same). As this foliation is diffeomorphic to 
the final one, it is endowed with an integrable almost CR-structure by our 
assumption. We may take the pull-back of this structure and obtain 
thus an integrable almost CR-structure on $\R^2/\Z^2 \times \R \times \Bbb 
D_4$ such that the projection map is CR. In other words, it becomes a 
product foliated covering. Lemma 4 implies that the leaves of this covering are 
all biholomorphic to $\Bbb C\times\Bbb E_\beta$. Moreover, since it is CR-isomorphic 
as a covering space to the covering used in the proof of Lemma 4 (in restriction to $\Bbb D_4\setminus\{0\}$), the leaves can be compactified 
uniformly as $\Bbb P^1\times\Bbb E_\beta$. This compactified foliated covering is CR-trivial by Proposition 2.
Since the compactification is uniform and since the automorphism group of $\Bbb P^1\times\Bbb E_\beta$ is transitive on
the $\Bbb P_1$-factor, the initial product foliated covering is CR-trivial. It follows now from Proposition 3 that the
monodromy of this covering is well-defined as a biholomorphism of $\Bbb C\times\Bbb E_\beta$. Since the smooth monodromy is given by $(Id,A)$, this means that $\Bbb C\times\Bbb E_\beta$
should admit a biholomorphism smoothly isotopic to $(Id, A)$, and in particular which acts as the matrix $A$ on the first homology group $H_1(L,\Bbb Z)\simeq \Bbb Z^2$. However, it is an easy matter
to check that every automorphism of $\Bbb C\times\Bbb E_\beta$ is isotopic to $(Id, S)$ where $S$ is an automorphism of $\Bbb E_\beta$ and
that this is not the case of $(Id,A)$. 
$\square$ 
\enddemo

\head {\bf 6. The set of integrable almost CR-structures of a smooth 
foliation} 
\endhead 
 
Let $(M,\Cal F^{diff})$ be the pair consisting of a smooth manifold 
$M$ without boundary and a smooth codimension-one foliation $\Cal 
F^{diff}$ on it. We call {\it complex structure} on $(M,\Cal F^{diff})$ the data of an integrable almost CR operator $J_{\Cal F}$ on the tangent bundle
to the foliation $T\Cal F^{diff}$. Of course a complex structure corresponds to a foliation $\Cal F$ by complex manifolds of $M$
whose underlying smooth foliation is $\Cal F^{diff}$. 

Two complex structures $J_{\Cal F}$ and $J_{\Cal F'}$ of $(M,\Cal F^{diff})$ are {\it equivalent} if there exists a foliated diffeomorphism of
$(M,\Cal F^{diff})$ whose differential commutes with $J_{\Cal F}$ and $J_{\Cal F'}$. Two equivalent complex structures give rise to CR-isomorphic foliations
by complex manifolds $\Cal F$ and $\Cal F'$.

 
 
\proclaim{Definition} We say that two complex 
structures $J_{\Cal F}$ and $J_{\Cal F'}$ are strongly equivalent if there 
exists a CR-isomorphism between $(M, J_{\Cal F})$ and $(M, J_{\Cal 
F'})$ which does not exchange the leaves, i.e. which descends to the 
identity of the leaf space $M/\Cal F^{diff}$. 
\endproclaim 
 
We are now in position to give the principal definitions of this 
Section. Assume from now on that $M$ is compact. 
 
\proclaim{Definition} We call space of complex structures of 
$(M,\Cal F^{diff})$ the set of classes of strongly equivalent complex structures on $(M,\Cal F^{diff})$. 
 
We denote it by $\Cal C(M,\Cal F^{diff})$. 
\endproclaim 
 
\remark{Remark} In the same way, we may define a set $\Cal C(M, \Cal 
F^{top})$ by considering, modulo CR-isomorphisms which do not 
exchange the leaves, the set of classes of foliations by complex 
manifolds of $M$ {\it homeomorphic} to a fixed model $\Cal F^{top}$. 
We will see that it can be completely distinct from $\Cal C(M,\Cal 
F^{diff})$. 
\endremark 
 
\remark{Remark} We could also define $\Cal C(M,\Cal F^{diff})$ 
(respectively $\Cal C(M, \Cal F^{top})$) as the set of classes of 
foliations by complex manifolds of $(M,\Cal F^{diff})$ modulo 
CR-isomorphisms, that is allowing also CR-isomorphisms which 
exchange the leaves. We will see in the second example of this 
Section that it may lead to different spaces. 
\endremark 
\medskip 
 
Note that $\Cal C(M,\Cal F^{diff})$ can be empty even if each leaf 
of $\Cal F^{diff}$ can be endowed separately with a complex 
structure: by Theorem A, this is the case for the Lawson foliation. 
 
Note that when non empty, $\Cal C(M,\Cal F^{diff})$ does not have in 
general a structure of a complex manifold; if it has, it may be 
finite or infinite dimensional. 
 
We turn now to the definitions of deformation families and moduli 
spaces for foliations by complex manifolds (compare with \cite{Su, 
p.138--139}). 
 
\proclaim{Definition} Let $X$ be a complex manifold. Let $(V,J)$ be 
a smooth manifold endowed with a codimension-one integrable and 
Levi-flat almost CR-structure. Let $\pi : V\to X$ be a CR-submersion. Then, 
we say that $\pi$ is a deformation family of $(M,\Cal F^{diff})$ if, 
for every $x\in X$, 
\medskip 
\noindent (i) the manifold $\pi^{-1}(\{x\})$ is diffeomorphic to 
$M$. 
 
\noindent (ii) the almost CR-structure $J$ of $V$ defines a 
codimension-one, integrable and Levi-flat almost CR-structure $J_x$ on the 
fiber $\pi^{-1}(\{x\})$. 
 
\noindent (iii) there is a foliated diffeomorphism between $(V,\Cal F)$ and $M\times X$ endowed 
with the foliation $\Cal F^{diff}$ on each $M\times\{pt\}$ such that we have a diagram
$$
\CD 
(V,\Cal F) \aro >\simeq>> M\times X\cr
\aro V \pi VV \aro VV \text {2nd projection}V \cr
X \aro > Id >> X
\endCD
$$
\endproclaim 
 
In the case of a foliation having a compact leaf, there is a direct 
relationship between deformation families of the foliation and 
deformation families of the compact leaf. 
 
\proclaim{Proposition 6} Let $\pi : V\to X$ be a deformation family 
of $(M,\Cal F^{diff})$. Let $L$ be a compact leaf of $\Cal 
F^{diff}$. Then $\pi$ induces a holomorphic deformation family $W\to 
X$ of $L$. 
\endproclaim 
 
\demo{Proof} Consider the foliation by complex manifolds $\Cal G$ of 
$V$ induced by $J$. For $x$ varying in $X$, the union of the compact 
submanifolds of $\pi^{-1}(\{x\})$ corresponding to the compact leaf 
$L$ via point (iii) of the Definition is a leaf $W$ of $\Cal G$, and therefore a complex manifold. The 
restriction of $\pi$ to $W$ is a holomorphic submersion onto $X$ and 
defines thus a holomorphic deformation family of $L$. $\square$ 
\enddemo 
 
Given a deformation family $\pi : V\to X$ of $(M,\Cal F^{diff})$, 
there is a natural map $\alpha_{\pi}$ from $X$ to $\Cal C(M,\Cal 
F^{diff})$: just define $\alpha_{\pi}(x)$ to be the point of $\Cal 
C(M,\Cal F^{diff})$ corresponding to the CR-structure of 
$\pi^{-1}(\{x\})$. 
 
\proclaim{Definition} Assume that $\Cal C(M,\Cal F^{diff})$ can be 
endowed with the structure of a complex manifold $\Cal C_{\C}$. Let 
$i$ be the identification map between $\Cal C(M,\Cal F^{diff})$ and 
$\Cal C_{\C}$. Then $\Cal C_{\C}$ is called a coarse moduli space if 
\medskip 
\noindent (i) given any deformation family $\pi : V\to X$ of 
$(M,\Cal F^{diff})$, the natural map $i\circ \alpha_{\pi}$ is 
holomorphic as map into $\Cal C_{\C}$. 
 
\noindent (ii) the pair $(\Cal C_{\C},i)$ is unique up to 
composition with a biholomorphism of $\Cal C_{\C}$. 
\medskip 
If a coarse moduli space exists, we denote it by $\Cal M(M,\Cal 
F^{diff})$. 
\endproclaim 
 
\remark{Remark} In the classical case of moduli spaces of compact 
complex manifolds, a coarse moduli space as well as the base of a 
deformation family are usually not assumed to be a manifold but only 
a complex analytic space (with some special properties). Here, we 
may easily modify our Definition of coarse moduli space in this way. 
However, the exact version of the singular deformation families in 
our case is not very clear: if we take in the definition of a 
deformation family a singular base $X$, then $(V,J)$ should also be 
singular. But $V$ is real, and {\it cannot be taken real analytic 
but only smooth}, since real analytic codimension-one foliations do 
not exist on a simply-connected compact manifold by a classical 
result of Haefliger \cite{Hae}. 
\endremark 
\medskip 
 
Given a deformation family $\pi : V\to X$ and a holomorphic map 
between complex manifolds $f : Y\to X$, we may of course define a 
pull-back family $f^*\pi$.

\proclaim{Definition} Assume that $(M,\Cal F^{diff})$ has a coarse 
moduli space $\Cal M(M,\Cal F^{diff})$. Then it is a fine moduli 
space if there exists a deformation family $\Pi : W\to \Cal M(M,\Cal 
F^{diff})$ of $(M,\Cal F^{diff})$ such that every deformation family 
$\pi : V\to X$ of $(M,\Cal F^{diff})$ coincides with 
$\alpha_{\pi}^*\Pi$. 
\endproclaim 
 
In the classical case of compact complex manifolds, it is too much 
to expect to have a fine moduli space and the corresponding 
Definition is not pertinent. Indeed, if a fine moduli space exists 
for a smooth compact manifold $X$, then every locally trivial 
holomorphic bundle with fiber diffeomorphic to $X$ has to be 
holomorphically trivial, since it defines a holomorphic family of 
deformations of the fiber which is the pull-back of a point in the 
moduli space. But this is very restrictive, and {\it we do not know 
of any smooth compact manifold admitting a (non-empty) fine moduli 
space}. 
 
In our case, similarly, if a fine moduli space exists for $(M,\Cal 
F^{diff})$, then every deformation family whose fibers are all 
biholomorphic must be trivial, that is biholomorphic to a product. 
For this reason, the examples we present in this article do not have 
a fine moduli space. This leads to the question. 
 
\proclaim{Question} Does there exist a pair $(M,\Cal F^{diff})$ 
which admits a (non-empty) fine moduli space? 
\endproclaim 
 
Here are two examples. 
 
\example{Example} Consider the classical Reeb foliation $\Cal 
F^{Reeb}$ of $\Bbb S^3$ (see \cite{Go}). Fix a Riemannian metric 
$\mu$ on $\Bbb S^3$ and an orientation on $\Cal F^{Reeb}$. The 
restriction of $\mu$ to a leaf $L$ defines a Riemannian metric on 
the oriented manifold $L$, that is a structure of Riemann surface on 
$L$. The foliation becomes in this way a foliation by Riemann 
surfaces. In particular $\Cal C(\Bbb S^3,\Cal F^{Reeb})$ is not 
empty and there is a map between the set of conformal classes of 
Riemannian metrics on $\Bbb S^3$ and $\Cal C(\Bbb S^3,\Cal 
F^{Reeb})$. Conversely, fix a foliation by complex manifolds on 
$(\Bbb S^3,\Cal F^{Reeb})$. The integrable almost complex operator 
defines a conformal class of Riemannian metrics on each leaf and the 
property of transverse smoothness of the operator means that these 
conformal classes fit into a conformal class of Riemannian metrics 
on the whole $\Bbb S^3$. Therefore this map is surjective. 
 
Notice that the choice of the orientation of $\Cal F^{Reeb}$ is not 
important. Indeed, let $J$ be an integrable almost CR-structure on $\Cal 
F^{Reeb}$ respecting a fixed orientation. Then $-J$ defines an 
integrable almost CR-structure on $\Cal F^{Reeb}$ respecting the other 
orientation. So the situation here is different from that of $\Bbb S^{5}$.
 
The Reeb foliation is constituted by one compact leaf diffeomorphic 
to a real $2$-torus and by non compact leaves, which are all 
diffeomorphic to $\Bbb R^2$. More precisely, the sphere $\Bbb S^3$ 
is decomposed into the union of two solid tori 
$$ 
\Bbb S^3=\Bbb S^1\times\overline{\Bbb D} \cup \overline{\Bbb D}\times\Bbb S^1 
$$ 
each copy of these being endowed with a foliation (Reeb component) 
which can be described as follows (cf \cite{M-V, Lemma 2}). Let 
$$ 
X=\R^2 \times [0,\infty)\setminus \{(0,0,0)\} 
$$ 
foliated by the level sets of the projection onto the second factor. 
Let 
$$ 
h\ :\ (x,y,t)\in X \longmapsto (\alpha\cdot x,\beta\cdot y, d(t))\in 
X 
$$ 
where $0<\alpha^2+\beta^2<1$ and where $d$ is a smooth function of 
$\R$ into $\R^+$, which is a contracting diffeomorphism of $\Bbb 
R^+$ and which is the identity on $\R^-$. Then the quotient of $X$ by the 
group generated by $h$ is a Reeb component. 
 
Assume now that $(\Bbb S^3,\Cal F^{Reeb})$ is endowed with a 
foliation by Riemann surfaces $\Cal F$ and consider the induced 
integrable, Levi flat, almost CR-structure on each Reeb component. Then the 
previous covering $X$ becomes, by pull-back by the covering map, a 
product foliated covering. The compact leaf of $\Cal F$ is an 
elliptic curve $\Bbb E_{\tau}$ and thus the boundary leaf of $X$ is 
$\Bbb C^*$. We claim that the non-compact complex leaves of $\Cal F$ 
are all biholomorphic to $\Bbb C$. Indeed, by Corollary 3, we may assume that the previous structure is tame. Then the claim
is proved in the first Example of Section 4.


\remark{Remark}
Here is an alternative proof for this fact.
Assume that there exists 
$t>0$ such that $L_t$ is not biholomorphic to $\C$. By 
Riemann-Poincar\'e-Koebe Uniformization Theorem, it is thus 
biholomorphic to the open unit disk $\Bbb D$. In particular, there 
does not exist a sequence of holomorphic functions of the disk of 
radius $n$ into $L_t$ with derivatives in $0$ bounded by below. Now, 
take an increasing sequence $\Cal D=(D_n)$ of disks of radius $n$ in 
the boundary leaf of $X$ seen as $\Bbb C^*$. Note that $L_{d^p(t)}$ 
is biholomorphic to the unit disk for all $p$. As $d$ is 
contracting, this means that there exists a copy of $L_t$ as close 
to the boundary leaf as wanted. Note that the map 
$$ 
i_s \ :\ (x,y,0)\in (\R^2\times\{0\}\setminus 
\{(0,0,0)\})\longmapsto (x,y,s)\in L_s 
$$ 
injects smoothly the boundary leaf in the leaf $L_s$ for all $s\geq 
0$. Using this map, we may embed {\it smoothly} the family $\Cal D$ 
into $L_{d^p(t)}$ for all $p$ by embeddings with derivatives in $0$ 
bounded by below. Although $i_s$ is not a priori holomorphic, as 
$i_0$ is the identity map of $\Bbb C^*$, then for $s$ sufficiently 
small, $i_s$ is quasi-conformal with a distortion factor which tends 
to one as $s$ tends to $0$. As a consequence, the family $\Cal D$ in 
$L_{d^p(t)}$ is a family of quasi-conformal disks with a distortion 
factor which tends to one when $p$ tends to $\infty$. Passing to the 
limit in $p$, we obtain holomorphic embeddings of each disk of $\Cal 
D$ in $L_t\simeq\Bbb D$ with derivatives in $0$ bounded by below. 
Contradiction. 
\endremark
\medskip

 Moreover, we know from the second Example of Section 4 that the interior part of
 the product foliated covering $X$ is biholomorphic to $\Bbb C\times (0,\infty)$.

We want now to prove that $X$ is CR-isomorphic to
$\Bbb C\times\Bbb R^+\setminus\{(0,0)\}$. Let $A$ be an (open) annulus in $\partial X\simeq\Bbb C^*$. Observe that there exists a relatively compact open set $\Cal A$, a real number $\epsilon$ and a map
$\phi : \Cal A\to A\times [0,\epsilon)$ such that

\medskip
\noindent (i) The intersection of $\Cal A$ with the leaf $X_t$ is a topological annulus $A_t$ (with $A_0=A$) for $t<\epsilon$ and is empty for $t\geq \epsilon$.

\noindent (ii) The map $\phi$ is a CR-isomorphism.
\medskip

This can be done by hand in this particular case, or may be deduced directly from \cite{A-V, Proposition 2}. Let $D$ be the disk of $\Bbb C$ whose boundary is the exterior boundary of $A$. In the same way, let $D_t$ be the topological
disk of $X_t$ whose boundary is the exterior boundary of $A_t$. Call $\Cal K$ the union of $D\setminus \{0\}$ and $D_t$ for all $t$ in $(0,\epsilon)$. 
Note that the set $\Cal K$ is open in $X$ and contains $\Cal A$. Recall that we may uniformly compactify the family $X$ by adding one point at infinity to each leaf (cf Proposition 5). This gives a family $X^c$ whose interior leaves are all 
biholomorphic to $\Bbb P^1$ and whose boundary leaf is $\Bbb C$.

We perform now the following surgery on $X^c$: cut
$\Cal K\setminus\Cal A$ and glue $D\times [0,\epsilon]$ along $\Cal A$ by $\phi$. We obtain in this way a new CR-manifold $\hat X^c$ such that

\medskip

\noindent (i) There is a CR-injection from $X^c\setminus \Cal K$ in $\hat X^c$.

\noindent (ii) The family $\hat X^c$ is a deformation family of $\Bbb P ^1$ parametrized by $[0,\infty)$.
\medskip

By Proposition 2, $\hat X^c$ is a trivial family. Hence there exists a CR-isomorphism $\psi$ between $X^c\setminus(\Cal K\setminus\Cal A)$ and $\Bbb P^1\times [0,\infty)\setminus ((D\setminus A)\times [0,\epsilon))$. We may transfer the action through $\psi$.   
To be more precise, if, in a model $X=(\Bbb R^2\times [0,\infty)\setminus \{(0,0)\}, J)$, the action defining the Reeb foliation is given by
$$
(z,t)\in\Bbb R^2\times [0,\infty)\setminus \{(0,0)\}\longmapsto (A_t(z),d(t))\in\Bbb R^2\times [0,\infty)\setminus \{(0,0)\}
$$
as usual, then let 
$$
B_t=\psi_{d(t)}\circ A_t\circ \psi_t^{-1}
$$
Through the conjugation by $\psi$, the action on $\Bbb P^1\times [0,\infty)\setminus ((D\setminus A)\times [0,\epsilon))$ 
is given by $(z,t)\mapsto (B_t(z),d(t))$ (where it is well-defined). Now, we may also see $(B_t)$ as a 1-parameter family
of biholomorphisms defined in a neighborhood of $\infty$ in $\Bbb P^1$. Observe that these biholomorphisms fix $\infty$ and that this fixed point is contracting. Hence by Koenigs' Theorem with parameters \cite{Mi3, p.74-75}, it may be
linearized and the linearization map may be assumed smooth in the parameter. Hence, there exists a neighborhood $\Cal U$ of the infinite section $X^c\setminus X$ in $X^c$ and a CR-isomorphism $\chi$ sending $\Cal U$
to an open neighborhood of $\infty\times [0,\infty)$ in $\Bbb P^1\times [0,\infty)$ which conjugates the previous action to a linear action
$$
(z,t)\in\Cal U'\subset\chi(\Cal U)\longmapsto (\lambda(t)\cdot z, d(t))\in\Cal U
$$
where $(\lambda(t))_{t\in [0,\infty)}$ is a smooth family of complex numbers satisfying $0<\vert\lambda(t)\vert <1$, and $\Cal U'$ a well-chosen open set.

Since the action on $\Bbb P^1\times [0,\infty)\setminus ((D\setminus A)\times [0,\epsilon))$ extends to $\Bbb P^1\times [0,\infty)$, we may use it to extend the CR-isomorphism $\chi$ as a CR-isomorphism from 
$(X^{diff}\setminus \{0\}\times [0,\infty),J)$ to $\Bbb C^*\times [0,\infty)$ (we drop the infinite section; since the compactification is uniform and all the maps fix the section, this is not a problem), and finally
from $X$ to  $\Bbb C\times\Bbb R^+\setminus\{(0,0)\}$ by Riemann's Theorem.

\medskip

As a consequence of Corollary 2, two ``complex'' Reeb components are biholomorphic 
if and only if their compact leaves are biholomorphic (notice that 
the previously described CR-isomorphisms do not exchange the 
leaves). Therefore, to determine the set of almost CR-structures of $(\Bbb 
S^3,\Cal F^{Reeb})$, we just have to determine which complex numbers 
can appear as modulus of the compact leaf. 
 
We claim that any modulus $\tau$ can be obtained: to do this apply 
twice Lemma 2 of \cite{M-V}, once to obtain a complex Reeb component 
with boundary leaf biholomorphic to $\Bbb E_{\tau}$ and then to 
obtain a complex Reeb component with boundary leaf biholomorphic to 
$\Bbb E_{\tau^{-1}}$. Glue together these two components via the 
biholomorphism switching meridians and parallels \cite{M-V, 
Corollary of Lemma 1}. 
 
  Observe that, if we set $\Cal F_{\tau}=(\Cal F^{Reeb}, J)$, then the "complex conjugate", that is
 $(\Cal F^{Reeb}, -J)$, is CR-isomorphic to $\Cal F_{-\bar\tau}$. This is because the complex conjugation
 on the leaves of the product foliated covering $\Bbb C\times\Bbb R^+\setminus\{(0,0)\}$ descends to a biholomorphism
 between $\Cal F_{\tau}$ and $\Cal F_{-\bar\tau}$.
 
As a conclusion of all that preceeds, the set of complex 
struct\-ures of the Reeb foliation can be identified with $\Bbb 
H/\text{PSL}_2(\Bbb Z)\simeq\Bbb C$. 
 
 We claim that it is a coarse moduli space. To see this, take any 
deformation family $\pi : M\to X$ of $(\Bbb S^3,\Cal F^{Reeb})$ and 
consider the natural map $\alpha_{\pi} : X\to \Bbb C$. Restricting 
$\pi$ to the union of the compact leaves of each fiber, we obtain a 
complex analytic deformation family of complex tori $\Cal E\to X$ by 
Proposition 6. The map $\alpha_{\pi}$ can thus be seen as the 
modular function along the fibers of this family and is therefore 
holomorphic in the base space $X$. In other words, $(\Bbb S^3,\Cal 
F^{Reeb})$ admits a coarse moduli space isomorphic to $\C$. 
 
We also claim that there does not exist a fine moduli space. This is 
because the translations of the compact leaf extend to 
CR-isomorphisms of the foliation, \cite{M-V, Lemma 2}. As a 
consequence, every holomorphic principal elliptic fiber bundle gives 
rise to a holomorphic family of deformations of $(\Bbb S^3,\Cal 
F^{Reeb})$. Indeed, let $X\to B$ be such a bundle and let $\Bbb 
E_\tau$ be its fiber. Let $(U_\alpha)_{\alpha\in A}$ be an open 
covering of $B$ and 
$$ 
g_{\alpha\beta} \: \ U_\alpha\cap U_\beta \longrightarrow \Bbb 
E_\tau 
$$ 
be a cocycle representing $X$. Then seeing now $\Bbb E_\tau$ not as 
the group of translations of the fiber $\Bbb E_\tau$ but as a 
subgroup of the group of CR-isomorphisms of $(\Bbb S^3,\Cal 
F^{Reeb})$ endowed with the complex structure $\Cal F_\tau$, 
we may 
construct from the cocycle $(g_{\alpha\beta})$ a holomorphic family 
of deformations of $(\Bbb S^3,\Cal F^{Reeb})$ whose fibers are all 
isomorphic to $\Cal F_\tau$. 
 
Assume that there exists a fine moduli space for $(\Bbb S^3,\Cal 
F^{Reeb})$. As observed before, this implies that the previous 
family of deformations is CR-trivial. Then the bundle $X\to B$ has 
to be holomorphically trivial. But there exist principal holomorphic 
elliptic bundles which are not trivial, for example the surface 
$\Cal S$ of this article. Contradiction.

\endexample 
 
\example{Example} Let $\Bbb T^2$ denote the real $2$-torus and let 
$X=\Bbb T^2\times\Bbb S^1$. Consider on $X$ the smooth foliation 
$\Cal F^{diff}$ by $\Bbb T^2$ given by the level sets of the 
projection onto the $\Bbb S^1$-factor. We claim that the set of 
integrable almost CR-structures of $(X,\Cal F^{diff})$ is an 
infinite-dimensional Fr\'echet space. 
 
To see that, recall the following construction of a complete 
deformation family for $\Bbb T^2$ (see \cite{M-K, p.18}). Let $\Bbb 
H$ denote the Poincar\'e half-plane and consider the group of 
transformations of $\Bbb H\times\Bbb C$ given by 
$$ 
G=\{(\tau, z)\in\Bbb H\times\Bbb C\longmapsto (\tau, z+m+n\tau)\in 
\Bbb H\times\Bbb C \quad\vert\quad (m,n)\in\Bbb Z^2\}\ . 
$$ 
This group acts freely and properly discontinuously on $\Bbb 
H\times\Bbb C$ so the quotient is a well-defined complex manifold 
$\Cal M$. It can be checked that the natural projection $\Cal M\to 
\Bbb H$ is a holomorphic submersion. This implies that $\Cal M$ is 
diffeomorphic to $\Bbb T^2\times\Bbb R^2$. The fiber over $\tau$ is 
the elliptic curve of modulus $\tau$. 
 
Let $\alpha$ be a smooth map from $\Bbb S^1$ to $\Bbb H$. We can 
take the pull-back of the submersion $\Cal M\to\Bbb H$ by $\alpha$. 
As this submersion is diffeomorphically trivial, the total space of 
the pull-back is diffeomorphic to $X$. 
$$ 
\CD 
X &\aro >>> &\Cal M \\ 
\aro VVV &&\aro VVV \\ 
\Bbb S^1 &\aro >\alpha>> &\Bbb H 
\endCD 
$$ 
This pull-back construction endows $(X,\Cal F^{diff})$ with a 
compatible foliation by elliptic curves, let us denote it by $\Cal 
F_{\alpha}$. 
 
Recall now that $\Bbb H/\text{PSL}_2(\Bbb Z)\simeq \Bbb C$ 
parametrizes the modulus of the elliptic curves. Let $\pi$ be the 
natural projection from $\Bbb H$ to $\Bbb H/\text{PSL}_2(\Bbb Z)$. 
Let $V=\Bbb T^2\times [0,1]$ and let $\Cal F$ be a foliation by 
elliptic curves of $V$ compatible with the natural projection onto 
$[0,1]$. Then, there is a smooth map $\bar f$ from $[0,1]$ to $\Bbb 
H/\text{PSL}_2(\Bbb Z)\simeq \Bbb C$ given by the modular function 
on each fiber. And $\bar f$ lifts to a smooth map $f$ from $[0,1]$ 
to $\Bbb H$. Consider the pull-back foliation by elliptic curves on 
$V$ defined as $f^*\Cal M$. 
 
\proclaim{Lemma} The foliations $\Cal F$ and $f^*\Cal M$ are 
biholomorphic. 
\endproclaim 
 
\demo{Proof} Since the family $\Cal M\to\Bbb H$ is {\it complete}, the foliations 
$f^*\Cal M$ and $\Cal F$ are locally biholomorphic. Hence $f^{*}\Cal M$ is locally biholomorphic to the quotient
of $\Bbb C\times [0,1]$ by a smooth family of lattices $\Cal L_{t}$. This implies that the universal covering of $\Cal F$ is a locally trivial
CR bundle with fiber $\Bbb C$ over the unit interval. But such a bundle is globally trivial and the previous biholomorphism is also global.
$\square$
\enddemo
 
 
As a consequence of this Lemma, every foliation of $(X,\Cal 
F^{diff})$ is biholomorphic to some $\Cal F_{\alpha}$. On the other 
hand, take two maps $\alpha_1$ and $\alpha_2$ from $\Bbb S^1$ to 
$\Bbb H$. If $\Cal F_{\alpha_1}$ is biholomorphic to $\Cal 
F_{\alpha_2}$, we have 
$$ 
\pi\circ\alpha_2\equiv\pi\circ\alpha_1 
$$ 
Conversely, if this equality is satisfied, it lifts to a 
diffeomorphism $f$ of $X$ which maps the pull-back foliation $\Cal 
F_{\alpha_1}$ to the pull-back foliation $\Cal F_{\alpha_2}$ without 
exchanging the leaves. But, by the previous equality, the modulus of 
the leaf of $\Cal F_{\alpha_1}$ above $z\in\Bbb S^1$ and the modulus 
of the leaf of $\Cal F_{\alpha_2}$ above $z$ are the same and thus 
$f$ is a  biholomorphism.

Therefore $\Cal C(X,\Cal F^{diff})$ can be identified with the space 
of smooth maps from $\Bbb S^1$ to $\Bbb H$, up to action of 
PSL$_2(\Bbb Z)$ on the target space. This is an infinite-dimensional 
Fr\'echet space \cite{P-S}. Note that in this case there is no 
natural complex structure on this set of integrable almost CR-structures, 
and no coarse moduli space. 
 
\remark {Remark} The space of smooth maps from $\Bbb S^1$ to $\Bbb 
H$, up to action of PSL$_2(\Bbb Z)$ on the target space is {\it 
different from} the loop space of $\Bbb H/\text{PSL}_2(\Bbb Z)\simeq 
\Bbb C$. In fact, a map from $\Bbb S^1$ to $\Bbb H/\text{PSL}_2(\Bbb 
Z)\simeq \Bbb C$ lifts to a map from $[0,1]$ to $\Bbb H$ such that 
the images of $0$ and $1$ belong to the same PSL$_2(\Bbb Z)$-orbit. 
By the pull-back construction, this defines a foliation by elliptic 
curves on the suspension of $\Bbb T^2$ by some element of 
PSL$_2(\Bbb Z)$, and not on $X$ (cf the Example after Proposition 3). 
\endremark 
\medskip 
 
Note also that if we define $\Cal C(X,\Cal F^{diff})$ as the set of 
foliations by complex manifolds of $X$ diffeomorphic to $\Cal 
F^{diff}$ {\it up to CR-isomorphisms}, then it can be identified 
with the space of smooth maps from $\Bbb S^1$ to $\Bbb H$, up to 
action of PSL$_2(\Bbb Z)$ on $\Bbb H$, and {\it up to 
reparametrization}. 
\endexample 
 
\head {\bf 7. The set of integrable almost CR-structures of $\Cal F_{\Bbb 
C}^{diff}$} 
\endhead 
 
Let $X=K\times\Bbb S^1$. By Lemma 1, every complex structure on $X$ 
is that of a primary Kodaira surface. It is obtained from a $\Bbb 
C^*$-bundle over an elliptic curve $\Bbb E_{\alpha}$ of Chern number 
$-3$ by taking the quotient of the fibers by a fixed homothety. 
Therefore, the set of pairwise non-biholomorphic complex structures 
on $X$ -let us denote it by $\Cal C(X)$- identifies naturally with 
$$ 
\{(\alpha, \beta,x) \quad\vert\quad \alpha\in\Bbb 
H/\text{PSL}_2(\Z), \ \beta\in\Bbb H/\text{PSL}_2(\Z), \ 
x\in\text{Pic}_{-3}(\Bbb E_{\alpha})\} 
$$ 
where $\text{Pic}_{-3}(\Bbb E_{\alpha})$ denotes the subset of the 
Picard group of $\Bbb E_{\alpha}$ constituted by elements 
corresponding to line bundles of Chern number $-3$. Recall that it 
has a natural structure of an elliptic curve \cite{Gu, \S 7-8}. 
 
\remark{Remark} The minus sign of the Chern number is not important. 
In fact, if $X\to\Bbb E_\alpha$ is an elliptic fibration of Chern 
number $-3$, the automorphism $z\to -z$ on the fibers sends $X$ 
biholomorphically onto an elliptic fibration of Chern number $3$. As 
a consequence, we could replace $\text{Pic}_{-3}(\Bbb E_{\alpha})$ 
by $\text{Pic}_{3}(\Bbb E_{\alpha})$ in the previous identification 
of $\Cal C(X)$. 
\endremark 
\medskip

We will denote the corresponding manifolds by $\Cal 
S(\alpha,\beta,x)$ and will say that such a complex $X$ is of type 
$(\alpha,\beta,x)$.  We denote by $W(\alpha,x)$ the corresponding 
$\Bbb C^*$-bundle. By abuse of notation, $\alpha$ (respectively 
$\beta$) will be considered as an element of $\Bbb H$ and not only 
as a class of $\Bbb H/$PSL$_2(\Bbb Z)$. Let $\alpha\in\Bbb H$. 
Consider an embedding $i_{\alpha}$ of $\Bbb E_\alpha$ into $\Bbb 
P^2$ as a cubic curve. To this embedding is associated the 
$\C^*$-bundle obtained by pull-back by $i_{\alpha}$ of the bundle 
$\C^3\setminus\{(0,0,0)\}\to \Bbb P^2$. This bundle is independent 
of the embedding. Indeed, two 
distinct embeddings are conjugated by an element of $\text{PGL}_{3}(\Bbb C)$, which induces an isomorphism of the tautological bundles. We call this 
bundle {\it the natural $\C^*$-bundle} of $\Bbb E_\alpha$ and we 
denote it by $\tilde \alpha$. 
 
Recall (cf Section 1.3) that $\Cal F_{\Bbb C}$ is obtained by gluing 
a foliation of an open set $\Cal N$ and a foliation of its 
complement. We call {\it interior part} this set $\Cal N$ with its 
foliation and denote by $\Cal S_1$ the compact leaf which bounds it. 
On the other hand, the foliation of $\Bbb S^5\setminus \Cal N$ 
contains another compact leaf. We denote it by $\Cal S_2$. The leaf 
$\Cal S_2$ separates $\Bbb S^5\setminus \Cal N$ into an open set 
foliated by leaves diffeomorphic to $W$ and an open set foliated by 
Milnor fibers. We call the first open set the {\it collar} and the 
second one {\it the exterior part}.

We turn now to the description of the set of integrable 
almost CR-structures of the foliation of [M-V]. 
 
\proclaim{Theorem B (Rigidity Theorem)} 
 
\noindent (i) Let $\Cal F$ be a foliation of $\Bbb S^5$ by complex 
surfaces diffeomorphic to $\Cal F_{\Bbb C}$. Then, the two compact 
leaves are of respective type $(\alpha,\beta,\tilde\alpha)$ and 
$(\alpha,\beta',\tilde\alpha)$, where $\alpha$, $\beta$ and $\beta'$ 
are any classes of $\Bbb H/\text{\rom{PSL}}_2(\Z)$. 
 
\noindent (ii) Let $\Cal F$ and $\Cal F'$ be two foliations of $\Bbb 
S^5$ by complex surfaces diffeomorphic to $\Cal F_{\Bbb C}$. Then 
$\Cal F$ and $\Cal F'$ are biholomorphic if and only if the compact 
leaf $\Cal S_1$ of $\Cal F$ is biholomorphic to the compact leaf 
$\Cal S'_1$ of $\Cal F'$ and the compact leaf $\Cal S_2$ of $\Cal F$ 
is biholomorphic to the compact leaf $\Cal S'_2$ of $\Cal F'$. 
\endproclaim 
 
The following Theorem is a realization result. 
 
\proclaim{Theorem C (Realization Theorem)} The set $\Cal C(\Bbb S^5, 
\Cal F_{\Bbb C}^{diff})$ identifies with 
$$ 
\Bbb H/\text{\rom{PSL}}_2(\Z)\times \Bbb H/\text{\rom{PSL}}_2(\Z) 
\times \Bbb H/\text{\rom{PSL}}_2(\Z)\simeq \C\times\C\times \C \ . 
$$ 
A point of $\Cal C(\Bbb S^5, \Cal F_{\Bbb C}^{diff})$ is entirely 
determined by the images by the modular function of the moduli of 
the following elliptic curves: 
 
\noindent (i) the common base of the two compact leaves. 
 
\noindent (ii) the fiber of $\Cal S_1$. 
 
\noindent (iii) the fiber of $\Cal S_2$. 
\endproclaim 
 
In other words, the deformations of $\Cal F_{\Bbb C}$ are fixed by 
the deformations of the two compact leaves, and only three 
``parameters'' among the six of $\Cal C(X)\times \Cal C(X)$ may 
vary. 
 
\remark{Remark} As noticed at the end of Section 1.3, the foliation 
of $\Cal F_{\Bbb C}^{diff}$ is not unique up to smooth isomorphisms, 
but only up to topological isomorphisms. As a consequence, it is 
easy to check that the set $\Cal C(\Bbb S^5, \Cal F_{\Bbb C}^{top})$ 
is not finite dimensional. 
\endremark 
 
\demo{Proof of Theorem B} Let $\Cal F$ be a foliation of $\Bbb S^5$ 
by complex surfaces diffeomorphic to $\Cal F_{\Bbb C}$. 
\medskip 
\noindent $\underline{\text{Step 1 : the interior part}}$. We will 
prove that the non-compact leaves of the interior part are all 
biholomorphic and that their biholomorphism type is fixed by the 
type of $\Cal S_1$. Let us denote by $(\alpha, \beta, x)$ the type 
of $\Cal S_1$. 
 
The general argument, which will be used here and in steps 2 and 3, 
is the following. Take a non-compact leaf $L$. As it spirals over 
$\Cal S_1$, it (or an open set of it) admits a 
$C^{\infty}$-submersion over $\Cal S_1$ (see Section 4). So we may 
define $L^{pb}$ and compare it to $L$. {\it A priori} they are 
different as complex manifolds. Now, we determine the exact 
biholomorphism type of $L^{pb}$ (notice that it depends only on the 
complex structure of $\Cal S_1$) and show that it admits a 
holomorphic compactification. The compactification Lemma ensures us 
that $L$ admits also a holomorphic compactification. Then, making 
use of the Enriques-Kodaira classification, it is possible to give 
the biholomorphism type of this compactification and to prove that 
it is the same as the compactification of $L^{pb}$. Going back to 
$L$, we conclude that $L$ and $L^{pb}$ are biholomorphic and 
therefore that the complex structure of the non-compact leaves is 
fixed by the structure of the compact leaf.

Lemmas 1 and 2 are still valid. In particular, a non-compact leaf $L$ 
is diffeomorphic to a line bundle over $\Bbb E_{\alpha}$. Let 
$L^*\subset L$ be the complex manifold diffeomorphic to the 
associated $\C^*$-bundle. Let also $(L^{*})^{pb}$ denote the associated pull-back complex
structure. It follows from what preceeds and from the 
construction that, topologically, $L^{*}$ is the $\Bbb Z$-covering of 
$\Cal S_1$ obtained by unrolling the elliptic fibers so that 
$$ 
(L^*)^{pb}=W_{(\alpha,x)}\ . 
$$ 
 
We may partially compactify $(L^*)^{pb}$ as the line bundle 
associated to $W_{(\alpha,x)}$ by adding an elliptic curve $\Bbb 
E_{\alpha}$. We denote this line bundle by $L_{(\alpha,x)}$. It is 
diffeomorphic to $L$. The compactification Lemma shows now that $L$ 
can be compactified by adding a copy of $\Bbb E_{\alpha}$ and is 
diffeomorphic to a $\Bbb P^1$-bundle over $\Bbb E_{\alpha}$. As 
usual, let $L^c$ be the corresponding compact surface. 
 
\proclaim{Lemma 5} The compact surface $L^c$ is a ruled surface of 
genus one. 
\endproclaim 
 
\demo{Proof} From the exact sequence in homotopy of the smooth $\Bbb 
S^2$-fibration $L^c\to\Bbb S^1\times\Bbb S^1$, we know immediately 
that the second homotopy group of $L^c$ is generated by the fiber 
and that the fundamental group of $L^c$ is $\Bbb Z^2$. Therefore 
$L^c$ is minimal with $b_1$ equal to two. Recall that the group of orientation preserving diffeomorphisms of the 2-sphere 
retracts onto $\text{SO}_{2}$ \cite{Sm}; hence we may assume that our fiber bundle has structural group $\text{SO}_{2}$ and that it is
the boundary of a disk bundle. Therefore $L^{c}$ is cobordant 
to zero and its Chern numbers and signature are zero. The 
Enriques-Kodaira classification proves thus that $L^c$ is a ruled 
surface of genus one, or a minimal properly elliptic surface or an 
hyperelliptic surface. This last case is impossible since the 
universal covering of an hyperelliptic surface is contractible 
whereas the universal covering of $L^c$ has non-zero second homotopy 
group. 
 
By \cite{B-H-P-V, Chapter IV, Theorem 2.7}, we have $h^{1,0}(L^c)=1$, i.e. 
there exists a global holomorphic $1$-form on $L^c$. This implies 
that the Albanese torus is an elliptic curve $\Bbb E_{\tau}$ and the 
Albanese map is a surjective holomorphic map $\pi : L^c\to \Bbb 
E_{\tau}$. 
 
Assume now that $L^c$ is elliptic, that is that there exists a 
holomorphic map $p$ from $L^c$ to some Riemann surface $\Sigma$ 
whose generic fiber is an elliptic curve. Observe that the smooth vector fields of the base of the smooth fibration $L^c\to\Bbb S^1\times\Bbb S^1$ can be lifted to
smooth vector fields of $L^{c}$ by means of a connection. Hence the Euler 
characteristic of $L^c$ is zero, so this elliptic fibration is obtained 
from a torus bundle over $\Sigma$ by performing logarithmic 
transformations along some fibers (see \cite{F-M, Chapter II, 
Proposition 7.2}). In particular, the fibers of $p$ are non-singular 
and form a family $\Cal F$ of smooth elliptic curves which covers 
$L^c$. 
 
Let $C\in\Cal F$. The restriction of $\pi$ to $C$ is a holomorphic 
map between complex tori, that is a constant or a unramified 
covering of degree $d>0$. By continuity of the degree of $\pi$ in 
the family $\Cal F$, this does not depend on the curve $C$, that is 
one the following statements is verified. 
\medskip 
\noindent (i) For every $C\in\Cal F$, the map $\pi$ restricted to 
$C$ is constant. 
 
\noindent (ii) For every $C\in\Cal F$, the map $\pi$ restricted to 
$C$ is an unramified covering of fixed degree $d$. 
\medskip 
Assume (i). Then the fibers of $\pi$ contain elliptic curves, that 
is the generic fibers of $\pi$ are elliptic curves. In other words, 
$\pi$ is also an elliptic fibration, with base a complex torus. But 
this would imply that the second homotopy group of $L^c$ is zero and 
would give, as before, a contradiction. Indeed, let $i$ be a map 
from $\Bbb S^2$ to $L^c$. The projection $\pi\circ i$ is homotopic 
to a constant since the target space is a torus. Now, we may define 
a connection on the elliptic fibration by taking a transverse field 
to the foliation it defines. This allows us to lift the previous 
homotopy to a homotopy between $i$ and a map whose image lies in a 
fiber. Since the fiber is a torus we may then homotope $i$ to a 
constant. The statement (ii) holds. 
 
\remark{Remark} The same proof shows that, for any elliptic Seifert 
fibration $X\to B$ (that is with only multiple fibers), the $k$-th 
homotopy groups of $X$ and $B$ are equal (for $k>1$). 
\endremark 
\medskip 
 
As a consequence of (ii), the map $\pi$ is a holomorphic submersion, 
so, by Ehresmann's Lemma, a locally trivial smooth fiber bundle. The 
exact sequence in homotopy of the bundle implies that the fiber is 
$\Bbb P^1$. The Theorem of Fischer-Grauert \cite{F-G} implies that 
$\pi$ is a locally trivial holomorphic fiber bundle over $\Bbb 
E_{\tau}$ with fiber $\Bbb P^1$, i.e. a ruled surface of genus one. 
$\square$ 
\enddemo 
 
A non-compact leaf $L$ is obtained from the ruled surface $L^c$ by 
removing an elliptic curve $\Bbb E_{\alpha}$. This elliptic curve is 
a section of the $\Bbb P^1$-bundle since it is a section for its 
topological model and a holomorphic submanifold of $L^c$. In other 
words the ruled surface $L^c$ is a P-bundle over $\Bbb E_{\alpha}$ 
with one holomorphic section, P-bundle meaning locally trivial 
holomorphic bundle with fiber $\Bbb P^1$ and structural group 
PSL$_2(\C)$; whereas $L$ is a $A$-bundle over $\Bbb E_{\alpha}$, 
that is a locally trivial holomorphic bundle with fiber $\C$ and 
structural group $Aff(\C)$, the one-dimensional affine group. We 
want to prove that $L$ is a {\it line bundle} over $\Bbb 
E_{\alpha}$. This is equivalent to showing that the P-bundle $L^c$ 
admits another holomorphic section disjoint from the first one. 
Notice that $L$ is diffeomorphic to a line bundle of Chern number 
$-3$ over $\Bbb E_{\alpha}$. Consider an open covering $U_{\alpha}$ 
of $\Bbb E_{\alpha}$ and a cocycle 
$$ 
g_{\alpha\beta} \ : \ z\in U_{\alpha}\cap U_{\beta} \longmapsto 
(w\to a(z)w+b(z))\in Aff(\C) 
$$ 
representing the affine bundle $L$. Then this cocycle is smoothly 
homotopic to the cocycle 
$$ 
z\in U_{\alpha}\cap U_{\beta} \longmapsto a(z)\in\C^* 
$$ 
that is, the affine bundle is equivalent, as a smooth bundle, to the 
line bundle defined by the previous cocycle. Call this bundle the 
{\it associated line bundle}. Notice that this line bundle has also 
Chern number $-3$. 
 
The following result is a weak version of a Theorem of Atiyah. 
 
\proclaim{Theorem (see [At], Theorem 6.1)} Let $B$ be a A-bundle 
over an elliptic curve. If the degree of the associated line bundle 
is different from $0$ and from $-1$, then $B$ is projectively 
equivalent to a $\Bbb C^*$-bundle. 
\endproclaim 
 
Thanks to this Theorem, a non-compact leaf $L$ is a line bundle 
$L(\alpha, y)$. Nevertheless, step 1 is not yet finished, since {\it 
it is not a priori clear that $y$ is independent of $L$ and that $y$ 
is equal to $x$}. To prove this fact, we make use of the Dumping Lemma. 
 
\proclaim{Lemma 6} Every non-compact leaf is biholomorphic to 
$L_{(\alpha,x)}$. 
\endproclaim 
 
\demo{Proof of Lemma 6} The argument used here will be referred 
hereafter as the {\it dumping trick}. It uses of course the Dumping Lemma. Let 
$$ 
\tilde X^{diff}=\Bbb R^2\setminus \{(0,0)\}\times (\Bbb R^2\times 
[0,\infty))\setminus \{(0,0,0)\} 
$$ 
and define 
$$ 
\eqalign{ T(x,y,u,v,t)&=(x,y,(u,v)\cdot A,d(t))\cr 
S(x,y,u,v,t)&=((x,y)\cdot B, (u,v)\cdot (\Psi(x,y))^{-3},t) } 
$$ 
for $A$, $B$ and $\Psi$ matrices chosen so that the previous model 
forms exactly the smooth version of the complex $\Bbb Z^2$-covering 
used in \cite{M-V, p.922} to foliate the interior part $\Cal N$ 
(see Section 1.3). Notice that Int $Y=\text{Int }(\tilde X^{diff})/\langle S\rangle$ 
is diffeomorphic to $L^{diff}\times (0,\infty )$ and that $\partial Y$ is diffeomorphic to $W^{diff}$. Let $J(t)$ be the 
CR-structure on $\tilde X^{diff}$ coming from the foliation of $\Cal N$. The leaves are the level sets of 
the projection onto the $t$-factor. Let $J_0$ be 
the CR-structure induced on $Y$. Let $L_t$ be the complex leaf of 
$Y$ corresponding to the level set $t$. Remark that, for $t>0$, the 
leaf $L_t$ is a line bundle $L_{(\alpha, x(t))}$, whereas $L_0$ is 
$W_{(\alpha,x)}$. Indeed, the $\langle T\rangle$-action onto $Y$ 
sends the leaf $L_t$ onto the leaf $L_{d(t)}$. As a consequence, 
there is no $\langle T\rangle$-action on a leaf $L_t$ for $t>0$, 
thus the leaves of the interior of $\Cal N$ and the leaves of the 
interior of $Y$ are biholomorphic. Finally, the $\langle 
T\rangle$-action onto $L_0$ defines exactly the covering 
$W_{(\alpha,x)}\to \Cal S_1$.

The uniform compactification Lemma tells us that we may compactify Int $Y$ as a deformation family of $\Bbb P^1$-bundles by adding
$\Bbb E_{\alpha}\times (0,\infty)$ smoothly in $t$. Moreover, by Proposition 5, this uniform compactification can be extended to the boundary leaf.
Consider now the union of the normal bundles of this smooth family of elliptic curves.
We obtain a smooth deformation family of line bundles (in the sense of Section 3) over $\Bbb E_{\alpha}$ of fixed topological degree. 
Now this family satisfies the hypotheses of the Dumping Lemma, hence all these bundles are isomorphic. Since these bundles are indeed isomorphic to the bundle 
$L_t$, and since $L_0=L_{(\alpha,x)}$, we have that every $L_t$ is biholomorphic to $L_{(\alpha,x)}$.
$\square$ 
\enddemo


\medskip 
\noindent $\underline{\text{Step 2 : the collar}}$. The argument is 
very similar to that of the interior part. A non-compact leaf $L$ is 
diffeomorphic to $W$, that is to a $\C^*$-bundle over an elliptic 
curve. Such a leaf has two ends, corresponding to the two ends of 
$\C^*$, its fiber. It spirals over $\Cal S_1$ at one end and over 
$\Cal S_2$ at the other end. Therefore $L^{pb}_1$ (respectively $L^{pb}_2$)
is the $\Z$-covering of $\Cal S_1$ (respectively of $\Cal S_2$) obtained by 
unrolling the fibers. Denote by $(\alpha',\beta',x')$ the type 
of $\Cal S_2$. 
 
As in step 1, the use of the compactification Lemma on both ends ensures us that 
$L$ admits a holomorphic compactification by adding two elliptic 
curves $\Bbb E_{\alpha}$ and $\Bbb E_{\alpha'}$. This compactification $L^c$ is {\it 
diffeomorphic} to a ruled surface of genus one. By Lemma 5, it {\it 
is} a ruled surface of genus one and the two elliptic curves are holomorphic sections of this bundle, hence $\alpha=\alpha'$. Besides, $L$ is biholomorphic to 
some $W(\tilde\alpha,\tilde x)$. Notice that, in this case, we do 
not need the result of Atiyah, since we know that $L^c$ is a ruled 
surface with two disjoint holomorphic sections. 
 
Finally, we use the dumping trick of Lemma 6 for the product 
foliated covering described in \cite{M-V, p.925--926} to prove that 
each non-compact leaf is in fact biholomorphic to $W(\alpha,x)$. The 
smooth model of this covering can be described as follows. Let 
$$ 
Z=W^{diff}\times [-1,0]\subset \Bbb R^6 \times [-1,0] 
$$ 
and let 
$$ 
g(x,t)=((x_1,x_2)\cdot A, (x_3,x_4)\cdot A, (x_5,x_6)\cdot A, h(t)) 
$$ 
where $A$ is a fixed matrix and where $h$ is a smooth diffeomorphism 
from $\Bbb R$ to $\Bbb R$ described in \cite{M-V, p.925} whose only 
fixed points in $[-1,0]$ are $-1$ and $0$. 
 
Set $Z_0=W^{diff}\times [-1,0)$ and $Z_1=W^{diff}\times (-1,0]$.
Then $Z_0\to Z_0/ \langle g\rangle$ and $Z_1\to Z_1/ \langle g\rangle$ are product foliated coverings of 
the collar minus $\Cal S_1$ or minus $\Cal S_2$ to which we apply the dumping trick. As a consequence, $L=W_{(\alpha,x)}=W_{(\alpha,x')}$ and $x=x'$. The two 
compact leaves differ only by the modulus of the fibers. 
 
\medskip 
\noindent $\underline{\text{Step 3 : the exterior part}}$. We 
consider the product foliated covering of the exterior part coming 
from the Milnor fibration. More precisely, setting 
$$ 
Y=P^{-1}(\Bbb R^+)\setminus \{(0,0,0\} \subset \Bbb R^6 
$$ 
and considering $Y$ as a smooth manifold with boundary, we have a 
$\Z$-covering 
$$ 
s \ :\  Y \longrightarrow \Bbb S^5\setminus \text{Int } (\Cal N) \ . 
$$ 
satisfying the following commutative diagram 
$$ 
\CD 
Y \aro >s>> \Bbb S^5\setminus \text{Int }(\Cal N) \\ 
\aro V p VV  \aro VVV \\ 
\Bbb R^+ \aro > \Bbb Z-\text{covering}>> \Bbb S^1 
\endCD 
$$ 
 
The deck transformation group is generated by 
$$ 
x\in Y\longmapsto ((x_1,x_2)\cdot A, (x_3,x_4)\cdot A, 
(x_5,x_6)\cdot A) 
$$ 
where $A$ is the same fixed matrix as in step 2 (to be more precise, it is the real matrix form of the complex number $\lambda\omega$, see \cite{M-V, p.924}). 

\remark{Remark}
It is important to keep in mind that this product foliated covering is not tame.
\endremark
 \medskip

We may thus equip $Y$ with the pull-back CR-structure induced from 
this covering. The associated smooth foliation is given by the level 
sets of $P$ (cf \cite{M-V, p.923--924}). 
 
Let $L_t$ be the complex leaf diffeomorphic to $P^{-1}(\{t\})$. 
Assume $t\not = 0$. Then $s(L_t)$ is still biholomorphic to $L_t$ 
and descends as a leaf in $\Bbb S^5\setminus\Cal N$ whose unique end 
spirals over $\Cal S_2$. There exists thus a closed set $F_t$ in 
$L_t$ such that $L_t\setminus F_t$ is diffeomorphic to a covering of 
$\Cal S_2\simeq K\times \Bbb S^1$. Indeed, by \cite{Mi1, Lemma 6.1 
and Theorem 5.11}, it is diffeomorphic to $K\times\Bbb R$. This 
implies that $(L_t\setminus F_t)^{pb}$ is biholomorphic to $W$. By 
use of the compactification Lemma, we may thus compactify $L_t$ by 
adding an elliptic curve $\Bbb E_{\alpha}$ and obtain a compact 
surface $L^c_t$. 

\remark{Remark}
The use of the Compactification Lemma is possible even if the covering is not tame because it satisfies the weaker condition
detailed in the remark after the Compactification Lemma.
\endremark
\medskip

On the other hand, $L_t$ is diffeomorphic to the 
affine cubic surface $P^{-1}(\{t\})$ of $\Bbb C^3$, so that $L^c_t$ 
is diffeomorphic to a cubic surface in $\Bbb P^3$. 


\proclaim{Lemma 7} The compact surface $L_t^c$ is biholomorphic to the blow-up of $\Bbb P^2$ at $6$ points or to the blow-up of
the Hirzebruch surface $\Bbb F_2$ at $5$ points.. 
\endproclaim 
 
\demo{Proof} Since $L_t^c$ is diffeomorphic to a cubic surface, it 
is diffeomorphic to $\Bbb P^2\sharp (6) \overline{\Bbb P^2}$, the 
blow-up of $\Bbb P^2$ at $6$ points. We infer from this description 
that $L_t^c$ has {\it at most} $6$ {\it disjoint} embedded rational 
curves with self-intersection $-1$. The minimal model of $L_t^c$, 
let us denote it by $X$, satisfies 
$$ 
c_1^2(X)\geq 3 \qquad c_2(X)\leq 9 \qquad \tau\geq -5 
$$ 
where $c_1^2$ and $c_2$ denote the Chern numbers of $X$ and $\tau$ 
the signature. 
 
The Enriques-Kodaira classification tells us that $X$ is a rational 
surface or that $X$ is of general type. This last case is 
impossible, due to the following deep result of Friedmann and Qin. 
 
\proclaim{Theorem [F-Q]} A surface of general type cannot be 
diffeomorphic to a rational surface. 
\endproclaim 
 
Therefore, $X$ is rational. There are two cases: 
\medskip 
\noindent (i) the surface $X$ is $\Bbb P^2$ and $L_t^c$ is the 
blow-up of $\Bbb P^2$ at $6$ points. 
 
\noindent (ii) the surface $X$ is a Hirzebruch surface $\Bbb F_a$ 
and $L_t^c$ is the blow-up of $\Bbb F_a$ at $5$ points. 
\medskip 
 
Assume (ii). Then, since $\Bbb F_1$ is not minimal and since the 
blow-up of $\Bbb P^2$ at six points and the blow-up of $\Bbb 
F_0=\Bbb P^1\times\Bbb P^1$ at five points are biholomorphic, we 
assume that $a\geq 2$. 
 
Assume that $a=2r$ is even. Then $\Bbb F_a$ is diffeomorphic to $\Bbb S^2\times \Bbb S^2$, hence a basis of topological cycles of dimension $2$ of $\Bbb F_a$ is given by $C_1$
diffeomorphic to $\{pt\}\times \Bbb S^2$ and $C_2$ diffeomorphic to $\Bbb S^2\times \{pt\}$. On the other hand, we may assume that $C_1$ is a fibre of the bundle $\Bbb F_a\to \Bbb P^1$ and that
a basis of analytical cycles on $\Bbb F_a$ is given by $C_1$ and by its zero section $\tilde S$ homologous to $C_2-rC_1$. Consider the blow-up map $p : L_t^c\to \Bbb F_a$. 
A basis of analytic cycles of dimension $2$ of $L_t^c$ is then given by $S$, the pull-back by $p$ 
of $\tilde S$, by $E_1$, ..., $E_5$, the exceptional divisors,
and by $F$ the pull-back by $p$ of a fiber of $\Bbb F_a$, not containing any of the blown up points. Let $k$ be the number of points to blow up which belongs to $p(S)$. Assume that $E_1$, ..., $E_k$
are the corresponding exceptional divisors. Notice that $k$ may take any value between $0$
and $5$. One obtains easily the following table of intersection numbers.
$$\eqalign{
S\cdot S&=-a-k, \quad S\cdot F=1, \quad S\cdot E_i=1 \ (i\leq k), \quad S\cdot E_i=0 \ (i>k) \cr
F\cdot S&= 1, \quad F\cdot F=0, \quad F\cdot E_i=0 \ (\text{for all }i) \cr
}
$$

We know that $L_t^c$ contains an elliptic curve $\Bbb E_\alpha$ such that the pair $(L_t^c, \Bbb E_\alpha)$ is diffeomorphic to the pair constituted by a non-singular cubic surface of $\Bbb P^3$ and by its section at infinity.
Using this model, we see that $L_t^c$ is diffeomorphic to the blow-up of $\Bbb P^2$ at six points. Then, another basis of topological cycles for $L_t^c$ is given by $H$ diffeomorphic to a line in $\Bbb P^2$ not intersecting the exceptional divisors and to $D_1$, ..., $D_6$ diffeomorphic to the $6$ exceptional divisors. In this basis, we also have that (see \cite{Ha, V.4.8})
$$
\Bbb E_\alpha=3H-D_1-...-D_6
$$
where $=$ stands here and in the sequel of the proof for ``is homologous''.

On the other hand, one checks easily that we may assume that
$$
D_1=C_1-E_1, \quad D_2=C_2-E_1, \quad H=C_1+C_2-E_1, \quad E_{i}=D_{i+1} \text{ (for } 2\leq i\leq 5)
$$
where we still call $C_1$ and $C_2$ in $L_t^c$ the pull-back by $p$ of $C_1$ and $C_2$ in $\Bbb F_a$.

Now, the elliptic curve $\Bbb E_\alpha$ must have non-negative intersection with the basis of analytic cycles of $L_t^c$. Taking into account that
$$
F=C_1 \qquad S=C_2-rC_1-E_1-\hdots -E_k
$$
this gives the following inequality
$$
\Bbb E_\alpha\cdot S=2-2r-k \geq 0
$$
Hence $r=0$ or $r=1$ and $k=0$. But $r=0$ means $a=0$, a case that we have already excluded.
So we have $r=1$ and $k=0$,
i.e. $L_t^c$ is the blow-up of $\Bbb F_2$ at five points not belonging to the zero-section $S$.

Assume now that $a=2r+1$ is odd. Exactly the same line of arguments leads to the inequation
$$
2-(2r+1)-k\geq 0
$$
hence $r=0$ and $k=0$. So $a=1$, a case that we have already excluded.
$\square$
\enddemo

 
We would like to thank Lucy Moser for the following observation.

\remark{Remark}
There really exists such an elliptic curve in the blow-up of $\Bbb F_2$ at five points, so we cannot exclude this case
at this stage. Indeed, consider the curve
$z^2y^2-x^4-z^4$ in $\Bbb P^3$. This is an elliptic curve with a unique singularity: a tacnode at infinity. The resolution of this singularity
requires two successive blow-ups. In this way, we obtain a non-singular elliptic curve in the blow-up of $\Bbb F_1$ at one point belonging to the zero section. But this is the same
as the blow-up of $\Bbb F_2$ at a point not belonging to the zero section this time. Observe that, in this last description, the elliptic curve does not intersect
the exceptional divisor. Hence we may blow down to $\Bbb F_2$ keeping the non-singular elliptic curve. Finally blow-up this curve at five points.
\endremark
\medskip

We want to conclude that $L^c_t$ is the blow-up of $\Bbb P^2$ at six points, and even more that it embeds in $\Bbb P^3$ as a non-singular cubic surface.
We first need to collect one more fact about this surface.

\proclaim{Lemma 8}
The surface $L_t^c$ admits an automorphism of order three.
\endproclaim

\demo{Proof}
Consider the product foliated covering Int $Y\to (0,\infty)$. By the uniform Compactification Lemma, the compactification of the leaves $L_t$ as $L_t^c$ can be assumed 
uniform. This gives a deformation family $\bar Y$ over $[0,\infty)$ whose interior has compact leaves. For each $t>0$, the corresponding leaf is $L_t^c$ and contains an elliptic
curve $E_t$ biholomorphic to $\Bbb E_\alpha$.

The deck transformation group identifies leaves of Int $Y$. In the presentation given at
the beginning of step 3, the leaves $L_t$ and $L_{f(t)}=L_{\lambda^{3}t}$ are CR-isomorphic and it extends as a CR-isomorphism between $L_t^c$ and $L_{f(t)}^c$.
On the other hand, Int $\bar Y$ is {\it diffeomorphic} to the suspension of $(L_t^c)^{diff}$ by a diffeomorphism of order three (this is the monodromy of the Milnor's fibration used in the initial construction of the foliation
$\Cal F_{\Bbb C}$). More precisely, on the smooth model, representing $L_t$ as the affine cubic surface $P^{-1}(t)$ for fixed $t$, the monodromy in restriction to $L_t$ is given by the map
$$
(z_1,z_2,z_3)\in P^{-1}(t)\longmapsto (\omega\cdot z_1, \omega\cdot z_2,\omega\cdot z_3) \in P^{-1}(t)
$$
which extends to the projectivization of $P^{-1}(t)$ as the identity on the elliptic curve.

We want to prove that this diffeomorphism of order three is in fact an automorphism of order three of each leaf. This could be easily deduced from Proposition 3 if we knew that all the leaves are biholomorphic. But here we
follow the inverse way: we want to use the existence of such an automorphism to deduce that all leaves of Int $\bar Y$ are biholomorphic. To achieve our goal, we first need a better model
for Int $\bar Y$: we will construct a CR-map from Int $\bar Y$
to $\Bbb P^3$ which is an embedding when restricted to an open and dense subset of any leaf $L_t^c$.

To avoid cumbersome repetitions in the sequel, we say that a leaf $L_t^c$ is of type (i) if it is the blow-up of $\Bbb P^2$ at $6$ points, and of type (ii) if it is the blow-up of $\Bbb F_2$ at five points.
We thus a priori have in our deformation family Int $\bar Y$ leaves of type (i) and leaves of type (ii).

From the computations of intersection numbers held in Lemma 7, it follows that in case (i), the six points belong to a cubic curve of $\Bbb P^2$, whereas, in case (ii), the five points belong 
the elliptic curve of $\Bbb F_2$ described in the previous remark. In both cases, it follows that it is in the linear system of the anticanonical divisor $-K_t$ of $L_t^c$.
This family of divisors gives rise to a family of line bundles over $L_t^c$ that we still call $-K_t$ or $E_t$.

We claim that, for all $t$, we have
$$
\dim H^0(L_t^c,\Cal O(E_t))=4\ .
$$
(we owe this computation to Laurent Bonavero).

Consider the short exact sequence \cite{B-H-P-V, p.62}
$$
0 \aro>>> \Cal O_{L_t^c} \aro >>> \Cal O_{L_t^c}(E_t) \aro>>> \Cal O_{E_t}(E_t) \aro>>> 0\ .
$$

Since $H^1(L_t^c,\Cal O_{L_t^c})$ is zero by 1-connectedness of $L_t^c$, and since $H^2(L_t^c,\Cal O_{L_t^c})$
is zero by Serre duality, we obtain from the long exact sequence
$$
H^1(L_t^c,\Cal O_{L_t^c}(E_t))=H^1(E_t,\Cal O(E_t))
$$
and by Serre duality on $E_t$,
$$
H^1(L_t^c,\Cal O_{L_t^c}(E_t))=H^0(E_t, \Cal O(-E_t))
$$
But $-E_t ^2$ is $-3$ (cf the proof of Lemma 7) hence this last group is zero.
On the other hand, we have also that $H^2(L_t^c,\Cal O_{L_t^c}(E_t))$ is zero by Serre duality.
Finally, by Riemann-Roch formula \cite{G-H, p.472},

$$
\dim H^0(L_t^c,\Cal O(E_t))=\dfrac{E_t\cdot E_t-E_t\cdot K_t}{2}+\chi(L_t^c)=\dfrac{E_t\cdot E_t+E_t\cdot E_t}{2}+1=4 
$$
and we are done.

By \cite{K-S, Theorem 2.1}, it follows from the claim that there exists a smooth (in $t$) family $\Sigma(t)=(\sigma_1(t),\hdots,\sigma_4(t))$
of holomorphic sections of $-K_t$ which form a basis of $ H^0(L_t^c,\Cal O(E_t))$ for each $t$. Remark that, when $L_t^c$ is of type (i), then $E_t$ is ample by Naka\"\i's criterion; and when it is
of type (ii), it is ample outside the zero section. In both cases, $\Sigma$ defines a CR map from Int $\bar Y$ to $\Bbb P^3$ which is an embedding for fixed $t$ when restricted to
an open and dense subset of $L_t^c$.

Now, fix $t>0$. Recall that $L_t^c$ and $L_{f(t)}^c$ are biholomorphic. We notice that since the biholomorphism between these two manifolds comes from a biholomorphism between
$L_t$ and $L_{f(t)}$ it sends $E_t$ onto $E_{f(t)}$. Hence $\Sigma(L_t^c)$ and $\Sigma(L_{f(t)}^c)$ are holomorphic basis of sections of the same divisor, hence their image in $\Bbb P^3$ are isomorphic and differ from a
global automorphism of $\Bbb P^3$, let us call it $A$. Choose a smooth path of automorphisms $(A_s)$ between $A_t=A$ and $A_{f(t)}=Id$.

We now replace the family $\Sigma (s)$ for $s\in [t, f(t)]$ by $T (s)=A_s\Sigma(s)$. It has the property that it maps an open set of the deformation family Int $\bar Y$ in $\Bbb P^3$ in such a way that
$T(t)$ and $T(f(t))$ are equal. Observe that this implies that the CR-isomorphism of order three is in fact isotopic to a biholomorphism of $T(t)$. Since $T(t)$ is an embedding of an open and dense subset of
$L_t^c$, we conclude that $L_t^c$ admits an automorphism isotopic to a self-map of order three. 


Let us call $\phi_t$ this automorphism. Then $\phi_t^3$ is isotopic to the identity. But there is no non-zero holomorphic vector fields on $L_t^c$ (this is proved in \cite{K, p.225} for $L_t^c$ of type (i), but the same arguments
work for $L_t^c$ of type (ii)). Hence $\phi_t^3$ is the identity and $L_t^c$ admits an automorphism of order three.
$\square$
\enddemo

We may finally state.

\proclaim{Lemma 9} Every surface $L_t^c$ embeds as a non-singular cubic surface $S_t$ of $\Bbb P^3$ and the restriction of the embedding $L_t^c\to 
S_t\subset\Bbb P^3$ embeds $\Bbb E_{\alpha}$ as a hyperplane section 
$H$ of $S_t$. 
\endproclaim 
 
\demo{Proof}
We know that there exists an automorphism $\phi_t$ of order three on 
$L_t^c$ whose fixed point set is an elliptic curve $\Bbb 
E_{\alpha}$. Assume that $L_t^c$ is of type (ii). Then this automorphism sends the five exceptional divisors 
of $L_t^c$ onto five disjoint rational curves with self-intersection 
$-1$. So it projects onto an automorphism $\psi_t$ of order three of 
$\Bbb F_2$. Now $\psi_t$ preserves the unique ruling of $\Bbb F_2$ 
and induces thus an automorphism $\chi_t$ of order three of $\Bbb 
P^1$, the base of the bundle $\Bbb F_2\to\Bbb P^1$. A 
straightforward computation shows that $\chi_t$ is conjugated to the 
identity, the multiplication by $\omega$ or the multiplication by 
$\omega^2$ on $\Bbb P^1=\Bbb C\cup \{\infty\}$. Such a map has two 
fixed points ($0$ and $\infty$) or the whole $\Bbb P^1$ of fixed 
points. 
 
Now the fixed points of $\psi_t$ are exactly the fixed points of its 
restriction to the fibers of $\Bbb F_2\to\Bbb P^1$ which are fixed 
by $\chi_t$. And, as before, the restriction to a fiber is 
conjugated to the identity, the multiplication by $\omega$ or the 
multiplication by $\omega^2$ on $\Bbb P^1=\Bbb C\cup \{\infty\}$. We 
infer from this description that the fixed point set of $\psi_t$ is 
four distinct points, or two points and a disjoint $\Bbb P^1$ or two 
disjoint $\Bbb P^1$ or the whole $\Bbb F_2$. In particular, if 
$\psi_t$ is not the identity, this set is not connected. Blowing-up 
at five points and lifting $\psi_t$ will not give an automorphism 
with an elliptic curve as fixed point set. Contradiction. The 
surface $L_t^c$ is the blow-up of $\Bbb P^2$ at six points. 
 
To conclude that $L^c_t$ is a non-singular cubic surface $S_t$ of 
$\Bbb P^2$, it is enough to prove that the six points are in general 
position. 

Assume the contrary. The automorphism $\phi_t$ of $L^c_t$ 
defines an automorphism $\psi_t$ of order three of $\Bbb P^2$ with 
an elliptic curve $C$ as fixed point set. The six points must belong 
to $C$. Assume first that three of them lie on the same projective 
line $D$. Then $\psi_t$ leaves $D$ globally invariant and fixes 
three distinct points of $D$. This implies that the whole $D$ is 
fixed. Contradiction with the fact that the fixed point set of 
$\psi_t$ is an elliptic curve. Assume now that the six points lie on 
the same conic. Then $\psi_t$ leaves this conic globally invariant 
and fixes six distinct points of it. This implies that the whole 
conic is fixed. Once again, contradiction. The surface $L^c_t$ is a 
non-singular cubic surface $S_t$ of $\Bbb P^2$. 

We would like to thank J. Kollar for the following proof. 

 The curve $\Bbb E_{\alpha}$ is homologous to a 
hyperplane section of $S_t$. By \cite{Ha, V.4.8}, the linear class 
of a divisor of a cubic surface is entirely determined by its 
intersection numbers with the $6$ exceptional curves and with a 
hyperplane section. Therefore $\Bbb E_{\alpha}$ embeds in $S_t$ as a 
divisor which is linearly equivalent to the hyperplane section. 
Since the surface is embedded in $\Bbb P^3$ by its hyperplane class, 
a divisor linearly equivalent to the hyperplane section is a 
hyperplane section. $\square$ 
\enddemo 

As a consequence, $L_t$ is an affine cubic surface. 

\remark{Remark}
The proof that $L_t^c$ is a non-singular cubic surface could be shortened if we knew that the CR-structure of the family Int $\bar Y$ extends smoothly to the singular point of
$\partial Y\simeq P^{-1}(\{0\})\setminus \{0\}$. Since it also extends at infinity by Proposition 5, we would obtain a compactified deformation family whose boundary $L_0^c$ is a singular cubic surface of
$\Bbb P^3$. Taking into account that $E_0\subset L_0^c$ is very ample with space of holomorphic sections of dimension $4$, it would follow directly from the claim proved in Lemma 8 and from a result of Schneider \cite{Sc, p.174}
(which generalizes \cite{K-S, Theorem 2.1} used above) that the deformation family embeds in $\Bbb P ^3$ near the boundary leaf. This would imply that every $L_t^c$ for small $t$ (and thus for all $t$) embeds
in $\Bbb P^3$ through the linear system $E_t$ and satisfies the conclusion of Lemma 9.
\endremark
\medskip

We have now to prove that the complex structure on $L_t$ does not depend on $t$.

Let us say that $S_t$ has equation $P_t=0$ in $\Bbb P ^3$.
The automorphism 
of order three of $L^c_t$ defines an automorphism of order three of 
$S_t$ which is the identity on the hyperplane section at infinity $H$. By 
\cite{G-H, p.178}, such an automorphism is the restriction to $S_t$ 
of an automorphism of $\Bbb P^3$. Straightforward computations show 
that it is projectively conjugated to 
$$ 
\pmatrix & \omega && &0 \cr && \omega && \cr &&& \omega & \cr &0 &&& 
1 
\endpmatrix 
\qquad\text{or}\qquad \pmatrix & \omega^2 && &0 \cr && \omega^2 && 
\cr &&& \omega^2 & \cr &0 &&& 1 
\endpmatrix 
$$ 
on $\C^4$. 
 
Therefore, the cubic surface $S_t$ has the form 
$$ 
P_t=Q_t(x,y,z) +a_tw^3=0 
$$ 
where $Q_t$ is a homogeneous polynomial of degree $3$ and $a_t$ a non-zero constant. Remark that we may assume that $a_t$ is one, making a change of coordinate in $w$ if necessary. It follows that
this equation is entirely determined by the modulus of the elliptic 
curve at infinity. But this curve is $\Bbb E_{\alpha}$ for all $t$. 
Therefore, $Q_t$ is a cubic equation of $\Bbb E_{\alpha}$ and all the non-compact leaves are biholomorphic.

Moreover, the deformation family Int $Y\to (0,\infty)$ is CR-isomorphic to 
$$
(z,t)\in  P_h^{-1}(t)\times (0,\infty)\longmapsto t\in (0,\infty)
$$ 

Consider once again the compactified deformation family Int $\bar Y\to (0,\infty)$. As said before, we have a smooth injection
of $\Bbb E_\alpha\times (0,\infty)$ in Int $\bar Y$ by the uniform compactification Lemma. And it extends smoothly in $0$ by Proposition 5. Take the union of the normal bundles
of this family of elliptic curves. This gives a deformation family of line bundles over $\Bbb E_\alpha$ (in the sense of Section 3) of fixed topological degree. It satisfies the hypotheses of
the Dumping Lemma, hence they are all isomorphic. Of course, we already know that, except in $0$, all these bundles are biholomorphic to the normal bundle of the hyperplane section $H$ in $S$, that is
biholomorphic to the natural bundle associated to $\Bbb E_\alpha$. We know now that it is also the case at $0$. Since the (non-compactified) boundary leaf is
$W(\alpha,x)$, we conclude that $x$ is in fact $\tilde\alpha$, the 
natural $\C^*$-bundle associated to $\Bbb E_{\alpha}$. This finishes 
the third step. 
\medskip 
 
Let us sum up the first three steps. Let $\Cal F$ be a foliation by 
complex surfaces diffeomorphic to $\Cal F_{\Bbb C}$. Then 
\medskip 
\noindent (i) the compact leaves $\Cal S_1$ and $\Cal S_2$ are of 
respective type $(\alpha,\beta,\tilde\alpha)$ and 
$(\alpha,\beta',\tilde\alpha)$. 
 
\noindent (ii) the non-compact leaves of the interior part 
(respectively of the exterior part, of the collar) are all 
biholomorphic. We may thus talk of {\it the} non-compact leaf of the 
interior part (respectively of the exterior part, of the collar). 
 
\noindent (iii) if $\Cal F'$ is another such foliation with $\Cal 
S'_1$ biholomorphic to $\Cal S_1$, and $\Cal S'_2$ biholomorphic to 
$\Cal S_2$, then the non-compact leaf of the interior part of $\Cal 
F'$ (respectively of the exterior part, of the collar) is 
biholomorphic to the non-compact leaf of the interior part of $\Cal 
F$ (respectively of the exterior part, of the collar). 
\medskip 
To finish with the proof of Theorem B, we have to prove that, under 
the assumptions of (iii), there really exists a biholomorphism 
between $\Cal F$ and $\Cal F'$. 
 
\medskip 
\noindent $\underline{\text{Step 4 : construction of the 
biholomorphism}}$. Consider now the exterior part of $\Cal F$, 
without the compact leaf $\Cal S_2$. Let us call it $\Cal E$. It 
admits a product foliated covering $Y$ fibering over $(0,\infty)$. As 
shown before, we may uniformly compactify the total space of this covering 
into a family of deformations of cubic surfaces by adding the same 
elliptic curve to each fiber. All these cubics are biholomorphic and 
the family is parametrized by $(0,\infty)$. It follows from Proposition 2 
that this compactified covering is biholomorphic to a product $S\times \R$, for $S$ the fixed cubic surface described in the previous step. 
Therefore the 
product foliated covering $C$ is biholomorphic to $A\times \R$, 
where $A$ is the affine part of $S$. We obviously have a 
commutative diagram of CR-maps 
$$ 
\CD 
C \aro >>> A \times \R \\ 
\aro VVV \aro VVV \\ 
\Cal E \aro >>> \Cal E 
\endCD 
$$ 
i.e. $\Cal E$ is biholomorphic to a fixed model. This allows us to 
construct a biholomorphism $\phi$ between $\Cal E$ (corresponding to $\Cal 
F$) and $\Cal E'$ (corresponding to $\Cal F'$). 
 
We need an extension Lemma.

\proclaim{Lemma 10}
The CR-isomorphism $\phi$ extends as a biholomorphism between $\Cal S_2$ and $\Cal S'_2$.
\endproclaim

\demo{Proof}
Consider the compactified deformation family $f : \bar Y\to [0,\infty)$ coming from $\Cal F$. The interior leaves are biholomorphic to a fixed cubic surface $S$ say
$$
P_h(x,y,z,w)=Q_h(x,y,z)+w^3=0\text{ in } \Bbb P^3\ .
$$
As a consequence, the interior Int $\bar Y$ can be embedded in $\Bbb P^3$ as the family
$$
Q_h(x,y,z,w)+tw^3=0 \leqno t\in (0,\infty)
$$
The deck transformation group of this compactified product foliated covering is generated by 
$$
[x,y,z,w]\in\Bbb P^3\longmapsto [\omega\cdot x, \omega\cdot y,\omega\cdot z,\lambda^{-1}w]\in\Bbb P^3
$$
It follows from this description that the complex structure on the whole family (that is including $t=0$) is uniquely determined by the complex structure of one leaf. In fact, it is obtained by stretching the complex structure of $S$
through the previous automorphism. The effect of this transformation is to make bigger and bigger a neighborhood of the infinite hyperplane section $\Bbb E_\alpha$ and has at end-point ($t=0$) the normal
bundle of $\Bbb E_\alpha$ in $S$. This description characterizes uniquely the deformation family up to CR-isomorphism. It follows that the previous embedding extends to $t=0$.
Now, the same is true from the family $\bar Y'$ coming from $\Cal F'$, hence the CR-isomorphism $\phi$ extends to the boundary.
$\square$
\enddemo

We claim now that the product foliated covering with boundary $Z_1$ of step 2 is CR-trivial (recall that we already know that its {\it interior} is CR-trivial). Since the boundary fiber $\partial Z_1=W$ is pseudo-convex,
this foliated covering, as a deformation family, is pseudo-trivial by \cite{A-V, Proposition 2}. This means that, given a relatively compact open set $A\in\partial Z_1$,
there exists a relatively compact open set $\Cal A$ in $Z_1$ and a CR-isomorphism between $\Cal A$ and $A\times [0,\epsilon)$ for some $\epsilon>0$.  

Consider $W$ as a $\Bbb C^*$ principal bundle 
over $\Bbb E_\alpha$ and let $L$ be the associated line bundle. Let $D$ be a tubular neighborhood of the zero section of $L$ in $L$. We may take  a holomorphic disk bundle over 
$\Bbb E_\alpha$ as $D$. Let $A$ be an annulus bundle included in $D$ (that is the difference between $D$ and a smaller disk bundle included in $D$). It is a relatively compact open subset of $W$. 
Use the property of pseudo-triviality with this $A$. We obtain a relatively compact open set $\Cal A$ in $Z_1$ and a CR-isomorphism $\chi$ sending $\Cal A$ to
$A\times [0,\epsilon)\subset L\times [0,\epsilon)$. Recall that the leaves of the interior of $Z_1$ may be compactified as copies of $L$ by adding a zero section (this is also true at infinity, but we do not need this here). Call
$Z_1^c$ the corresponding family. Remark that $\Cal A$ intersects each fiber of each leaf of $Z_1^c$ in an annulus. Let $\Cal K$ be the open set of $Z_1^c$ whose intersection with a fiber of a leaf of $Z_1^c$ is the corresponding
disk (or punctured disk if the leaf is the boundary leaf). So $\Cal K$ intersects each interior leaf $L_t$ as a disk bundle $K_t$ over $\Bbb E_\alpha$, and the boundary leaf $W$ as a punctured disk bundle $(K_0)^*$. 

We claim that $\chi$ extends as a CR-isomorphism from $\Cal K$ to $D\times [0,\epsilon)$. The extension is made in two steps. First, fix $t$ and consider $(K_t)^*$ ($K_t$ minus the zero section).
Since $L$ has negative first Chern number, 
this means that the zero section of $L$ admits a strictly pseudo-convex neighborhood. Hence we may assume that the boundary of $K_t$ is strictly pseudo-convex. This allows to extend $\chi$ to $(K_t)^*$ (indeed,
$(K_t)^*$ admits a Stein completion by one point, which is obtained by blowing down the zero section of $K_t$; 
and every CR-function on the strictly pseudo-convex boundary of $K_t$ extends to this Stein completion). Notice that the extension is 
unique and is defined using the Bochner-Martinelli kernel. As a consequence, this extension is smooth in the parameter $t$. Now, this extension must fix each fiber of $(K_t)^*$, since it fixes them in restriction to $A_t$. This
implies that $\chi$ is locally bounded near the zero-section of each $(K_t)$ hence extends to $K_t$ by Riemann's Theorem. So finally, we obtain a CR-injection of $\Cal K$ into $L\times [0,\epsilon)$. But this obviously means
that the complex structure of the boundary leaf $W$ extends to $L$ in the family (that is in such a way that this extension is smooth in the transverse parameter). Using this, and recalling that we may also uniformly compactify
{\it all} the leaves of $Z_1^c$ at infinity by adding an elliptic curve, we CR-embed $Z_1^c$ and thus $Z_1$ in a deformation family of $L^c$, the compact ruled surface associated to $L$. But by Proposition 2, such a family is CR-trivial, so, taking into account that the compactifications are uniform, we conclude that $Z_1$ is a CR-trivial product foliated covering.

 Notice that the same line of arguments proves that $Z_0$ as well as the product foliated covering of step 1 (with its boundary) are CR-trivial.


But now, it follows from Corollary 2 that $\phi$ extends to the product foliated covering $Z_1$ of step 2; and then as a biholomorphism between $\Cal S_1$ and $\Cal S'_1$, also from Corollary 2 applied to the product
foliated covering $Z_0$. And finally, applying once again Corollary 2, but this time on the product foliated covering of the interior part (with boundary), $\phi$ extends to a global biholomorphism
between $\Cal F$ and $\Cal F'$.
 
 
 
 
This concludes the proof of Theorem B. $\square$ 
\enddemo 
 
We pass now to the proof of Theorem C. 
 
\demo{Proof of Theorem C} Due to Theorem B, we only have to realize 
a foliation by complex surfaces diffeomorphic to $\Cal F_{\Bbb C}$ 
for every triple $(J(\alpha),J(\beta),J(\beta'))\in \C\times\C\times\C$. Here 
is the construction. Fix such a triple and fix $\alpha$, $\beta$ and $\beta'$ in $\Bbb H$. The foliation of $\Cal 
F_{\Bbb C}$ restricted to the interior part can easily be adapted: 
take for $P$ a homogeneous polynomial such that the projectivization 
of the set $P^{-1}(\{0\})\setminus \{(0,0,0)\}$ is $\Bbb E_{\alpha}$ 
and take as new definitions of $T$ and $S$ (compare with \cite{M-V, 
p.922--923}) 
$$ 
\eqalign{ T(z,u,t)=&(z,\lambda\omega\cdot u,d(t))\cr S(z,u,t)=&(\exp 
(2i\pi\alpha)\cdot z, (\psi(z))^{-1}\cdot u, t) } 
$$ 
where $\psi$ is the automorphism factor associated to the $\Bbb 
C^*$-bundle 
$$ 
P^{-1}(\{0\})\setminus \{(0,0,0)\}\to\Bbb E_{\alpha} 
$$ 
and where $\lambda\omega=\exp (2i\pi\beta)$.

\remark{Remark}
With this convention, we have $\vert \lambda \vert <1$. Notice that it is different from the convention of \cite{M-V}, where we have
$\lambda\omega=-\exp (2i\pi\beta)$, hence $\vert \lambda\vert >1$.
\endremark
\medskip

 Proceed now 
exactly as in \cite{M-V, p.922-923}. 
\medskip 
 
We modify now the foliation of the exterior part and of the collar. 
\medskip 
Set $\lambda'\omega=\exp (2i\pi\beta')$. Let $\mu$ and $\theta$ be smooth real maps satisfying 
\medskip 
\noindent (i) For $t\leq -1$, we have $\mu(t)\exp 
(i\theta(t))=\lambda$. 
 
\noindent (ii) For $t\geq 0$, we have $\mu(t)\exp 
(i\theta(t))=\lambda'$. 
\medskip 
Notice that $\mu(0)=\vert \lambda'\vert <1$. Let $h : \Bbb R\to\Bbb R$ 
be a diffeomorphism satisfying 
\medskip 
 
\noindent (i) $h(t)=t$ for $t\leq -1$. 
 
\noindent (ii) $h$ coincides with the map $t\mapsto \dfrac 
{t}{3t+1}$ on $[0,\infty)$. 
 
\noindent (iii) The only fixed points of $h$ are $0$ and $(-\infty, 
-1]$. 
\medskip 
 
Such a map is constructed in \cite{M-V, p.925}. Define 
$$ 
\phi \ : \ t\in\Bbb R\longmapsto \left\{\eqalign{ &\mu(0)^{1/t}\exp 
(i\theta(0)/t) \text{ if } t>0 \cr &0 \text{ else.} } \right . 
$$ 
 
Notice that $\phi$ is $C^{\infty}$. Let 
$$ 
g \ : \ (z,t)\in\Bbb C^3\times\Bbb R\longmapsto P(z)-\phi(t) \in\Bbb 
C\ . 
$$ 
 
Set 
$$ 
\Xi=g^{-1}(\{0\})\setminus (\{(0,0,0)\}\times\Bbb R) 
$$ 
and 
$$ 
L_t=\{(z,t)\in\Xi \quad\vert\quad P(z)=\phi(t)\}\ . 
$$ 
 
Finally, define 
$$ 
G\ :\ (z,t)\in\Xi \longmapsto (\mu(t)\exp (i\theta(t))\omega\cdot z, 
h(t))\in \Xi 
$$ 
 
Then $G$ sends the leaf $L_s$ onto the leaf $L_{h(s)}$ and the 
quotient of $\Xi$ by $\langle G\rangle$ gives a foliation by complex 
surfaces of $\Bbb S^5\setminus \Cal N$, as in \cite{M-V}. There are 
two compact leaves and it is straightforward to see that they are 
primary Kodaira surfaces of type $(\alpha,\beta,\tilde\alpha)$ and 
$(\alpha,\beta',\tilde\alpha)$.

 
\proclaim{Lemma 11} Consider the foliation by complex surfaces of 
$\Bbb S^5$ described just above. Then, it is diffeomorphic to $\Cal 
F_{\Bbb C}$. 
\endproclaim 
 
\demo{Proof} Due to the minor changes in the definitions of $T$ and 
$S$, the new foliation of the interior part is clearly diffeomorphic 
to $\Cal F_{\Bbb C}$ restricted to Int$(\Cal N)$. 
 
On the other hand, it follows from the previous proof that: 
 
\medskip 
\noindent (i) The (open) exterior part fibers over the circle with 
same fiber (up to diffeomorphism) and same monodromy (given by 
multiplication by $\omega$ on $P^{-1}(1)$) as $\Cal F_{\Bbb C}$. 
 
\noindent (ii) The collar is diffeomorphic to $W^{diff}\times 
(0,1)$. 
\medskip 
 
We deduce from these two facts that the new foliation is 
homeomorphic to $\Cal F_{\Bbb C}$ from one hand, and from the other 
hand diffeomorphic to it outside the two compact leaves. To finish 
with, it is enough to check that the holonomies of the two compact 
leaves are conjugated. But, in the two cases, they are given by the 
functions $d$ and $h$ (see also the remark at the end of Section 1). 
$\square$ 
\enddemo 
 
This completes the proof of Theorem C. $\square$ 
\enddemo 
 
\head {\bf 8. The moduli space of $\Cal F_{\Bbb C}^{diff}$} 
\endhead 
 
\proclaim{Theorem D} 
 
\noindent (i) The set $\Cal C(\Bbb S^5,\Cal F_{\Bbb C}^{diff})$ 
identified with $\C^3$ is a coarse moduli space for $(\Bbb S^5,\Cal 
F_{\Bbb C}^{diff})$. 
 
\noindent (ii) There does not exist a fine moduli space for $(\Bbb 
S^5,\Cal F_{\Bbb C}^{diff})$. 
\endproclaim 
 
\demo{Proof} The proof is very similar to that given in Section 5 
for the Reeb foliation. Take any deformation family $\pi : M\to X$ 
of $(\Bbb S^5,\Cal F_{\Bbb C}^{diff})$ and consider the natural map 
$\alpha_{\pi} : X\to \Bbb C\times\C \times \C$. The three components 
of this map may be thought of as the modular function on the common 
base of the compact leaves, on a fixed fiber of the first compact 
leaf and on a fixed fiber of the second compact leaf. Indeed, such a 
deformation family $\Cal M_1\to X$ (respectively $\Cal M_2\to X$) 
induces a deformation family $\Cal D_1$ (respectively $\Cal D_2$) of the first compact leaf $\Cal 
S_1^{diff}$ (respectively the second compact leaf $\Cal S_2^{diff}$) 
by Proposition 6. On the other hand, $\Cal D_1\to X$ (respectively 
$\Cal D_2\to X$) is a locally trivial smooth fiber bundle with 
structural group Diff$^+(\Cal S_1^{diff})$ (respectively 
Diff$^+(\Cal S_2^{diff})$ - the $+$ meaning orientation-preserving). 
Now, by Lemma 2, an orientation-preserving diffeomorphism of $\Cal 
S_1^{diff}$ (respectively $\Cal S_2^{diff}$) is in fact a bundle 
isomorphism. Therefore, $\Cal D_1\to X$ (respectively $\Cal D_2\to 
X$) induces a deformation family of the base of $\Cal S_1^{diff}$ 
(respectively $\Cal S_2^{diff}$). This allows us to identify 
$\alpha_{\pi}$ to a triple of modular functions and to conclude that 
$\C^3$ is a coarse moduli space for $(\Bbb S^5,\Cal F_{\Bbb 
C}^{diff})$. 
 
Assume now that there exists a fine moduli space. Then every 
deformation family of $(\Bbb S^5,\Cal F_{\Bbb C}^{diff})$ with all 
fibers biholomorphic to $\Cal F_{\Bbb C}$ is biholomorphic to a 
product. By restriction to the base of a fixed compact leaf, such a 
deformation family yields a family of complex structures on the 
smooth torus with all leaves biholomorphic to $\Bbb E_{\omega}$. Of 
course, this family has to be holomorphically trivial. Now, there 
exist non holomorphically trivial families of complex structures 
with all leaves biholomorphic to $\Bbb E_{\omega}$. This is due to 
the existence of a non-trivial automorphism $\alpha$ of order three 
on $\Bbb E_{\omega}$. For example, one may glue two copies of $\Bbb 
C\times \Bbb E_{\omega}$ along $\Bbb C^*  \times \Bbb E_{\omega}$ 
with 
$$ 
h \ : \ (z,w)\in \Bbb C^*  \times \Bbb E_{\omega}\longmapsto (1/z, 
\alpha(w))\in \Bbb C^*  \times \Bbb E_{\omega} 
$$ 
to obtain such an example. Assume now that $\alpha$ extends to a 
biholomorphism $A$ of $\Cal F_{\Bbb C}$. Then, in the same way, glue 
two copies of $\Bbb C\times \Bbb S^5$ endowed with the foliation 
$\Cal F_{\Bbb C}$ on each $\Bbb S^5$ along $\Bbb C^*\times \Bbb S^5$ 
by 
$$ 
H \ : \ (z,w)\in\Bbb C^*\times \Bbb S^5 \longmapsto 
(1/z,A(w))\in\Bbb C^*\times \Bbb S^5\ . 
$$ 
This gives a deformation family of $(\Bbb S^5,\Cal F_{\Bbb 
C}^{diff})$ with all fibers biholomorphic to $\Cal F_{\Bbb C}$ and 
such that its restriction to the base of a fixed compact leaf is 
non-trivial. Contradiction. 
 
To finish with, it is thus enough to prove the following Lemma. 
 
\proclaim{Lemma 12} The foliation $\Cal F_{\Bbb C}$ admits a 
biholomorphism $A$ of order three which extends an automorphism of 
order three $\alpha$ on the the base of the compact leaves. 
\endproclaim 
 
\demo{Proof} We just have to define this automorphism on the product 
foliated coverings used to foliate $\Cal N$ and $\Bbb 
S^5\setminus\Cal N$. We use the notations of \cite{M-V}. Start with 
the automorphism 
$$ 
\alpha \ : \ [z_1,z_2,z_3]\in\Bbb E_{\omega}=\{[z]\in\Bbb P^2 
\quad\vert\quad z_1^3+z_2^3+z_3^3=0\} \longmapsto [z_1,\omega\cdot 
z_2,\omega^2\cdot z_3]\in\Bbb E_{\omega} 
$$ 
 
Consider first the product foliated covering $\tilde X/\langle 
S\rangle$. By Proposition 2, it admits a compactification 
biholomorphic to $L^c\times [0,\infty)$, where $L^c$ is the 
compactification of $W$ as a ruled surface. Therefore, $\tilde 
X/\langle S\rangle$ is biholomorphic to $L\times [0,\infty ) 
\setminus \{s_0\times 0\}$, where $s_0$ denotes the zero section of 
the line bundle $L$. By definition of $L$, the automorphism $\alpha$ 
extends to an automorphism of order three of $L$ and of 
$W=L\setminus s_0$. Therefore, we may define a biholomorphismm of 
order three on $\tilde X/\langle S\rangle$. It is easy to see that 
it commutes with $T$ and thus descends to a biholomorphism of order 
three of $(\Cal N,\Cal F_{\Bbb C})$. For the exterior, just take the 
biholomorphism 
$$ 
(z_1,z_2,z_3,t)\in\C^3\times\R \longmapsto (z_1,\omega\cdot 
z_2,\omega^2 \cdot z_3,t)\in\C^3\times\R 
$$ 
and verify that it preserves $\Xi$ and commutes with $G$. It is now 
straightforward to check that this biholomorphism glue with the 
previous one and define $A$. This completes the proof of Theorem D. 
$\square$ 
\enddemo 
\enddemo 
 
\head {\bf 9. Concluding remarks }
\endhead 
 
In this paper, we do not speak about the foliated analogue of 
Kodaira-Spencer theory of small deformations of compact complex 
manifolds (see \cite{K-S}, \cite{M-K}). As noticed in Section 5, the 
notion of fine moduli space is not a good concept in the theory of 
compact complex manifolds. The right one is that of versal 
deformation space, that is of {\it local} moduli space. One may thus 
expect that the same is true for foliations by complex manifolds.

Let $\Cal F$ be a foliation by complex manifolds on $X$. Define 
$\Theta_{\Cal F}$ as the sheaf of germs of smooth vector fields on 
$X$ which are tangent to $\Cal F$ and holomorphic when restricted to 
any leaf of $\Cal F$. As in the classical case, the cohomology group 
$H^1(X,\Theta_{\Cal F})$ should describe the infinitesimal 
deformations of $(X,\Cal F)$. Of course, it can be 
infinite-dimensional as it is the case for the second example of 
Section 5. Nevertheless, in the case it is finite-dimensional, it is 
quite possible that the classical Theorems could be adapted with 
essentially the same proof. For example, the following result should 
be true: if $H^2(X,\Theta_{\Cal F})=0$, then there exists a germ of 
deformation family of $(X,\Cal F^{diff})$ of dimension 
$H^1(X,\Theta_{\Cal F})$ centered in $(X,\Cal F)$ which contains all 
``sufficiently small'' deformations of $(X,\Cal F)$. It should be 
interesting to develop such a theory and to apply it to deformations 
of CR-fiber bundles and to deformations of foliations given by 
suspension. 
 
\vfill \eject  
 
\enddocument